\documentclass[
10pt, 
letterpaper, 
oneside, 
headinclude,footinclude, 
BCOR5mm, 
]{scrartcl}

\usepackage{longtable} 
\usepackage[mmddyyyy]{datetime}

\usepackage{authblk}
\usepackage{amsthm}
\usepackage{fouriernc,amsmath}
\usepackage[T1]{fontenc}
\usepackage{indentfirst}
\usepackage{placeins}

\DeclareUnicodeCharacter{2061}{}

\usepackage{tocloft} 
 
\usepackage[
nochapters, 
beramono, 
eulermath,
pdfspacing, 
dottedtoc 
]
{classicthesis} 

\usepackage{arsclassica} 

\usepackage[T1]{fontenc} 

\usepackage[utf8]{inputenc} 

\usepackage{graphicx} 
\graphicspath{{Figures/}} 

\usepackage{enumitem} 

\usepackage{lipsum} 

\usepackage{subfig} 

\usepackage{amsmath,amssymb,amsthm} 

\usepackage{varioref} 

\theoremstyle{definition} 

\theoremstyle{plain} 

\theoremstyle{remark} 

\hypersetup{
colorlinks=true, breaklinks=true, bookmarks=true,bookmarksnumbered,
urlcolor=webbrown, linkcolor=RoyalBlue, citecolor=webgreen, 
pdftitle={}, 
pdfauthor={\textcopyright}, 
pdfsubject={}, 
pdfkeywords={}, 
pdfcreator={pdfLaTeX}, 
pdfproducer={LaTeX with hyperref and ClassicThesis} 
} 

\title{\normalfont\spacedallcaps{First Digit Probability and Benford's Law}} 

\author{\spacedlowsmallcaps{Irina Pashchenko}} 

\affil{\small \href{mailto:pashchenko.irina.16@gmail.com}{pashchenko.irina.16@gmail.com}}

\date{\small \today} 

\begin{document}

\renewcommand{\sectionmark}[1]{\markright{\spacedlowsmallcaps{#1}}} 
\lehead{\mbox{\llap{\small\thepage\kern1em\color{halfgray} \vline}\color{halfgray}\hspace{0.5em}\rightmark\hfil}} 

\pagestyle{scrheadings} 


\maketitle 

\setcounter{tocdepth}{3} 

\tableofcontents 

\listoffigures 

\listoftables 


\newpage
\section*{Abstract} 

The following work is written in easy language for college level students. It shows how the first digit probabilities of a group of continuous real-valued functions can be calculated. Thus, examples explaining how the probabilities are related to specific real-life situations, as well as the summary for all basic algebraic functions, were brought to the reader's attention. Besides, the Benford's formula was derived with no use of any additional guiding sources.

Moreover, a comprehensive analysis of a group of certain discrete functions was performed by approximating the functions to the above-mentioned continuous ones, taking limits, and other methods.

The work can be applied for calculating the first digit probabilities of more advanced functions as well while using the same approach. Furthermore, the technique can be useful while dealing with a large set of highly approximate numbers and a conclusion about their nature needs to be made.


\newpage
\section{Acknowledgements}

I would like to thank my son, Vitaliy Goncharenko, for his emotional support.

Moreover, I am grateful to my uncle, Mark Ayzenberg, for his professional guidance.

Finally, I would like to extend my sincere gratitude to the University of West Georgia Math department faculty members who accepted my invitation for a presentation where I had a chance to talk about my research work that was done for the paper.

Nevertheless, all the glory belongs to the Lord.


\newpage
\section{Preface}

Simon Newcomb (1835 – 1909) was a well-known American astronomer. \cite{r_Newcomb} He was the first one who had discovered the first digit phenomenon, called Benford’s law later, and had published a short paper about it in 1881. Newcomb noticed that in his logarithmic book, that showed logarithms of numbers beginning with the lowest digits, the first pages were more worn out than the last ones. Thus, he made a conclusion that first pages were used more often than the last ones; consequently, first digits do not appear with the same frequency. 

Newcomb managed writing a table of first two digits probabilities. \cite{r_Benford_1} However, his article was written with no mathematical background, that is why it did not get much attention at that time.

Frank Benford (1883 – 1948) was an American electrical  engineer and physicist. \cite{r_Benford} He discovered the same interesting detail about his logarithmic book pages in 1938. After that, Benford collected a large data set containing about $20,000$ numbers. The numbers were taken from different sources like algebraic functions values, addresses of people, scientific constants, areas of rivers, etc. Even though the groups of those numbers were not related to each other, Benford computed the first digit frequencies for each group and the average of all the groups for every digit. Then, he made a conclusion that the frequency $F_a$ of the digit $a$ as the first digit can be found by the following formula. \cite{r_Benford_2}

\begin{equation}
F_a = \log_{10}(\frac {a+1}{a})
\label{B_formula}
\end{equation}

Formula~\vref{B_formula} was named as Benford’s formula later. Nevertheless, Benford did not find any reasonable explanation of the phenomenon.

Many various mathematicians and other scientists have performed different kinds of research on Benford’s law since Newcomb and Benford had published their works. One of the major applications of the research is revealing of income tax evasion using a low of digital frequencies. \cite{r_Benford_3}

The paper presented below has a purpose of showing that Benford’s formula works for exponential functions regardless of their bases and for those functions whose limits approach to the exponential ones. All other functions do not follow Benford’s law.

In addition, this paper can shed some light on the data table used by Benford or other tables. If a group of numbers in a table represent an object growth rate or any object property’s change in time, and the data first digits follow a certain function’s first digit probabilities rule, there is a chance that the above-mentioned object grows or its property changes according to this specific function.


\newpage
\section{Introduction}

Since our goal is to see how a group of functions can describe real-life situations, all the functions will be graphed and analyzed in the first quadrant only since this quadrant provides us with positive $x$- and $y$-values.

Also, let us acknowledge that if we are talking about the first digit of a number, we mean the first digit of the whole portion of a number in a case if the number is a decimal. In other words, we assume that the number is bigger or equal to the number $1$. That is why the range of all the functions that will be discussed later should be restricted to $[1, \infty)$. 

The domain of each function will be restricted to match the corresponding restricted range. For instance, the domain of the exponential function will be restricted to $[0,\infty)$ to match the range mentioned above.

If a random real number is selected in the $[1, \infty)$ interval, the chance of getting a number which starts with a particular digit will be equal to the chance of getting a number which starts with any other digit. However, the vast majority of the functions that will be discussed in this work will have a different outcome. We will see how the probability of each first digit depends on the function’s shape.

In order to be able to calculate the first digit probability of each function’s y-value, a particular ruler is proposed. The ruler takes the interval $[1, \infty)$ and has infinitely many sections in it. Each section is exactly ten times longer than the previous one.

Moreover, each section is divided into nine equal parts. Each part is colored with a particular color. Connecting all the parts and sections together brings us back to the interval  $[1, \infty)$ without any overlapping or gaps. The same color is used to label all the intervals containing numbers that start with the same digit. Therefore, nine different colors were used for the infinitely long ruler. 

\begin{figure}[h]
\centering 
\includegraphics[width=0.5\columnwidth]{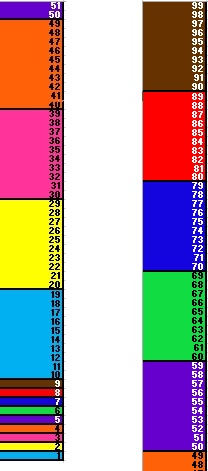} 
\caption{The measuring color-coded ruler for the range $[1, 100)$} 
\label{ruler} 
\end{figure}

The limited part of the ruler shown in Figure~\vref{ruler}  can be used to calculate the first digit probability in an interval $[1, 100)$. Each color will represent a total group of $y$-values that start with a particular digit. The complete ruler can be applied for any real-valued function on a specified range  $[1, \infty)$.

\begin{figure}[!h]
\centering 
\includegraphics[width=1\columnwidth]{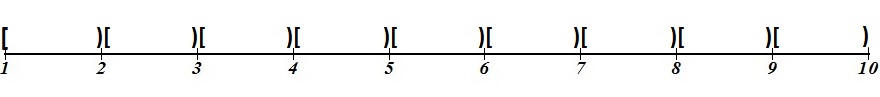} 
\caption{The measuring ruler for the range $[1, 10)$} 
\label{smaller_ruler} 
\end{figure}

The first section of the ruler is shown in Figure~\vref{smaller_ruler}. It starts at $1$ and ends at $10$. In order to calculate the first digits probability fairly, the number $1$ is included to the interval, but the number $10$ is not. Therefore, the interval is $[1, 10)$. A randomly chosen number from this interval will have the same chance of having any of the nine digits as the leading ones because the interval is broken into nine pieces of equal length. Each of them starts at a closed point representing a whole number from $1$ to $9$ and ends at an open point representing the next whole number from $2$ to $10$. This fact makes all the chances of picking a number, which start with a particular digit equal to each other. 

\begin{figure}[!h]
\centering 
\includegraphics[width=1\columnwidth]{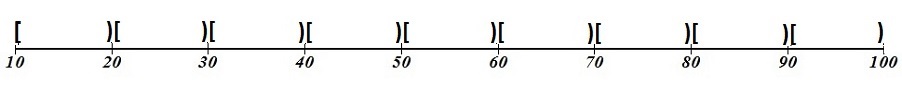} 
\caption{The measuring ruler for the range $[10, 100)$}   
\label{bigger_ruler}
\end{figure}

The same partition will be used for the rest of the sections of the ruler. The interval $[10, 100)$ is shown in Figure~\vref{bigger_ruler}. It is also broken into nine equal intervals in the same manner. As it was mentioned above, the section $[10, 100)$ is exactly ten times longer than the previous one.

\begin{figure}[!h]
\centering 
\includegraphics[width=1\columnwidth]{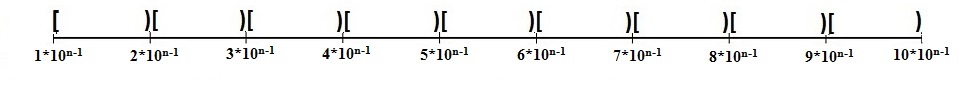} 
\caption{The measuring ruler for the range $[10^{n-1}, 10^n)$}   
\label{general_ruler}
\end{figure}

Any of the sections of the ruler can be expressed using a general form. It is $[10^{n-1}, 10^n)$, where $n$ is a natural number. In particular, all the parts of the section can be listed as the following intervals: 
\newline
$[1 \cdot 10^{n-1}, 2 \cdot 10^{n-1}), [2 \cdot 10^{n-1}, 3 \cdot 10^{n-1}), [3 \cdot 10^{n-1}, 4 \cdot 10^{n-1}), . . . , [9 \cdot 10^{n-1}, 10^n)$. Again, the chances of picking a number with a particular first digit from the interval $[10^{n-1}, 10^n)$ are equal to each other due to the equal lengths of the subintervals. The interval is shown in Figure~\vref{general_ruler}.

In order to compare each function's first digit probability to each other, a term of the first digit probability will be used. The abbreviation $fdp$ for discussion purposes and $P_k$ for formulas where $k$ is the first digit of a number will be utilized in this work. 
Let us recall that when any event takes place (rolling a die, etc.), the probability of each single outcome of the event is a fraction from $0$ to $1$ which may be equal to either $0$ or $1$ as well and the sum of the probabilities of all the outcomes equals to $1$. Since any number of any real-valued function’s range starts with a digit from $1$ to $9$, the probability of any of the digits to be the first one is located in an interval $0 \leq P_k \leq 1$ where $k = 1, 2, …, 9$ and $P_1 + P_2 + P_3 + P_4 + P_5 + P_6 + P_7 + P_8 + P_9 = 1$.

\FloatBarrier
In this work, a brief analysis on how all the $P_k$ relate to each other will be done for each function that will be discussed. A particular method will be used to derive a formula of the $fdp$ of each function. As it was mentioned above, a specific interval of the $y$-axis will be used as a range for each function, which is $[1, \infty)$. In addition, a certain partition of the interval into sections was proposed, which breaks it into subintervals $[10_{n-1}, 10_n)$ or $[1 \cdot 10_{n-1}, 10 \cdot 10_{n-1})$, where $n$ is a natural number. 

\begin{figure}[h]
\centering 
\includegraphics[width=1\columnwidth]{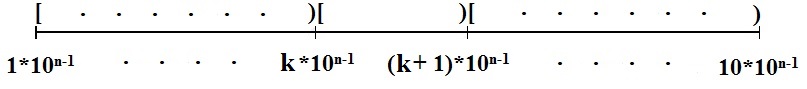} 
\caption{The specific measuring ruler for the range $[10^{n-1}, 10^n)$}   
\label{specific_ruler}
\end{figure}

Each interval $[10_{n-1}, 10 \cdot 10_{n-1})$ will be partitioned into nine subintervals of equal length, each of which can be expressed as $[k \cdot 10_{n-1}, (k + 1) \cdot 10_{n-1})$, as shown in Figure~\vref{specific_ruler}. Picking any number from the mentioned above subinterval will automatically result on having a digit $k$ as the leading digit of the number. 

\FloatBarrier
Let us use a general real-valued continuous function $y = f(x)$ on the range $[1, \infty)$ for the purpose of finding its $fdp$. We will derive a formula for $P_k$, which will show the total probability of getting any $y$-value that starts with a digit $k$. In other words, we need to find a formula for calculating the ratio of the probability of a particular $y$-value being located in the interval $[k \cdot 10^{n-1}, (k + 1) \cdot 10^{n-1})$ to the probability of the same value being located in the interval $[1 \cdot 10^{n-1}, 10 \cdot 10^{n-1})$ or $[10^{n-1}, 10^n)$ which is nine times longer. Let us first find the formulas for $P_k$ of each function on the above interval and then we can analyze all the results. 

\begin{figure}[!h]
\centering 
\includegraphics[width=1\columnwidth]{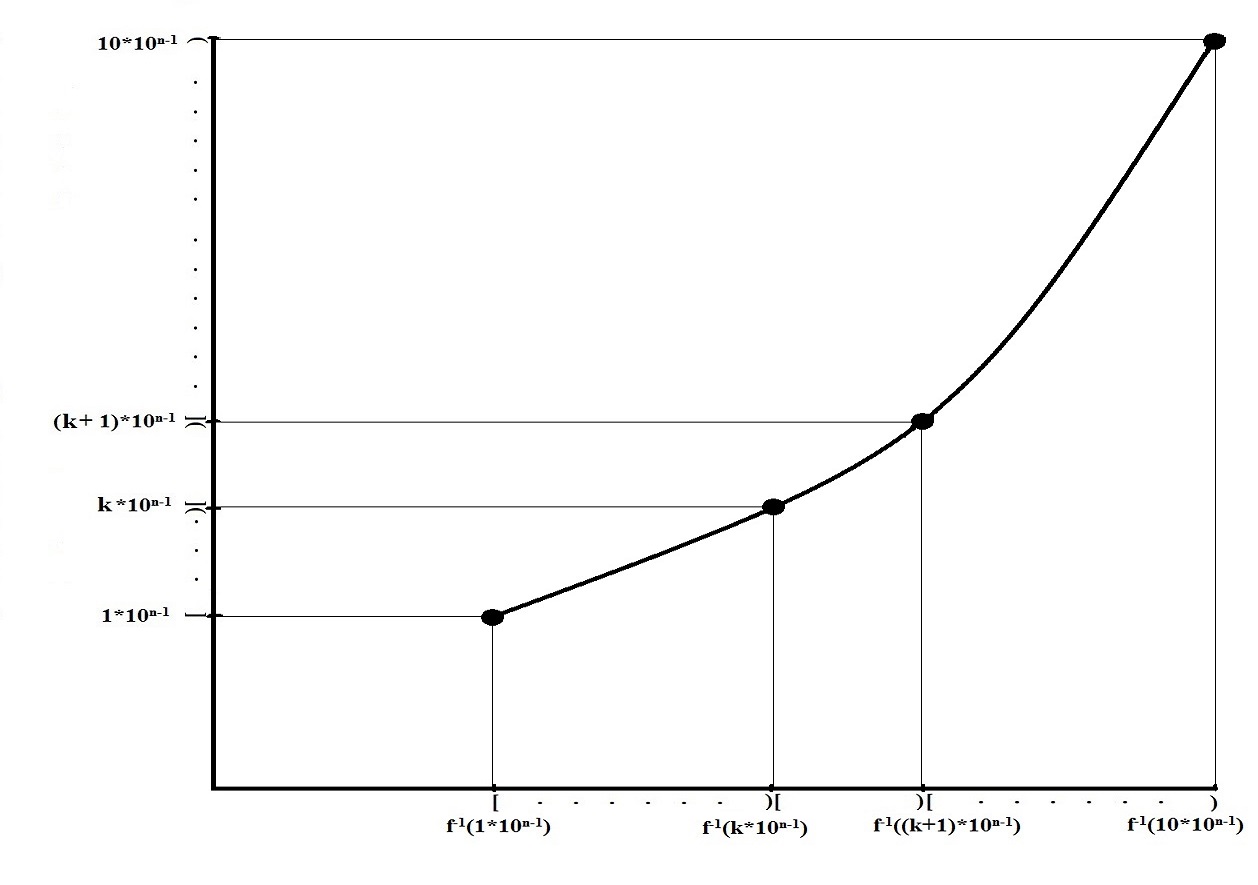} 
\caption{Continuous function $y = f(x)$}   
\label{XY_axis_ruler}
\end{figure}

\FloatBarrier
The part of the continuous function $y = f(x)$ that is shown in Figure~\vref{XY_axis_ruler} has the range $[1 \cdot 10^{n-1}, 10 \cdot 10^{n-1})$. Therefore, the domain of this part is $[f^{-1}(1 \cdot 10^{n-1}), f^{-1}(10 \cdot 10^{n-1}))$. Moreover, all the $x$-values of the interval $[f^{-1}(1 \cdot 10^{n-1}), f^{-1}(10 \cdot 10^{n-1}))$ are mapped into the $y$-values of the interval $[k \cdot 10^{n-1}, (k + 1) \cdot 10^{n-1})$.

Also, all nine subintervals of the $y$-axis interval $[1 \cdot 10^{n-1}, 10 \cdot 10^{n-1})$ have equal length. However, nine corresponding subintervals of the $x$-axis interval $[f^{-1}(1 \cdot 10^{n-1}), f^{-1}(10 \cdot 10^{n-1}))$ are expected to have equal length only in a case if $y = f(x)$ is a linear function. In all other cases, unless a coincidence takes place, all nine subintervals of this interval will have different length. 

The percentage of the $y$-values mapped into $[k \cdot 10^{n-1}, (k + 1) \cdot 10^{n-1})$ out of those that mapped into $[1 \cdot 10^{n-1}, 10 \cdot 10^{n-1})$ depends directly on the length of the interval $[f^{-1}(k \cdot 10^{n-1}), f^{-1}((k+1) \cdot 10^{n-1}))$. The longer this interval is, the more numbers out of it will be mapped into the corresponding interval $[k \cdot 10^{n-1}, (k + 1) \cdot 10^{n-1})$. 

The ratio of the length of the interval $[f^{-1}(k \cdot 10^{n-1}), f^{-1}((k+1) \cdot 10^{n-1}))$ to the interval $[f^{-1}(1 \cdot 10^{n-1}), f^{-1}(10 \cdot 10^{n-1}))$ equals to the ratio of a number of the $y$-values mapped into $[k \cdot 10^{n-1}, (k + 1) \cdot 10^{n-1})$ to a number of those that mapped into $[1 \cdot 10^{n-1}, 10 \cdot 10^{n-1})$. 

In other words,
\begin{equation}
P_k = \frac{f^{-1}((k+1)\cdot 10^{n-1}) - f^{-1}(k \cdot 10^{n-1})}{f^{-1}(10 \cdot 10^{n-1}) - f^{-1}(1 \cdot 10^{n-1})}
\label{l_formula}
\end{equation}
or
\begin{equation}
P_k = \frac{f^{-1}((k+1)\cdot 10^{n-1}) - f^{-1}(k \cdot 10^{n-1})}{f^{-1}(10^n) - f^{-1}(10^{n-1})}
\label{s_formula}
\end{equation}

\noindent
where $k$ is a digit from $1$ to $9$ and $n$ is a natural number. 
Since Formula~\vref{s_formula} is shorter, it will be used in this work to calculate $fdp$ of every continuous real-valued function that will be discussed later. 
\newline


\newpage
\section{First digits probabilities of basic continuous functions}

\subsection{$P_k$ values of the exponential function $y = a^x$}
\subsubsection{Deriving the $P_k$ formula}

\begin{figure}[h]
\centering 
\includegraphics[width=1\columnwidth]{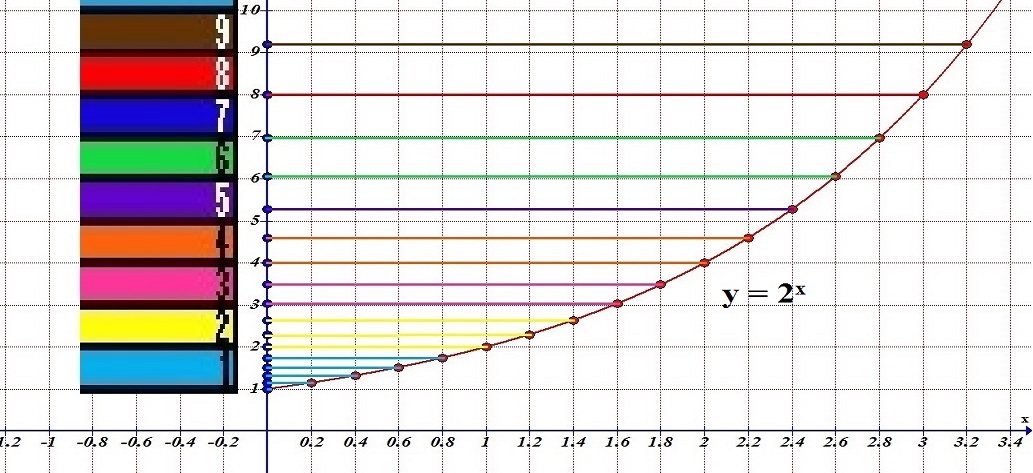} 
\caption{Exponential function $y = 2^x$ on the range $[1, 10)$}   
\label{Base_2_Smaller}
\end{figure}

Let us first recall what the exponential function $y = 2^x$ looks like. The first digit probabilities of the function, that is in the range $[1, 10)$, are explained in Figure~\vref{Base_2_Smaller}. The ruler that was proposed earlier is used to measure the probabilities.

\FloatBarrier
As we can see clearly from the graph, all nine digits have different probabilities in this range. After all the horizontal lines connecting the graph with its corresponding $y$-values on the particular range were included to the picture, the continuous exponential function might look like a discrete one with $x$-axis steps of $0.2$. Therefore, the conclusion that smaller digits have higher probabilities can be proposed, but not finalized. 

\begin{figure}[!h]
\centering 
\includegraphics[width=1\columnwidth]{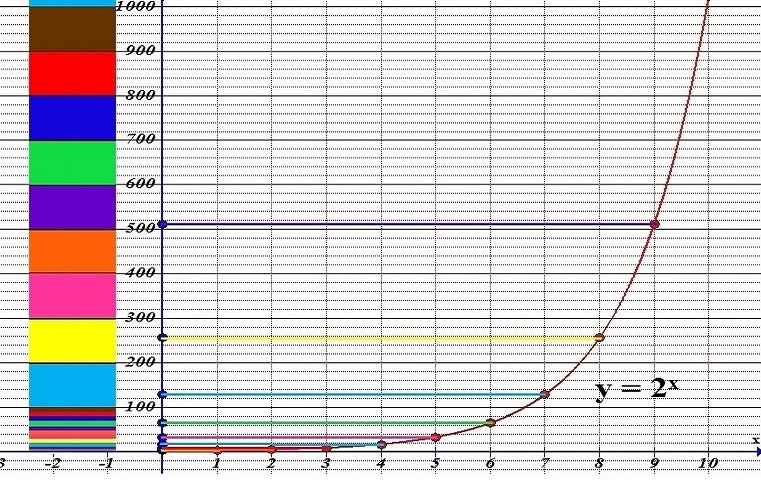} 
\caption{Exponential function $y = 2^x$ on the range $[1, 1000)$}   
\label{Base_2_Larger}
\end{figure}

\FloatBarrier
For instance, it might look like the digit $6$ has a higher probability than the digit $5$. However, it is not true. The problem is that steps of $0.2$ on that particular range are too big to show that in reality the smaller the digit is, the higher its probability. That is why the formula for $P_k$ is needed.

Before we derive the formula, let us look at other graphs of the exponential function. Figure~\vref{Base_2_Larger} shows the function $y = 2^x$ on the range $[1, 1000)$. As we can see, all the digits have different probabilities on that range as well as on the previous one. However, the $x$-axis steps are not small enough to make any specific conclusion about all $P_k$ values on that range.

\begin{figure}[h]
\centering 
\includegraphics[width=1\columnwidth]{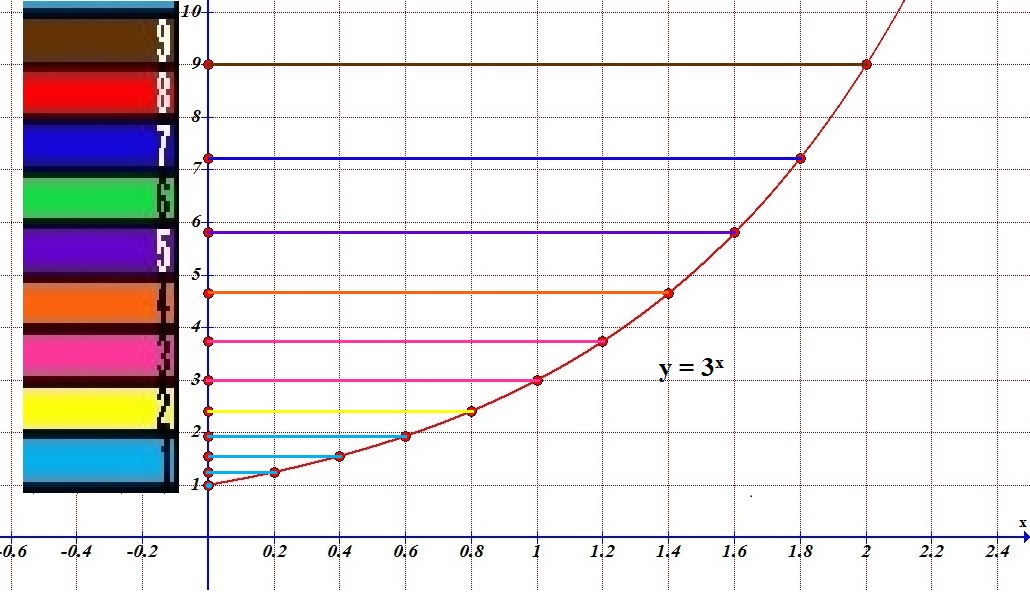} 
\caption{Exponential function $y = 3^x$ on the range $[1, 10)$}   
\label{Base_3_Smaller}
\end{figure}

\FloatBarrier
Let us look at the exponential function $y = 3^x$ now. Figure~\vref{Base_3_Smaller} shows the function on the range $[1, 10)$. Again, all the digits have different probabilities in that interval.

\begin{figure}[h]
\centering 
\includegraphics[width=1\columnwidth]{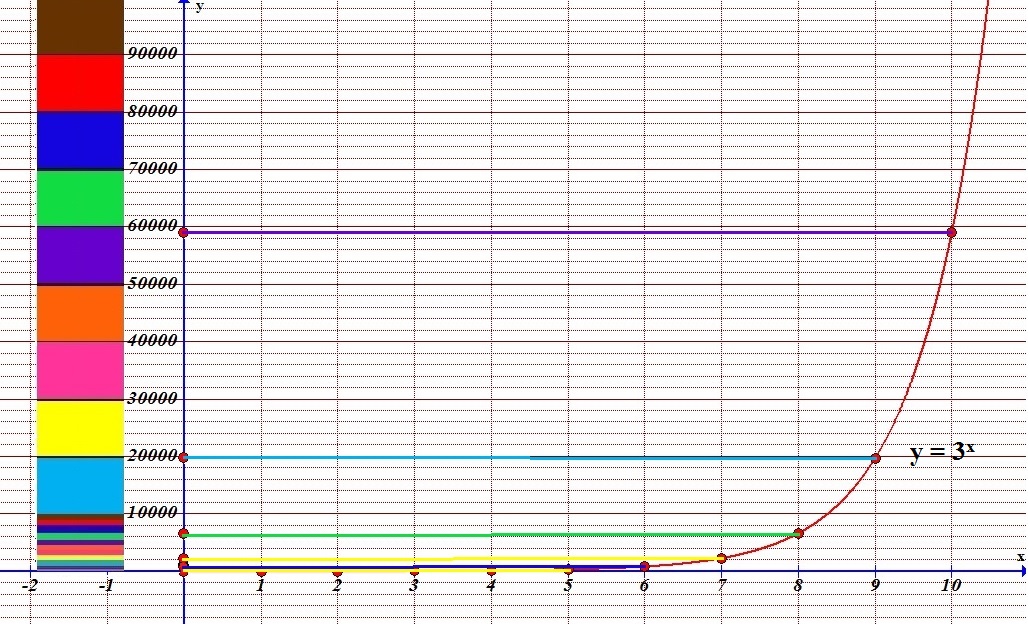} 
\caption{Exponential function $y = 3^x$ on the range $[1, 100000)$}   
\label{Base_3_Larger}
\end{figure}

\FloatBarrier
The same function is shown on the range $[1, 1000)$ in Figure~\vref{Base_3_Larger}. The new range will bring us to the same conclusion: different digits have different probabilities in this interval as well.

We had a chance to see the exponential functions with different bases on different ranges. Let us derive the $P_k$ formula now.

As it was mentioned above, Formula~\vref{s_formula} will be used for each particular function. We will derive the $P_k$ formula for the exponential function $y = a^x$ where $0 \leq x$ < $\infty$ and $a$ > $1$.
First, $f(x)=a^x$. Then, $f^{-1}(x)=\log_a⁡(x)$ according to the definition of the logarithmic function. Next, 

$$P_k = \frac{\log_a⁡((k+1) \cdot 10^{n-1}) - \log_a⁡(k \cdot 10^{n-1})}{\log_a(⁡10^n) - \log_a⁡(10^{n-1})} = 
\frac{\log_a⁡(\frac{(k+1) \cdot 10^{n-1}}{⁡k \cdot 10^{n-1}})}{\log_a (\frac {⁡10^n}{⁡10^{n-1}})}$$ 
$$= \frac{\log_a (⁡\frac {k+1}{⁡k})}{\log_a (⁡10)} = \frac{\frac {\log_{10} (⁡\frac {k+1}{⁡k})}{\log_{10}(a)}}{\frac {\log_{10} (⁡10)}{\log_{10} (a)}} = \frac{\log_{10} ⁡(\frac {k+1}{⁡k})}{\log_{10} (⁡10)} = \frac{\log_{10} ⁡(\frac {k+1}{⁡k})}{1} = \log_{10} ⁡(\frac {k+1}{⁡k});$$

$$ $$
Thus,
\begin{equation}
P_k = \log_{10} ⁡(\frac {k+1}{⁡k})
\label{e_formula}
\end{equation}
$$ $$

\begin{figure}[h]
\centering 
\includegraphics[width=1\columnwidth]{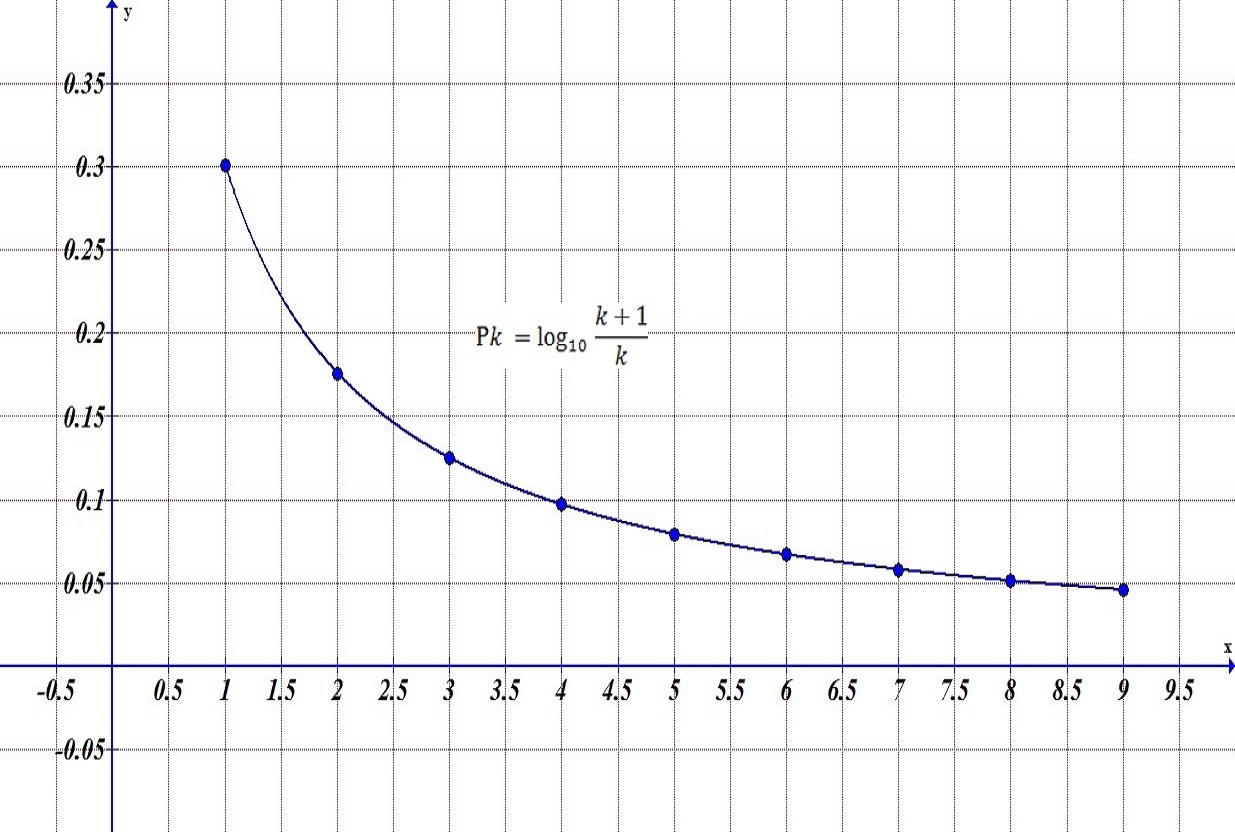} 
\caption{The function $P_k = \log_{10} ⁡(\frac {k+1}{⁡k})$}   
\label{Benfords_Probability_2}
\end{figure}

\begin{figure}[h]
\centering 
\includegraphics[width=1\columnwidth]{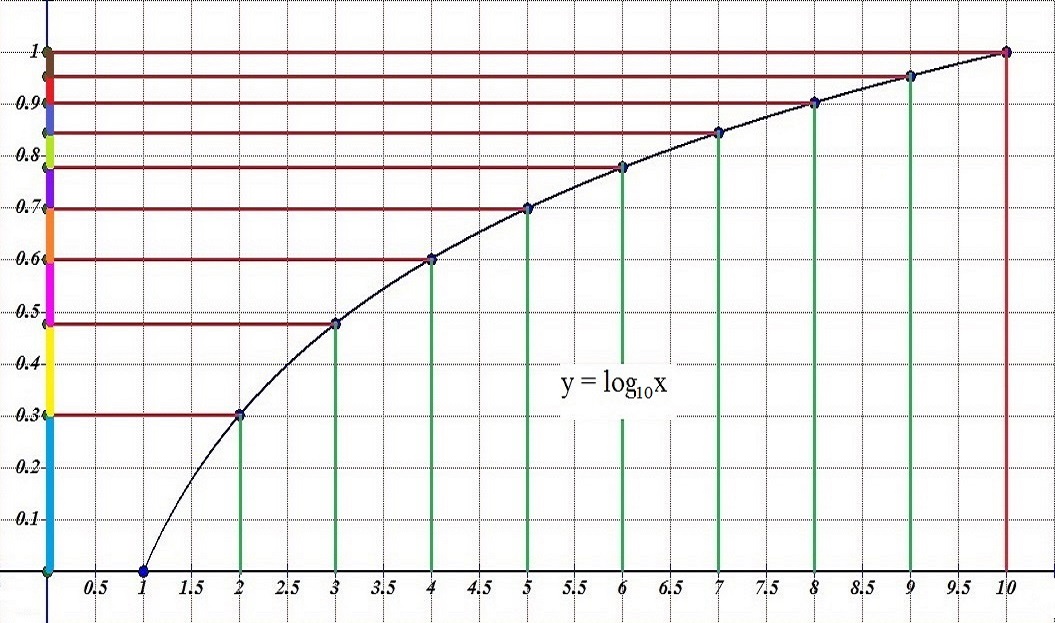} 
\caption{The function $y = \log_{10}(⁡x)$}   
\label{Benfords_Probability_1}
\end{figure}

Formula~\vref{e_formula} satisfies Benford's Formula~\vref{B_formula} \cite{r_Benford_2}. It tells us two very important facts. 

First, the formula does not contain the exponent $n$ as well as any other references to the range. It means that all nine $P_k$ values of the exponential function will be the same on any range. 

Second, the formula does not have any reference to the initial general base $a$ of the function. It means that all $P_k$ values do not depend on the base of the exponential function. In other words, all exponential functions with the base $a$ which is greater than $1$, for instance $y = 2^x, y = e^x, y = 3^x$, etc., have the same $fdp$ values.

Due to that fact that $k$ is a whole number from $1$ to $9$, the newly derived function should be considered as a discrete one on the domain $[1, 9]$. The graph of it is shown in Figure~\vref{Benfords_Probability_2}.

Now, it is clear that the smaller first digit is, the higher its probability. Moreover, all the probabilities are positive numbers less than $1$. This fact satisfies the probability rules. 

The only thing that we still need to check is whether the sum of all nine probabilities will total to $1$. A calculator could definitely be used for this purpose; however, there is a couple of more neat solutions to the problem.

First, let us modify the formula.
\begin{equation}
P_k = \log_{10} (⁡\frac {k+1}{⁡k}) = \log_{10}⁡(k+1) - \log_{10}(⁡k) 
\label{e_m_formula}
\end{equation}

Then, we will calculate the sum and prove that it equals to $1$: 

\begin{proof} 
$\Sigma_{k=1}^9 P_k = \Sigma_{k=1}^9 (\log_{10}⁡(k+1) - \log_{10}(⁡k))$\\
\\$= (\log_{10}⁡(1+1) - \log_{10}(⁡1)) + (\log_{10}⁡(2+1)- \log_{10}(⁡2)) + (\log_{10}⁡(3+1)- \log_{10}(⁡3))$\\
\\$+ (\log_{10}⁡(4+1) - \log_{10}(⁡4)) + (\log_{10}⁡(5+1)- \log_{10}(⁡5)) + (\log_{10}⁡(6+1)- \log_{10}(⁡6))$\\
\\$+ (\log_{10}⁡(7+1) - \log_{10}(⁡7)) + (\log_{10}⁡(8+1) - \log_{10}⁡(8)) + (\log_{10}⁡(9+1)- \log_{10}⁡(9))$\\ 
\\$= \log_{10}(2) - \log_{10}(⁡1) + \log_{10}(3) - \log_{10}(⁡2) + \log_{10}(4) - \log_{10}⁡(3)$\\
\\$+ \log_{10}(5) - \log_{10}(⁡4) + \log_{10}(6) - \log_{10}⁡(5) + \log_{10}(7) - \log_{10}⁡(6)$\\
\\$+ \log_{10}(8) - \log_{10}(⁡7) + \log_{10}(9) - \log_{10}⁡(8) + \log_{10}(10) - \log_{10}(⁡9)$\\  
\\$= \log_{10}⁡(10) - \log_{10}(1) = 1 - 0 = 1;$\\

Thus, $\Sigma_{k=1}^9 P_k = 1$, so we got one of the proofs. 
\end{proof} 

Then, let us graph the function $y = \log_{10}(⁡x)$ on the domain from $1$ to $10$. Figure~\vref{Benfords_Probability_1} shows the graph. The graph can be used to find the values of $\log_{10}(⁡k)$ and $\log_{10}⁡(k+1)$ where $k$ is a whole number from $1$ to $9$. Both values can be found on the $y$-axis. After picking a particular $k$, subtracting $\log_{10}⁡(k)$ from $\log_{10}⁡(k+1)$ will give us the $P_k$ value according to Formula~\vref{e_m_formula}. The graph shows the differences between the values $\log_{10}⁡(k+1)$ and $\log_{10}⁡(k)$ for all nine values of $k$. All the differences are equal to their corresponding $P_k$ values. They are labeled with the corresponding colors that are taken from the ruler described above. 

The graph explains that the sum of all nine $P_k$ values of the exponential function total to $1$. It is another proof of the fact that satisfies one more condition of probabilities that we were looking for earlier.

Let us now calculate the $P_k$ values of the exponential function and list them in Table~\vref{P_k_Exp}. As it was mentioned above, all the $P_k$ values are positive numbers less than $1$ and they total to the number $1$. In addition, smaller $k$ values have higher $P_k$ values.

\begin{table}[hbt]
\caption{$P_k$ values of the exponential function}
\label{P_k_Exp}
\centering
\begin{tabular}{cc}
\toprule
Digit & $P_k$ for $y = a^x$\\
\midrule
1 &    0.30103000  \\
2 &    0.17609126  \\
3 &    0.12493874  \\
4 &    0.09691001  \\
5 &    0.07918125  \\
6 &    0.06694679  \\
7 &    0.05799195  \\
8 &    0.05115252  \\
9 &    0.04575749  \\
\bottomrule
Sum & 1.00000000  \\
\bottomrule
\end{tabular}
\end{table}

Let us figure out how the facts of multiplying the variable $x$ by its coefficient and multiplying the entire function by another coefficient changes the $fdp$ function.

Let $y = ha^{mx}$  where $h$, $a$, and $m$ are constants and $h$ > $0$, $a$ > $1$, $m$ > $0$.
Next, we will switch $x$ and $y$.

$$x = ha^{my}; \: \frac{x}{h} = a^{my}; \: my = \log_a (\frac {x}{h}); \: y = \frac {1}{m}\log_a (\frac {x}{h}); \: f(x)^{-1} = \frac {1}{m}\log_a (\frac {x}{h})$$

We will use the same technique for finding the $P_k$ values for the function $y = ha^{mx}$ that was used before to find the $P_k$ values for the $y = a^x$ function. Thus, we will start with Formula~\vref{s_formula}.

$$ P_k = \frac{f^{-1}((k+1)\cdot 10^{n-1}) - f^{-1}(k \cdot 10^{n-1})}{f^{-1}(10^n) - f^{-1}(10^{n-1})}= \frac{\frac {1}{m}\log_a (\frac {(k+1)\cdot 10^{n-1}}{h}) - \frac {1}{m}\log_a (\frac {k \cdot 10^{n-1}}{h})}{\frac {1}{m}\log_a (\frac {10^n}{h}) - \frac {1}{m}\log_a (\frac {10^{n-1}}{h})}$$

$$= \frac{\frac {1}{m}(\log_a (\frac {(k+1)\cdot 10^{n-1}}{h}) - \log_a (\frac {k \cdot 10^{n-1}}{h}))}{\frac {1}{m}(\log_a (\frac {10^n}{h}) - \log_a (\frac {10^{n-1}}{h}))} = \frac{\log_a (\frac {(k+1)\cdot 10^{n-1}}{h}) - \log_a (\frac {k \cdot 10^{n-1}}{h})}{\log_a (\frac {10^n}{h}) - \log_a (\frac {10^{n-1}}{h})}$$

$$= \frac{\log_a ((k+1)\cdot 10^{n-1}) - \log_a (h) - \log_a (k \cdot 10^{n-1}) + \log_a (h)}{\log_a (10^n) - \log_a (h) - \log_a (10^{n-1}) + \log_a (h)}$$ 

$$= \frac{\log_a ((k+1)\cdot 10^{n-1}) - \log_a (k \cdot 10^{n-1})}{\log_a (10^n) - \log_a (10^{n-1})};$$
$$ $$

The last fraction shown above was used as the first step for getting the Benford's formula. It means that if we repeat the same steps as we already have done, we will derive the exact same Formula~\vref{e_formula} for our function $y = ha^{mx}$. Thus, both constants $h$ and $m$ have no effects on $fdp$ of the exponential function $y = ha^{mx}$. The $P_k$ values will stay the same which means they will satisfy Benford's law. \cite{r_Benford_2}

Let us think now about how the exponential function can model a real-life situation. We know that there is a group of problems like compound interest, population growth, bacteria growth, and others that could be used to fulfill our needs. However, the bacteria growth problem will be the most suitable for us due to that fact that bacteria grow much faster than money in a bank or population in a particular area. Dealing with an exponential function, which grows too slow, would prevent us from getting accurate results for the $P_k$ values, unless we want to analyze our function during an incredibly long time interval. 

\FloatBarrier
\subsubsection{A real-life example}
We are talking about bacterial growth. The formula that will be used is

\begin{equation}
N(t) = N_0 e^{rt}
\label{bac_growth}
\end{equation}

\noindent
where $N(t)$ is the final amount of bacteria, $N_0$ is the initial amount of bacteria, $e \approx 2.71$, $r$ is the rate of bacterial growth, and $t$ is time in hours. We should assign values to $N_0$ and $r$. For our bacterial growth problem, let the initial number of bacteria be $300$ and the rate of growth be $40\%$. Thus, $N_0 = 300$ and $r = 0.4$. 

\begin{equation}
N(t) = 300e^{0.04t}\label{bac_growth_example}
\end{equation}

\begin{figure}[h]
\centering 
\includegraphics[width=1\columnwidth]{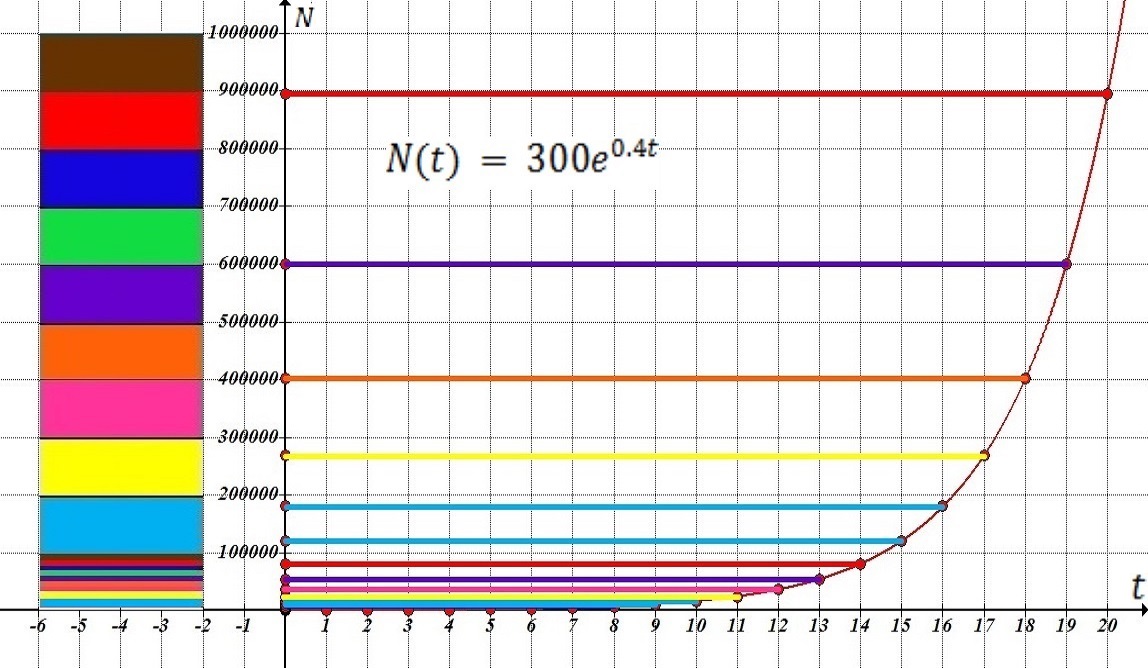} 
\caption{The function $N(t) = 300e^{0.04t}$}   
\label{Bacteria_Growth}
\end{figure}

Let us graph the function and then we can do some calculations. The graph describing bacterial growth during first $20$ hours is shown in Figure~\vref{Bacteria_Growth}.

If we want to calculate the number of bacteria regularly and to be able to make a correct conclusion about the numbers, we need to do the calculations with the same probability. We will choose the time interval of one hour and calculate the number of bacteria for the first $200$ hours. At least that many calculations are needed because we need to analyze as many as nine different probabilities. If we choose a relatively small number of time intervals, we will prevent ourselves from getting accurate probabilities. 

Let us use Formula~\vref{Bacteria_Growth} and record all the calculations of the value $N(t)$ into Table~\vref{N(t)}. Even though the caluclated $N(t)$ values are approximated, they all have the correct first digits. In order to make our results fair, we will not include the last two values to our calculations because both of them contain numbers that start a new section in the range and that section is not fully represented.

\newpage
\begin{longtable}[c]{llllll} 
 \caption{Calculated $N(t)$ values and their first digits}
 \label{N(t)}\\
\toprule
$t$- & $N$- & First & $t$- & $N$- & First \\
values & values & digit & values & values & digit \\ 
\midrule
\endfirsthead 

\toprule
$t$- & $N$- & First & $t$- & $N$- & First \\
values & values & digit & values & values & digit \\ 
\midrule
\endhead 
 1	&	4.48E+02	&	4	&	51	&	2.17E+11	&	2	\\
2	&	6.68E+02	&	6	&	52	&	3.24E+11	&	3	\\
3	&	9.96E+02	&	9	&	53	&	4.83E+11	&	4	\\
4	&	1.49E+03	&	1	&	54	&	7.21E+11	&	7	\\
5	&	2.22E+03	&	2	&	55	&	1.08E+12	&	1	\\
6	&	3.31E+03	&	3	&	56	&	1.60E+12	&	1	\\
7	&	4.93E+03	&	4	&	57	&	2.39E+12	&	2	\\
8	&	7.36E+03	&	7	&	58	&	3.57E+12	&	3	\\
9	&	1.10E+04	&	1	&	59	&	5.33E+12	&	5	\\
10	&	1.64E+04	&	1	&	60	&	7.95E+12	&	7	\\
11	&	2.44E+04	&	2	&	61	&	1.19E+13	&	1	\\
12	&	3.65E+04	&	3	&	62	&	1.77E+13	&	1	\\
13	&	5.44E+04	&	5	&	63	&	2.64E+13	&	2	\\
14	&	8.11E+04	&	8	&	64	&	3.94E+13	&	3	\\
15	&	1.21E+05	&	1	&	65	&	5.87E+13	&	5	\\
16	&	1.81E+05	&	1	&	66	&	8.76E+13	&	8	\\
17	&	2.69E+05	&	2	&	67	&	1.31E+14	&	1	\\
18	&	4.02E+05	&	4	&	68	&	1.95E+14	&	1	\\
19	&	5.99E+05	&	5	&	69	&	2.91E+14	&	2	\\
20	&	8.94E+05	&	8	&	70	&	4.34E+14	&	4	\\
21	&	1.33E+06	&	1	&	71	&	6.47E+14	&	6	\\
22	&	1.99E+06	&	1	&	72	&	9.66E+14	&	9	\\
23	&	2.97E+06	&	2	&	73	&	1.44E+15	&	1	\\
24	&	4.43E+06	&	4	&	74	&	2.15E+15	&	2	\\
25	&	6.61E+06	&	6	&	75	&	3.21E+15	&	3	\\
26	&	9.86E+06	&	9	&	76	&	4.78E+15	&	4	\\
27	&	1.47E+07	&	1	&	77	&	7.13E+15	&	7	\\
28	&	2.19E+07	&	2	&	78	&	1.06E+16	&	1	\\
29	&	3.27E+07	&	3	&	79	&	1.59E+16	&	1	\\
30	&	4.88E+07	&	4	&	80	&	2.37E+16	&	2	\\
31	&	7.28E+07	&	7	&	81	&	3.53E+16	&	3	\\
32	&	1.09E+08	&	1	&	82	&	5.27E+16	&	5	\\
33	&	1.62E+08	&	1	&	83	&	7.86E+16	&	7	\\
34	&	2.42E+08	&	2	&	84	&	1.17E+17	&	1	\\
35	&	3.61E+08	&	3	&	85	&	1.75E+17	&	1	\\
36	&	5.38E+08	&	5	&	86	&	2.61E+17	&	2	\\
37	&	8.03E+08	&	8	&	87	&	3.90E+17	&	3	\\
38	&	1.20E+09	&	1	&	88	&	5.81E+17	&	5	\\
39	&	1.79E+09	&	1	&	89	&	8.67E+17	&	8	\\
40	&	2.67E+09	&	2	&	90	&	1.29E+18	&	1	\\
41	&	3.98E+09	&	3	&	91	&	1.93E+18	&	1	\\
42	&	5.93E+09	&	5	&	92	&	2.88E+18	&	2	\\
43	&	8.85E+09	&	8	&	93	&	4.29E+18	&	4	\\
44	&	1.32E+10	&	1	&	94	&	6.41E+18	&	6	\\
45	&	1.97E+10	&	1	&	95	&	9.56E+18	&	9	\\
46	&	2.94E+10	&	2	&	96	&	1.43E+19	&	1	\\
47	&	4.38E+10	&	4	&	97	&	2.13E+19	&	2	\\
48	&	6.54E+10	&	6	&	98	&	3.17E+19	&	3	\\
49	&	9.76E+10	&	9	&	99	&	4.73E+19	&	4	\\
50	&	1.46E+11	&	1	&	100	&	7.06E+19	&	7	\\ \\ \\
101	&	1.05E+20	&	1	&	151	&	5.11E+28	&	5	\\
102	&	1.57E+20	&	1	&	152	&	7.62E+28	&	7	\\
103	&	2.34E+20	&	2	&	153	&	1.14E+29	&	1	\\
104	&	3.50E+20	&	3	&	154	&	1.70E+29	&	1	\\
105	&	5.22E+20	&	5	&	155	&	2.53E+29	&	2	\\
106	&	7.78E+20	&	7	&	156	&	3.78E+29	&	3	\\
107	&	1.16E+21	&	1	&	157	&	5.63E+29	&	5	\\
108	&	1.73E+21	&	1	&	158	&	8.40E+29	&	8	\\
109	&	2.58E+21	&	2	&	159	&	1.25E+30	&	1	\\
110	&	3.86E+21	&	3	&	160	&	1.87E+30	&	1	\\
111	&	5.75E+21	&	5	&	161	&	2.79E+30	&	2	\\
112	&	8.58E+21	&	8	&	162	&	4.16E+30	&	4	\\
113	&	1.28E+22	&	1	&	163	&	6.21E+30	&	6	\\
114	&	1.91E+22	&	1	&	164	&	9.26E+30	&	9	\\
115	&	2.85E+22	&	2	&	165	&	1.38E+31	&	1	\\
116	&	4.25E+22	&	4	&	166	&	2.06E+31	&	2	\\
117	&	6.34E+22	&	6	&	167	&	3.08E+31	&	3	\\
118	&	9.46E+22	&	9	&	168	&	4.59E+31	&	4	\\
119	&	1.41E+23	&	1	&	169	&	6.85E+31	&	6	\\
120	&	2.11E+23	&	2	&	170	&	1.02E+32	&	1	\\
121	&	3.14E+23	&	3	&	171	&	1.52E+32	&	1	\\
122	&	4.68E+23	&	4	&	172	&	2.27E+32	&	2	\\
123	&	6.99E+23	&	6	&	173	&	3.39E+32	&	3	\\
124	&	1.04E+24	&	1	&	174	&	5.06E+32	&	5	\\
125	&	1.56E+24	&	1	&	175	&	7.55E+32	&	7	\\
126	&	2.32E+24	&	2	&	176	&	1.13E+33	&	1	\\
127	&	3.46E+24	&	3	&	177	&	1.68E+33	&	1	\\
128	&	5.16E+24	&	5	&	178	&	2.51E+33	&	2	\\
129	&	7.70E+24	&	7	&	179	&	3.74E+33	&	3	\\
130	&	1.15E+25	&	1	&	180	&	5.58E+33	&	5	\\
131	&	1.71E+25	&	1	&	181	&	8.32E+33	&	8	\\
132	&	2.56E+25	&	2	&	182	&	1.24E+34	&	1	\\
133	&	3.82E+25	&	3	&	183	&	1.85E+34	&	1	\\
134	&	5.69E+25	&	5	&	184	&	2.76E+34	&	2	\\
135	&	8.49E+25	&	8	&	185	&	4.12E+34	&	4	\\
136	&	1.27E+26	&	1	&	186	&	6.15E+34	&	6	\\
137	&	1.89E+26	&	1	&	187	&	9.17E+34	&	9	\\
138	&	2.82E+26	&	2	&	188	&	1.37E+35	&	1	\\
139	&	4.21E+26	&	4	&	189	&	2.04E+35	&	2	\\
140	&	6.27E+26	&	6	&	190	&	3.04E+35	&	3	\\
141	&	9.36E+26	&	9	&	191	&	4.54E+35	&	4	\\
142	&	1.40E+27	&	1	&	192	&	6.78E+35	&	6	\\
143	&	2.08E+27	&	2	&	193	&	1.01E+36	&	1	\\
144	&	3.11E+27	&	3	&	194	&	1.51E+36	&	1	\\
145	&	4.64E+27	&	4	&	195	&	2.25E+36	&	2	\\
146	&	6.92E+27	&	6	&	196	&	3.36E+36	&	3	\\
147	&	1.03E+28	&	1	&	197	&	5.01E+36	&	5	\\
148	&	1.54E+28	&	1	&	198	&	7.47E+36	&	7	\\
149	&	2.30E+28	&	2	&	199	&	1.11E+37	&		\\
150	&	3.43E+28	&	3	&	200	&	1.66E+37	&		\\ 
\bottomrule \\ \\
\end{longtable}

Let us calculate the first digit probabilities now. Table~\vref{exp_analysis} has the analysis, which was done using Table~\vref{N(t)} and Benford's probabilities for making a comparison. As we can see, the probabilities from our table are very close to those of Benford. We can also check that the sum of all totals is $198$ and the sum of each column with probabilities is equal to $1$. 

\begin{table}[hbt]
\caption{Analysis of $N(t)$ values first digits}
\label{exp_analysis}
\centering
\begin{tabular}{cccc}
\toprule
Digit & Count & Frequency & $P_k$ values \\  
      &       & (Count / Sum) & for $y = a^x$ \\ 
\midrule
1	&	59	&	0.29797980	&	0.30103000	\\
2	&	34	&	0.17171717	&	0.17609126	\\
3	&	25	&	0.12626263	&	0.12493874	\\
4	&	19	&	0.09595960	&	0.09691001	\\
5	&	17	&	0.08585859	&	0.07918125	\\
6	&	13	&	0.06565657	&	0.06694679	\\
7	&	12	&	0.06060606	&	0.05799195	\\
8	&	10	&	0.05050505	&	0.05115252	\\
9	&	9	&	0.04545455	&	0.04575749	\\
\bottomrule
Sum	&	198	&	1.00000000	&	1.00000000	\\
\bottomrule
\end{tabular}
\end{table}

\newpage
\subsection{$P_k$ values of the power function $y = x^a$}
\subsubsection{Deriving the $P_k$ formula}

Let us talk about the function $y = x^a$ where $a$ is a natural number greater than $1$. Again, we will look at the function on the range $[1, \infty)$ for the reasons explained in the introduction. In order to choose the correct domain, we will pick the interval $[1, \infty)$, so it will match the range mentioned above for this function.
As we know, the functions of that kind with even exponents have shapes different from shapes of those with odd exponents. However, they all look more similar to each other on the domain that we will use for all of the functions $y = x^a$.

\begin{figure}[h]
\centering 
\includegraphics[width=1\columnwidth]{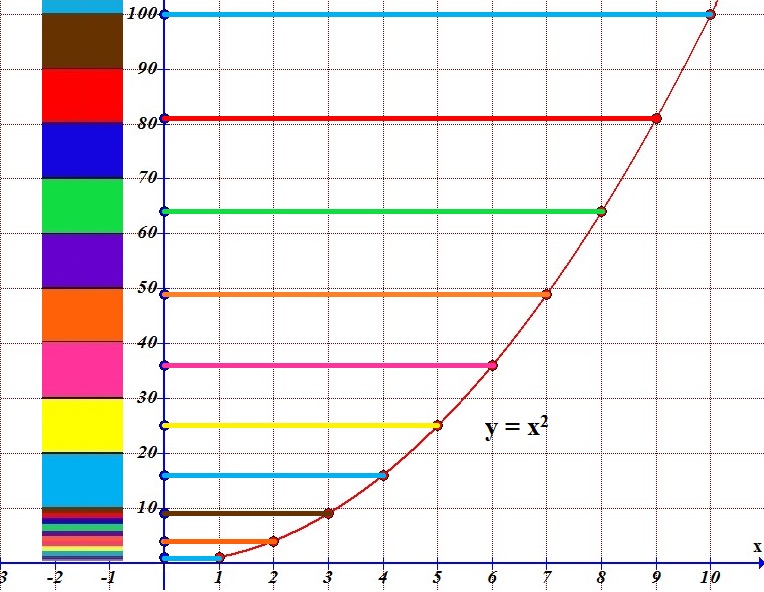} 
\caption{Quadratic function $y = x^2$ on the range  $[1, 100)$}   
\label{Quadratic_Smaller}
\end{figure}

First, let $a = 2$. The graph of the function $y = x^2$ on the limited range $[1, 100)$ is shown in Figure~\vref{Quadratic_Smaller}. The ruler that was used earlier is used again to measure the probabilities. 

\FloatBarrier
According to the graph, all nine digits have different probabilities in this range, which is correct. However, the continuous quadratic function was approximated to the discrete one with $x$-axis steps of one unit. That is why more work needs to be done before we can prove that the smaller digits have higher probabilities than the bigger ones. For instance, it might look like the digit $8$ has a higher probability than the digit $7$ and so on. However, this is not correct. The reason is that the $x$-axis steps of one unit for this particular domain and range are too big to show that the smaller the digit is, the higher its probability. 

\begin{figure}[h]
\centering 
\includegraphics[width=1\columnwidth]{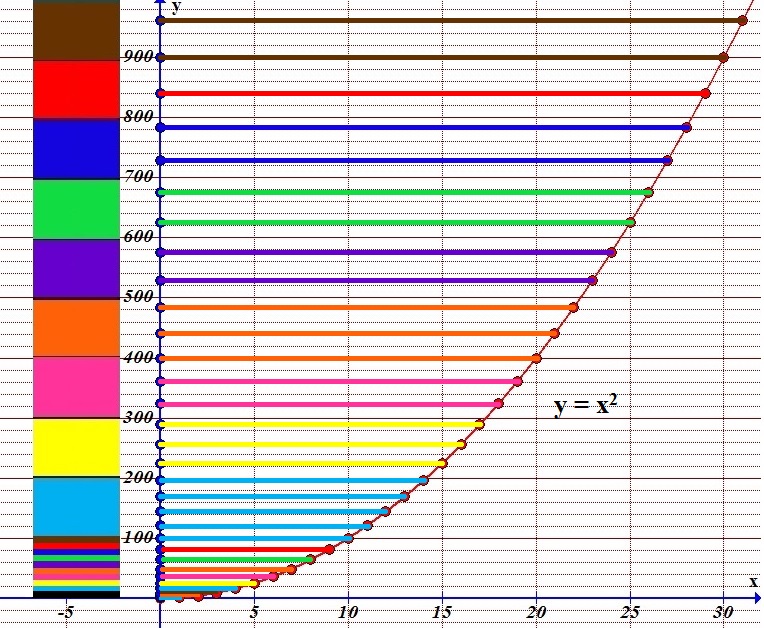} 
\caption{Quadratic function $y = x^2$ on the range  $[1, 1000)$}   
\label{Quadratic_Larger}
\end{figure}

Let us look at another graph of the function $y = x^2$. Figure~\vref{Quadratic_Larger} shows this function on the range $[1, 1000)$. As we can see, all the digits have different probabilities on that range as well as on the previous one. Unfortunately, the $x$-axis steps are still not small enough to allow us to make any specific conclusion about all $P_k$ values on that range as well. 

\begin{figure}[h]
\centering 
\includegraphics[width=1\columnwidth]{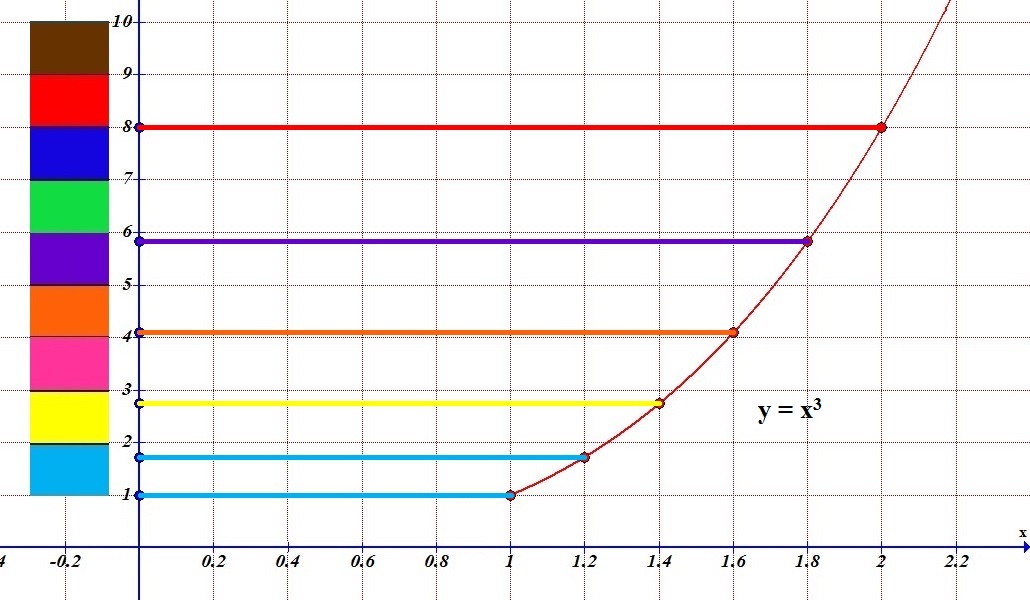} 
\caption{Cubic function $y = x^3$ on the range  $[1, 10)$}   
\label{Cube_Smaller}
\end{figure}

Let us look now at the function $y = x^3$. Figure~\vref{Cube_Smaller} shows the function on the range $[1, 10)$. Again, all the digits have different probabilities in that interval. 

\begin{figure}[h]
\centering 
\includegraphics[width=1\columnwidth]{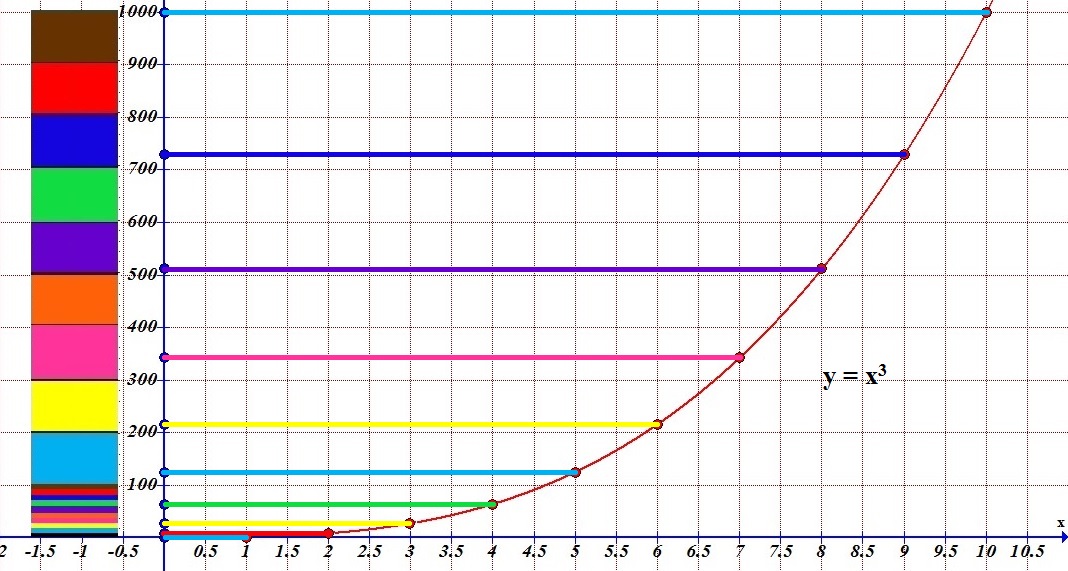} 
\caption{Cubic function $y = x^3$ on the range  $[1, 1000)$}   
\label{Cube_Larger}
\end{figure}

\FloatBarrier
The same function is shown on the range $[1, 1000)$ in Figure~\vref{Cube_Larger}. The new range will provide us with the same conclusion: different digits have different probabilities in this interval as well.

\FloatBarrier
After we saw the functions $y = x^a$ with different exponents on different ranges, let us derive the  $P_k$ formula for the function. The Formula~\vref{s_formula} should be used to accomplish our goal. We should also remind ourselves the restrictions:  $1 \leq x$ < $\infty$ and $a$ is a natural number greater than $1$.

Let $f(x) = x^a$. This function is a one-to-one function on the domain $[0, \infty)$ regardless of whether the exponent $a$ is even or odd. Thus, the function has an inverse, which is $f^{-1}(x) =x ^\frac{1}{a}$. Next, 

$$
P_k 
= \frac{((k+1)\cdot 10^{n-1})^\frac{1}{a} - (k \cdot 10^{n-1})^\frac{1}{a}}{(10^n)^{\frac{1}{a}} - (10^{n-1})^{\frac{1}{a}}}$$
$$ = \frac{(10^{n-1})^{\frac{1}{a}}((k+1)^\frac{1}{a} - k^{\frac{1}{a}})}{(10^{n-1})^{\frac{1}{a}}(10^{\frac{1}{a}}-1)}
= \frac{(k+1)^\frac{1}{a} - k^{\frac{1}{a}}}{10^{\frac{1}{a}}-1};
$$

Thus,
\begin{equation}
P_k = \frac{(k+1)^\frac{1}{a} - k^{\frac{1}{a}}}{10^{\frac{1}{a}}-1}\label{pow_formula}
\end{equation}
$$ $$

Let us analyze the result. As we can see, the final formula does not contain the exponent $n$ or any other references to the function’s range. It means that the function will have the same $fdp$ on any range. However, the formula contains the exponent $a$. Therefore, we can conclude that the $fdp$ of the quadratic function $y = x^2$ is different from the $fdp$ of the cubic function $y = x^3$ or any other function of that nature that has a different exponent. \\

Let us prove that the sum of all the $P_k$ of the function equals to $1$.
\begin{proof} 
$\Sigma^9_{k=1}P_k = \Sigma^9_{k=1}\frac{(k+1)^\frac{1}{a} - k^{\frac{1}{a}}}{10^{\frac{1}{a}}-1}$\\
\\$ = \frac{(1+1)^\frac{1}{a} - 1^{\frac{1}{a}}}{10^{\frac{1}{a}}-1} + \frac{(2+1)^\frac{1}{a} - 2^{\frac{1}{a}}}{10^{\frac{1}{a}}-1}+ \frac{(3+1)^\frac{1}{a} - 3^{\frac{1}{a}}}{10^{\frac{1}{a}}-1} + \frac{(4+1)^\frac{1}{a} - 4^{\frac{1}{a}}}{10^{\frac{1}{a}}-1} + \frac{(5+1)^\frac{1}{a} - 5^{\frac{1}{a}}}{10^{\frac{1}{a}}-1}\\ + \frac{(6+1)^\frac{1}{a} - 6^{\frac{1}{a}}}{10^{\frac{1}{a}}-1}$
$+ \frac{(7+1)^\frac{1}{a} - 7^{\frac{1}{a}}}{10^{\frac{1}{a}}-1} + \frac{(8+1)^\frac{1}{a} - 8^{\frac{1}{a}}}{10^{\frac{1}{a}}-1} + \frac{(9+1)^\frac{1}{a} - 9^{\frac{1}{a}}}{10^{\frac{1}{a}}-1}$\\
\\$= \frac{2^\frac{1}{a} - 1^{\frac{1}{a}} + 3^\frac{1}{a} - 2^{\frac{1}{a}}+ 4^\frac{1}{a} - 3^{\frac{1}{a}} + 5^\frac{1}{a} - 4^{\frac{1}{a}} + 6^\frac{1}{a} - 5^{\frac{1}{a}} + 7^\frac{1}{a} - 6^{\frac{1}{a}} + 8^\frac{1}{a} - 7^{\frac{1}{a}} + 9^\frac{1}{a} - 8^{\frac{1}{a}} + 10^\frac{1}{a} - 9^{\frac{1}{a}}}{10^{\frac{1}{a}}-1}$\\
\\$= \frac{10^\frac{1}{a} - 1^\frac{1}{a}}{10^\frac{1}{a}-1} = \frac{10^\frac{1}{a} - 1}{10^\frac{1}{a}-1} = 1;$\\

Thus, $\Sigma^9_{k=1}P_k = 1$, which satisfies one of the probability’s law.
\end{proof} 

Since the probabilities of both functions $y = x^2$ and $y = x^3$ are analyzed in this work, let us now derive the $P_k$ formulas for each of them.

First, let $y = x^2$. Then $a = 2$ and

\begin{equation}
P_k = \frac{(k+1)^\frac{1}{2} - k^{\frac{1}{2}}}{10^{\frac{1}{2}}-1} 
    = \frac{\sqrt{k+1} - \sqrt{k}}{\sqrt{10}-1} \approx \frac{\sqrt{k+1} - \sqrt{k}}{2.16}
\label{quad_formula} 
\end{equation}

\begin{figure}[h]
\centering 
\includegraphics[width=1\columnwidth]{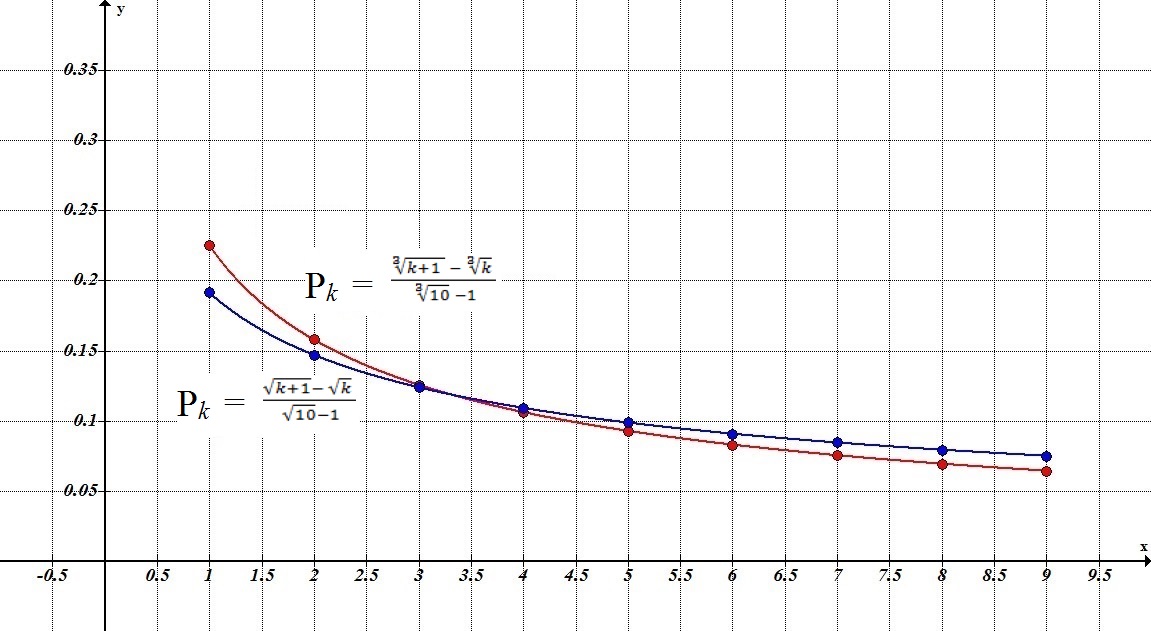} 
\caption{$P_k$ values of the functions $y=x^2$ and $y=x^3$}   
\label{Pk_Square_Cube}
\end{figure}

Next, let $y = x^3$. Then $a = 3$ and

\begin{equation}
P_k = \frac{(k+1)^\frac{1}{3} - k^{\frac{1}{3}}}{10^{\frac{1}{3}}-1} 
    = \frac{\sqrt[3]{k+1} - \sqrt[3]{k}}{\sqrt{10}-1} \approx \frac{\sqrt[3]{k+1} - \sqrt[3]{k}}{1.15}
      \label{cube_formula}
\end{equation}
\\
Let us recall that $k$ is a whole number from $1$ to $9$. As it was done in the previous section, we should graph both functions derived above as discrete ones on the domain $[1, 9]$. Both graphs are shown in Figure~\vref{Pk_Square_Cube}.

\FloatBarrier
The graphs show that the smaller the first digit is, the higher its probability and that all the probabilities are positive numbers less than $1$. The last fact satisfies the probability rules. 

We already had a chance to check whether the sum of all nine probabilities totals to $1$. Let us find one more way to prove the same fact; however, we will use graphs this time.

Since the Formula~\vref{pow_formula} tells us that the value of $P_k$ depends on the exponent $a$, we should sketch one graph proving that $\Sigma_{k=1}^9 P_k = 1$ for the function 
$y = x^2$ and another graph proving that $\Sigma_{k=1}^9 P_k = 1$ for the function $y = x^3$.

\begin{figure}[h]
\centering 
\includegraphics[width=1\columnwidth]{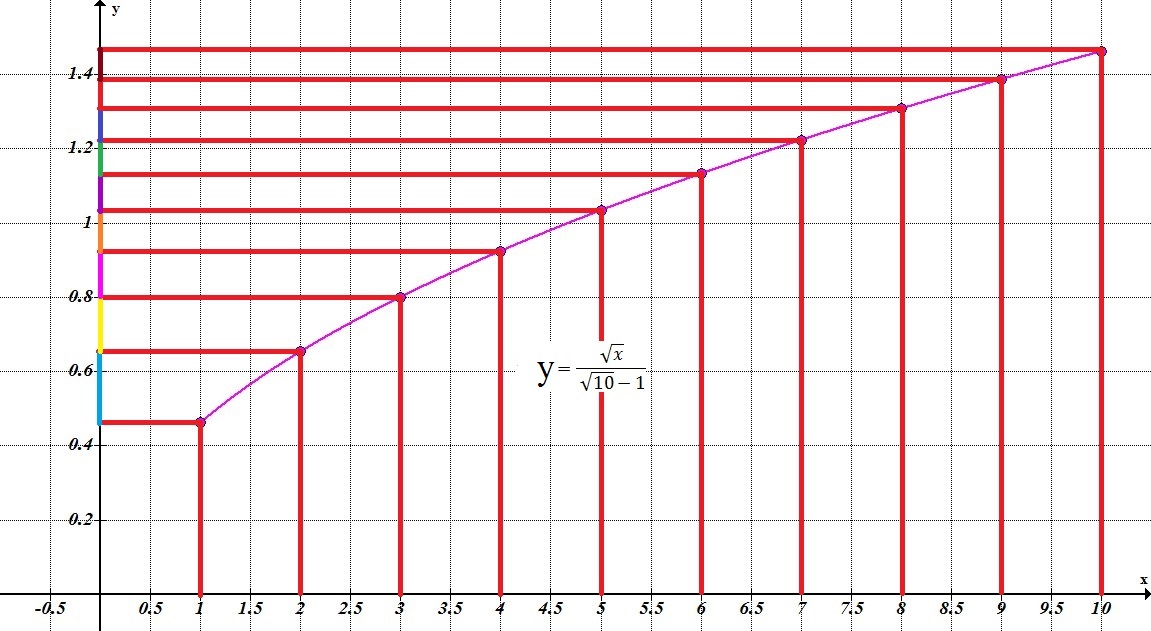} 
\caption{The function $f(x)=\frac{\sqrt{x}}{\sqrt{10} - 1}$}   
\label{Quadratic_Probability_2}
\end{figure}

First, let $a = 2$. We will use the Formula~\vref{quad_formula}.  The function $f(x)=\frac{\sqrt{x}}{\sqrt{10} - 1}$ on the domain $[1, 10]$ is shown in Figure~\vref{Quadratic_Probability_2}.

The graph can be used to find the values of $\frac{\sqrt{x}}{\sqrt{10} - 1}$ and $\frac{\sqrt{x + 1}}{\sqrt{10} - 1}$ where $k$ is a whole number from $1$ to $9$. Both of those values can be found on the $y$-axis. After picking a particular $k$, subtracting $\frac{\sqrt{x}}{\sqrt{10} - 1}$ from $\frac{\sqrt{x + 1}}{\sqrt{10} - 1}$ will give us the corresponding $P_k$ value. The graph shows the differences between the values $\frac{\sqrt{x + 1}}{\sqrt{10} - 1}$ and $\frac{\sqrt{x}}{\sqrt{10} - 1}$ for all nine values of $k$. All the differences are equal to the corresponding $P_k$ values. They are labeled with the corresponding colors that are taken from the ruler described above.

In addition, the graph clearly shows that the sum of all nine $P_k$ values of the function equals to the difference 

$$f(10) - f(1)  = \frac{\sqrt{10}}{\sqrt{10} - 1} - \frac{\sqrt{1}}{\sqrt{10} - 1} =  \frac{\sqrt{10} - 1}{\sqrt{10} - 1} = 1$$ 

Thus, we have proved once again that $\Sigma_{k=1}^9 P_k = 1$ for the function $y = x^2$. Let us now do the same for the function $y = x^3$.

\begin{figure}[h]
\centering 
\includegraphics[width=1\columnwidth]{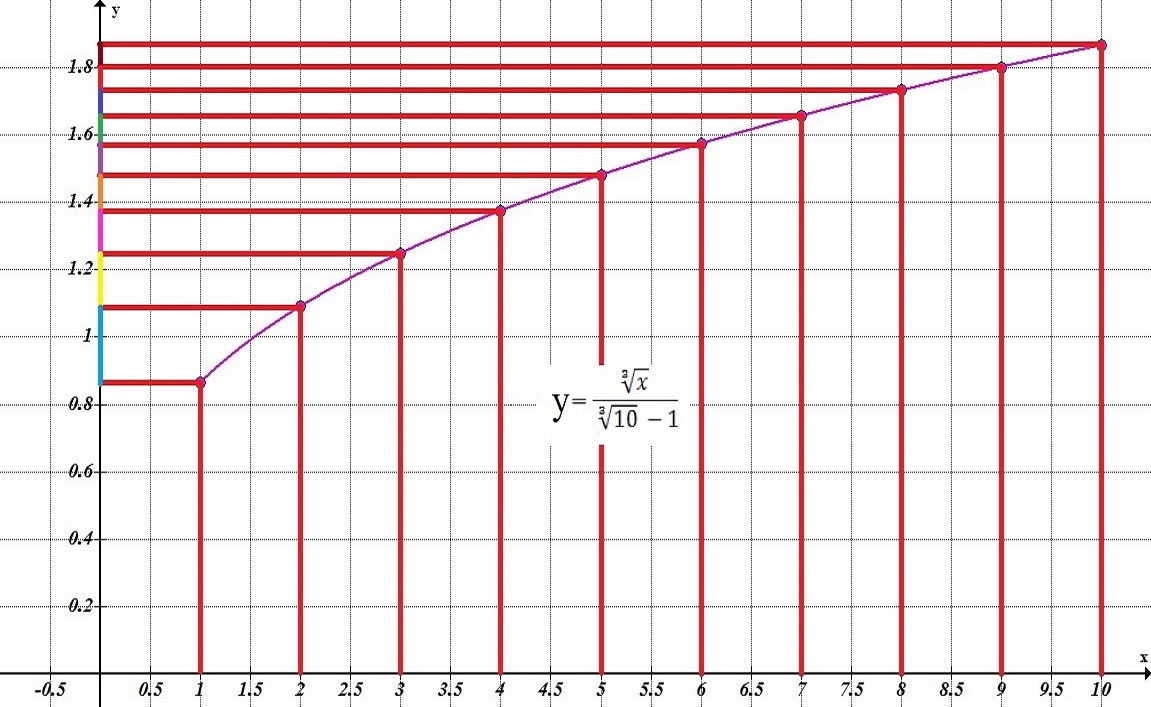} 
\caption{The function $f(x)=\frac{\sqrt[3]{x}}{\sqrt[3]{10} - 1}$}   
\label{Cube_Probability_2}
\end{figure}

\FloatBarrier
Let $a = 3$. The Formula~\vref{cube_formula} will be used. The function $f(x)=\frac{\sqrt[3]{x}}{\sqrt[3]{10} - 1}$ on the domain $[1, 10]$ is shown in Figure~\vref{Cube_Probability_2}.

Using the same technique that was used for the square function, we can see that the sum of all nine $P_k$ values of the cubic function equals to the difference 

$$f(10) - f(1)  = \frac{\sqrt[3]{10}}{\sqrt[3]{10} - 1} - \frac{\sqrt[3]{1}}{\sqrt[3]{10} - 1} =  \frac{\sqrt[3]{10} - 1}{\sqrt[3]{10} - 1} = 1$$

It proves once again that $\Sigma_{k=1}^9 P_k = 1$ for the function $y = x^3$.

Let us now calculate the $P_k$ values of the quadratic and cubic functions. They are listed in Table~\vref{P_k_Pow}. As we can see, all the $P_k$ values for both functions are positive numbers less than $1$ and they total to the number $1$. Moreover, smaller $k$ values have higher $P_k$ values.

\begin{table}[hbt]
\caption{$P_k$ values for the quadratic and cubic functions}
\label{P_k_Pow}
\centering
\begin{tabular}{ccc}
\toprule
Digit & $P_k$ for $y = x^2$ & $P_k$ for $y = x^3$\\
\midrule
1&	0.19156354&	0.22515007\\
2&	0.14699187&	0.15793749\\
3&	0.12391988&	0.12573382\\
4&	0.10917561&	0.10617742\\
5&	0.09870229&	0.09281135\\
6&	0.09076613&	0.08299351\\
7&	0.08448305&	0.07542117\\
8&	0.07934822&	0.0693706\\
9&	0.07504941&	0.06440457\\
\bottomrule
Sum & 1.00000000 & 1.00000000\\
\bottomrule
\end{tabular}
\end{table}

Before we talk about how the quadratic, cubic or other functions of this kind can model a real-life situation, we should check how a coefficient in front of the variable $x$ and another one in front of the whole function would affect the $P_k$ values of the function. We will do them one by one by considering two cases: $y = (mx)^a$ and $y = hx^a$ where $h$ and $m$ are positive constants.

First, let $y = (mx)^a$. Then, $y = m^ax^a$. As we can see, $m^a$ is a coefficient in front of the whole function and it falls into the $y = hx^a$ category. Therefore, we do not need to analyze them separately. Let us just talk about $y = hx^a$ because it looks easier.

Thus, $y = hx^a$. Next, we will switch $x$ and $y$:

$$x = hy^a; \: \frac{x}{h} = y^a; \: y = (\frac{x}{h})^{\frac{1}{a}}; \: 
f(x)^{-1} = (\frac{x}{h})^{\frac{1}{a}}$$

$$P_k = \frac{f^{-1}((k+1)\cdot 10^{n-1}) - f^{-1}(k \cdot 10^{n-1})}{f^{-1}(10 \cdot 10^{n-1}) - f^{-1}(1 \cdot 10^{n-1})} 
= \frac{(\frac{(k+1)\cdot 10^{n-1}}{h})^{\frac{1}{a}} - (\frac{k \cdot 10^{n-1}}{h})^{\frac{1}{a}}}{(\frac{10 \cdot 10^{n-1}}{h})^{\frac{1}{a}} - (\frac{1 \cdot 10^{n-1}}{h})^{\frac{1}{a}}}$$
$$= \frac{((\frac{10^{n-1}}{h})^{\frac{1}{a}})((k+1)^{\frac{1}{a}} - k^{\frac{1}{a}})}{((\frac{10^{n-1}}{h})^{\frac{1}{a}})(10^{\frac{1}{a}} - 1)} 
= \frac{(k+1)^{\frac{1}{a}} - k^{\frac{1}{a}}}{10^{\frac{1}{a}} - 1}; $$
$$ $$
In conclusion, we came back to the Formula~\vref{pow_formula}. It means that the $P_k$ values of the function $y = hx^a$ are the same as the $P_k$ values of the function $y = x^a$.

\subsubsection{A real-life example}

We are ready to talk about an actual real-life problem. Let us look at the motion of an object, released from a height of $1000$ feet. If we disregard the air resistance, the distance of the object from the initial point can be calculated by a formula 

\begin{equation}
D(t) =  \frac{1}{2}gt^2 
\label{free_fall}
\end{equation}

\noindent
where $g = 32.2 ft/s^2$ and $t$ is time in seconds. Thus, 

\begin{equation}
D(t) = 16.1t^2
\label{free_fall_example}
\end{equation}

It means that we will be dealing with the function $y = hx^2$, where $y = D$, $h = 16.1$, and $x = t$. 

\begin{figure}[h]
\centering 
\includegraphics[width=1\columnwidth]{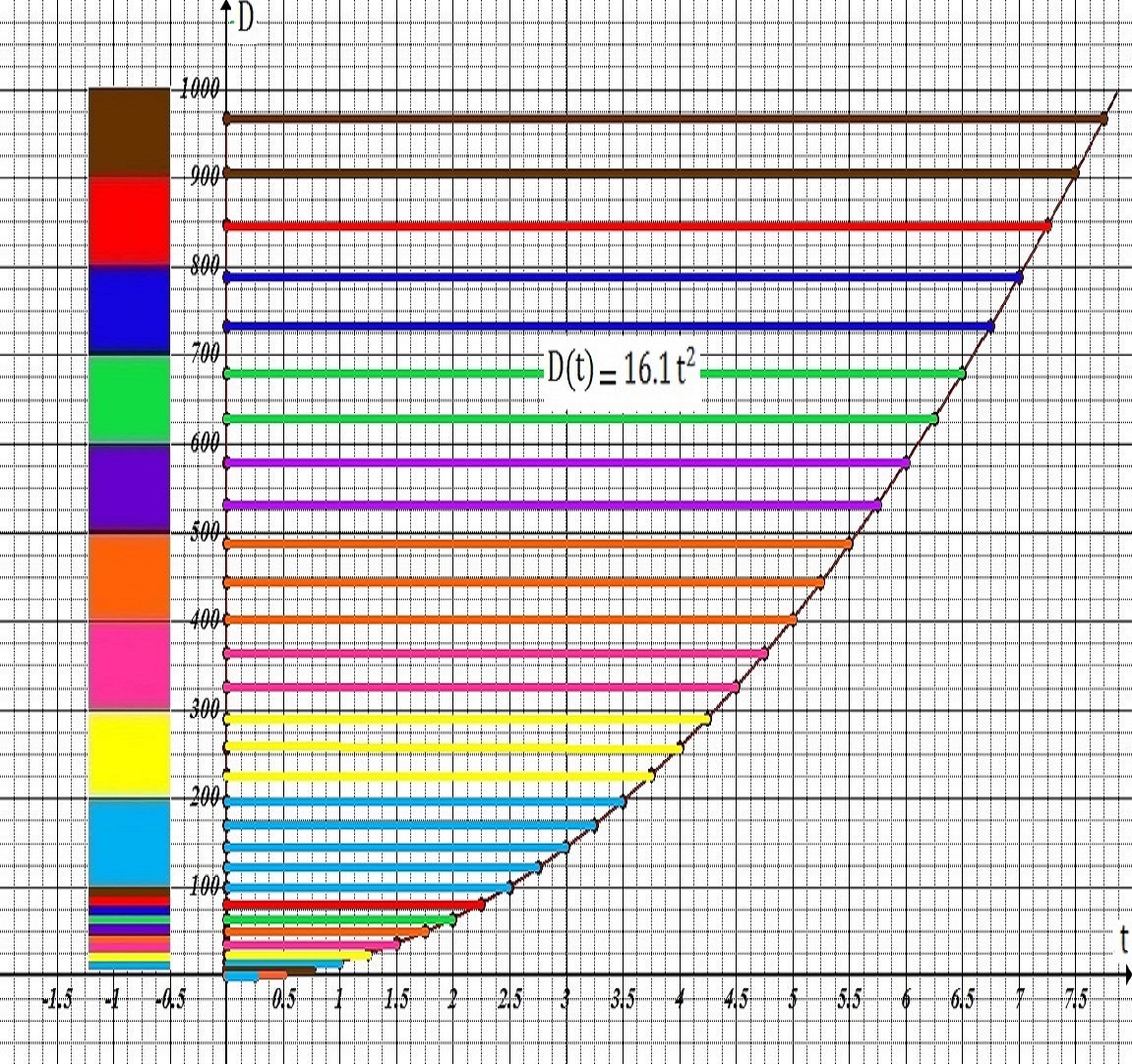} 
\caption{The function $D(t) = 16.1t^2$}   
\label{Free_Fall}
\end{figure}

The graph in Figure~\vref{Free_Fall} shows how the distance is changing through the time.

Let us calculate how many seconds it will take the object to reach the ground after falling from $1000$ feet height. First, $1000 = 16.1t^2$. Then, $t^2  \approx 62.1$. Thus,  $t \approx 7.9$. 

Since the object will touch the ground in almost $8$ seconds, we will calculate the distance of an object from the initial point every twentieth portion of a second. It will give us enough numbers for our table. 

\FloatBarrier
Let us look at Table~\vref{D(t)}. The first column shows how many seconds have passed after the object was released. The second column shows the distance that the object have gone through after the corresponding time period have passed. The third column shows the first digit of the number representing the distance. 

Since four first values in the distance column do not belong to the range of $[1, \infty)$ that we normally use, we will not do anything about the first digits of these values. In addition, the last calculated distance of $1004.801$ exceeds the limit of $1000$ feet that is why we will not do anything about the first digit of this number as well.

Now, let us analyze the first digit probabilities. Table~\vref{D(t)} calculations and the quadratic function probabilities for a comparison are shown in Table~\vref{dist_analysis}. 
The calculated free fall probabilities are very close to those of the quadratic function. If we had calculated the distance of the object more frequently than we did, we would get probabilities that are more accurate. 

\newpage
\begin{longtable}[c]{lllllllll} 
\caption{Calculated $D(t)$ values and their first digits}
\label{D(t)}\\
\toprule
$t$- & $D$- & First & $t$- & $D$- & First & $t$- & $D$- & First \\
values & values & digit & values & values & digit & values & values & digit \\ 
\midrule
\endfirsthead 
\hline
\toprule
$t$- & $D$- & First & $t$- & $D$- & First & $t$- & $D$- & First \\
values & values & digit & values & values & digit & values & values & digit \\ 
\midrule
\endhead 
0.25	&	1.01	&	1	&	2.8	&	126.22	&	1	&	5.35	&	460.82	&	4	\\
0.3	&	1.45	&	1	&	2.85	&	130.77	&	1	&	5.4	&	469.48	&	4	\\
0.35	&	1.97	&	1	&	2.9	&	135.40	&	1	&	5.45	&	478.21	&	4	\\
0.4	&	2.58	&	2	&	2.95	&	140.11	&	1	&	5.5	&	487.03	&	4	\\
0.45	&	3.26	&	3	&	3	&	144.90	&	1	&	5.55	&	495.92	&	4	\\
0.5	&	4.03	&	4	&	3.05	&	149.77	&	1	&	5.6	&	504.90	&	5	\\
0.55	&	4.87	&	4	&	3.1	&	154.72	&	1	&	5.65	&	513.95	&	5	\\
0.6	&	5.80	&	5	&	3.15	&	159.75	&	1	&	5.7	&	523.09	&	5	\\
0.65	&	6.80	&	6	&	3.2	&	164.86	&	1	&	5.75	&	532.31	&	5	\\
0.7	&	7.89	&	7	&	3.25	&	170.06	&	1	&	5.8	&	541.60	&	5	\\
0.75	&	9.06	&	9	&	3.3	&	175.33	&	1	&	5.85	&	550.98	&	5	\\
0.8	&	10.30	&	1	&	3.35	&	180.68	&	1	&	5.9	&	560.44	&	5	\\
0.85	&	11.63	&	1	&	3.4	&	186.12	&	1	&	5.95	&	569.98	&	5	\\
0.9	&	13.04	&	1	&	3.45	&	191.63	&	1	&	6	&	579.60	&	5	\\
0.95	&	14.53	&	1	&	3.5	&	197.23	&	1	&	6.05	&	589.30	&	5	\\
1	&	16.10	&	1	&	3.55	&	202.90	&	2	&	6.1	&	599.08	&	5	\\
1.05	&	17.75	&	1	&	3.6	&	208.66	&	2	&	6.15	&	608.94	&	6	\\
1.1	&	19.48	&	1	&	3.65	&	214.49	&	2	&	6.2	&	618.88	&	6	\\
1.15	&	21.29	&	2	&	3.7	&	220.41	&	2	&	6.25	&	628.91	&	6	\\
1.2	&	23.18	&	2	&	3.75	&	226.41	&	2	&	6.3	&	639.01	&	6	\\
1.25	&	25.16	&	2	&	3.8	&	232.48	&	2	&	6.35	&	649.19	&	6	\\
1.3	&	27.21	&	2	&	3.85	&	238.64	&	2	&	6.4	&	659.46	&	6	\\
1.35	&	29.34	&	2	&	3.9	&	244.88	&	2	&	6.45	&	669.80	&	6	\\
1.4	&	31.56	&	3	&	3.95	&	251.20	&	2	&	6.5	&	680.23	&	6	\\
1.45	&	33.85	&	3	&	4	&	257.60	&	2	&	6.55	&	690.73	&	6	\\
1.5	&	36.23	&	3	&	4.05	&	264.08	&	2	&	6.6	&	701.32	&	7	\\
1.55	&	38.68	&	3	&	4.1	&	270.64	&	2	&	6.65	&	711.98	&	7	\\
1.6	&	41.22	&	4	&	4.15	&	277.28	&	2	&	6.7	&	722.73	&	7	\\
1.65	&	43.83	&	4	&	4.2	&	284.00	&	2	&	6.75	&	733.56	&	7	\\
1.7	&	46.53	&	4	&	4.25	&	290.81	&	2	&	6.8	&	744.46	&	7	\\
1.75	&	49.31	&	4	&	4.3	&	297.69	&	2	&	6.85	&	755.45	&	7	\\
1.8	&	52.16	&	5	&	4.35	&	304.65	&	3	&	6.9	&	766.52	&	7	\\
1.85	&	55.10	&	5	&	4.4	&	311.70	&	3	&	6.95	&	777.67	&	7	\\
1.9	&	58.12	&	5	&	4.45	&	318.82	&	3	&	7	&	788.90	&	7	\\
1.95	&	61.22	&	6	&	4.5	&	326.03	&	3	&	7.05	&	800.21	&	8	\\
2	&	64.40	&	6	&	4.55	&	333.31	&	3	&	7.1	&	811.60	&	8	\\
2.05	&	67.66	&	6	&	4.6	&	340.68	&	3	&	7.15	&	823.07	&	8	\\
2.1	&	71.00	&	7	&	4.65	&	348.12	&	3	&	7.2	&	834.62	&	8	\\
2.15	&	74.42	&	7	&	4.7	&	355.65	&	3	&	7.25	&	846.26	&	8	\\
2.2	&	77.92	&	7	&	4.75	&	363.26	&	3	&	7.3	&	857.97	&	8	\\
2.25	&	81.51	&	8	&	4.8	&	370.94	&	3	&	7.35	&	869.76	&	8	\\
2.3	&	85.17	&	8	&	4.85	&	378.71	&	3	&	7.4	&	881.64	&	8	\\
2.35	&	88.91	&	8	&	4.9	&	386.56	&	3	&	7.45	&	893.59	&	8	\\
2.4	&	92.74	&	9	&	4.95	&	394.49	&	3	&	7.5	&	905.63	&	9	\\
2.45	&	96.64	&	9	&	5	&	402.50	&	4	&	7.55	&	917.74	&	9	\\
2.5	&	100.63	&	1	&	5.05	&	410.59	&	4	&	7.6	&	929.94	&	9	\\
2.55	&	104.69	&	1	&	5.1	&	418.76	&	4	&	7.65	&	942.21	&	9	\\
2.6	&	108.84	&	1	&	5.15	&	427.01	&	4	&	7.7	&	954.57	&	9	\\
2.65	&	113.06	&	1	&	5.2	&	435.34	&	4	&	7.75	&	967.01	&	9	\\
2.7	&	117.37	&	1	&	5.25	&	443.76	&	4	&	7.8	&	979.52	&	9	\\
2.75	&	121.76	&	1	&	5.3	&	452.25	&	4	&	7.85	&	992.12	&	9	\\
\bottomrule \\
\end{longtable}

\begin{table}[hbt]
\caption{Analysis of $D(t)$ values first digits}
\label{dist_analysis}
\centering
\begin{tabular}{cccc}
\toprule
Digit & Count & Frequency & $P_k$ values \\  
      &       & (Count / Sum) & for $y = x^2$ \\ 
\midrule
1	&	31	&	0.20261438	&	0.19156354	\\
2	&	22	&	0.14379085	&	0.14699187	\\
3	&	18	&	0.11764706	&	0.12391988	\\
4	&	18	&	0.11764706	&	0.10917561	\\
5	&	15	&	0.09803922	&	0.09870229	\\
6	&	13	&	0.08496732	&	0.09076613	\\
7	&	13	&	0.08496732	&	0.08448305	\\
8	&	12	&	0.07843137	&	0.07934822	\\
9	&	11	&	0.07189542	&	0.07504941	\\
\bottomrule
Sum	&	153	&	1.00000000	&	1.00000000	\\
\bottomrule
\end{tabular}
\end{table}

\newpage
\subsection{$P_k$ values of the linear function $y = mx$}
\subsubsection{Deriving the $P_k$ formula}

Let us look at the linear function $y = mx$ where $m$ > $0$.  The reason why we are choosing a positive slope is that fact that our initial goal is to apply each function in this work to a real-life situation. We should admit that a graph reflecting any such situation would take the first quadrant only. However, choosing a negative slope will prevent the graph from passing through the first quadrant. Thus, we are choosing $m$ > $0$. In addition, the slope cannot be equal to zero because otherwise we would get a horizontal line $y = 0$, which would prevent us from talking about any probability at all. All the first digits in this case would just be zeroes.

We will use that portion of the function, which does not go beyond the range $[1, \infty)$ since we agreed upon using this range for all our functions in the introduction to the work. A particular value of $m$ and the previously mentioned range will lead us to have the domain $[\frac{1}{m}, \infty)$.

\begin{figure}[h]
\centering 
\includegraphics[width=1\columnwidth]{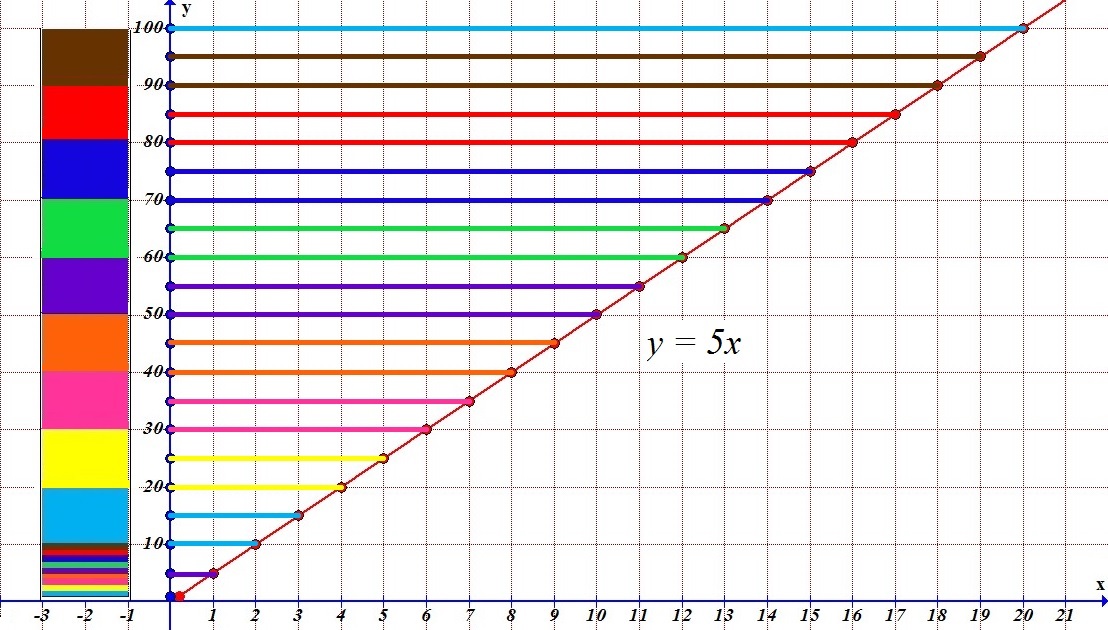} 
\caption{The function $y = 5x$}   
\label{Linear}
\end{figure}

Let us look at Figure~\vref{Linear}, which shows a graph of a linear function $y = 5x$ in the limited range from $1$ to $100$. As we can see, the graph’s direction will not be changed if we make the range smaller or bigger. Moreover, according to the colors of all the horizontal lines, an assumption that all the digits have equal probabilities can be made. We just have to prove it.

\FloatBarrier
Formula~\vref{s_formula} will be used. Let $f(x) = mx$. Then, $f^{-1}(x)= \frac{1}{m}x$. Next, 

$$P_k = \frac{\frac{1}{m}\cdot (k + 1) \cdot 10^{n-1} - \frac{1}{m}\cdot k \cdot 10^{n-1}}{\frac{1}{m}\cdot 10^n - \frac{1}{m}\cdot 10^{n-1}} $$ 
$$= \frac{\frac{1}{m}\cdot 10^{n-1}((k + 1)  - k)}{\frac{1}{m}\cdot 10^{n-1}(10  - 1)} =
\frac{k + 1 - k}{10  - 1} = \frac{1}{9};$$
Thus,
\begin{equation}
P_k =  \frac{1}{9}
\label{Linear_Prob}
\end{equation}
$$ $$

\begin{figure}[h]
\centering 
\includegraphics[width=1\columnwidth]{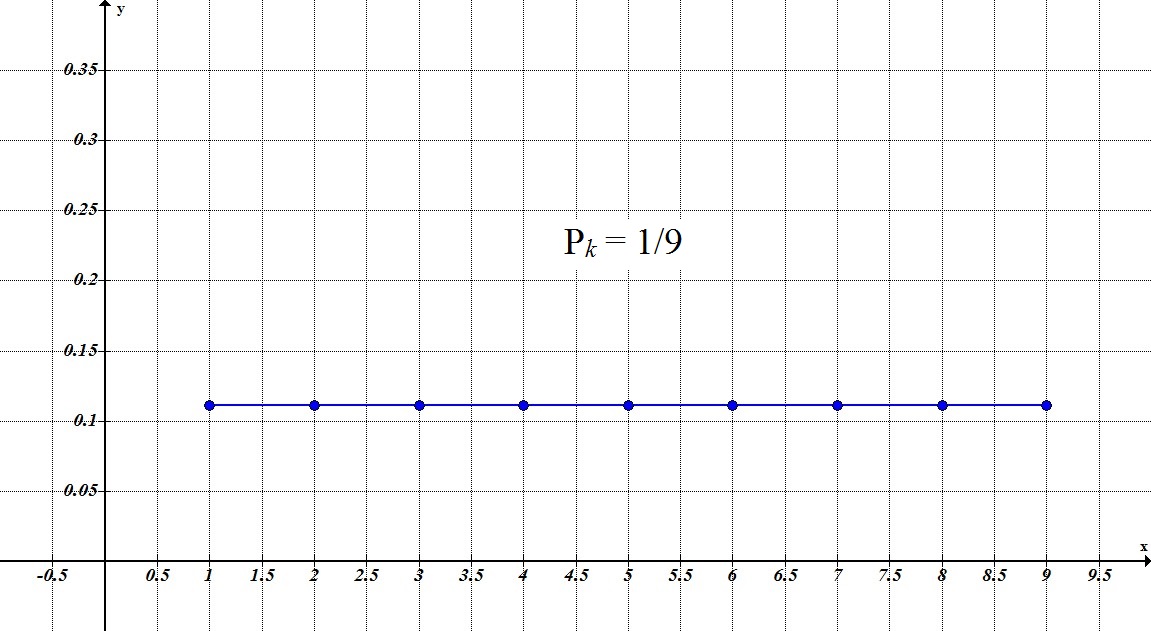} 
\caption{$P_k$ values of the function $y = mx$}   
\label{Linear_Probability}
\end{figure}

Formula~\vref{Linear_Prob} looks a bit unusual; however, it proves the assumption made by us before: all nine digits have the same probabilities regardless of the slope of the linear function or chosen range. In addition, it is very easy to check that the total of all nine $P_k$ values equals to $1$. The formula will lead us to a horizontal line if we decide to graph it. Figure~\vref{Linear_Probability} shows the probability line of a linear function.

Let us think about how a linear function can be applied to a real-life situation. 

\subsubsection{A real-life example}

Consider a swimming pool, which holds $1000 ft^3$ of water. We will fill it with water using a garden hose. After the hose is turned on, the volume of water in the pool before it is filled can be calculated with the formula

\begin{equation}
V(t)=rt
\label{Linear_Equation}
\end{equation}

\noindent
where $V(t)$ is the volume measured in $ft^3$, r is the flow rate, and $t$ is time in minutes. Let us suppose that the hose has a flow rate of $5 ft^3/min$. It means that we will be dealing with the function 

\begin{equation}
y = 5x
\label{Linear_Equation_Example}
\end{equation}

\noindent
where $y = V(t)$  and $x = t$. The graph in Figure~\vref{Linear} can be used to see how the amount of water is increasing during the first $20$ minutes.

Let us calculate how many minutes it will take to fill the entire pool with water. Let $1000 = 5t$.  Then $t = 1000/5 = 200$. Since the pool will be filled in $200$ minutes, we will calculate the amount of water in it every minute. Doing this should give us enough numbers for our table. 

Let us look at Table~\vref{V(t)}. The first column shows how many minutes have passed after the hose was turned on. The second column shows the volume of water in the pool in $ft^3$ after the corresponding number of minutes have passed. The third column shows the first digit of the number representing the volume. 

Since the table has only one $y$-value in the range $[1, 10)$, that value should be excluded from our calculations to keep them fair. In addition, the last $y$-value of our table exceeds the range of $[1, 1000)$.  Therefore, that value should be excluded from our calculations as well. We should use Table~\vref{V(t)} to analyze the first digit probabilities of the numbers representing the volume of water. Table~\vref{P_k_Linear} shows the calculations from Table~\vref{V(t)} and the linear function probabilities for making a comparison. 

\newpage
\begin{longtable}[c]{lllllllll} 
 \caption{Calculated $V(t)$ values and their first digits}
 \label{V(t)}\\
\toprule
$t$- & $V$- & First & $t$- & $V$- & First & $t$- & $V$- & First \\
values & values & digit & values & values & digit & values & values & digit \\ 
\midrule
\endfirsthead 
\toprule
$t$- & $V$- & First & $t$- & $V$- & First & $t$- & $V$- & First \\
values & values & digit & values & values & digit & values & values & digit \\ 
\midrule
\endhead 
1	&	5	&		&	51	&	255	&	2	&	101	&	505	&	5	\\
2	&	10	&	1	&	52	&	260	&	2	&	102	&	510	&	5	\\
3	&	15	&	1	&	53	&	265	&	2	&	103	&	515	&	5	\\
4	&	20	&	2	&	54	&	270	&	2	&	104	&	520	&	5	\\
5	&	25	&	2	&	55	&	275	&	2	&	105	&	525	&	5	\\
6	&	30	&	3	&	56	&	280	&	2	&	106	&	530	&	5	\\
7	&	35	&	3	&	57	&	285	&	2	&	107	&	535	&	5	\\
8	&	40	&	4	&	58	&	290	&	2	&	108	&	540	&	5	\\
9	&	45	&	4	&	59	&	295	&	2	&	109	&	545	&	5	\\
10	&	50	&	5	&	60	&	300	&	3	&	110	&	550	&	5	\\
11	&	55	&	5	&	61	&	305	&	3	&	111	&	555	&	5	\\
12	&	60	&	6	&	62	&	310	&	3	&	112	&	560	&	5	\\
13	&	65	&	6	&	63	&	315	&	3	&	113	&	565	&	5	\\
14	&	70	&	7	&	64	&	320	&	3	&	114	&	570	&	5	\\
15	&	75	&	7	&	65	&	325	&	3	&	115	&	575	&	5	\\
16	&	80	&	8	&	66	&	330	&	3	&	116	&	580	&	5	\\
17	&	85	&	8	&	67	&	335	&	3	&	117	&	585	&	5	\\
18	&	90	&	9	&	68	&	340	&	3	&	118	&	590	&	5	\\
19	&	95	&	9	&	69	&	345	&	3	&	119	&	595	&	5	\\
20	&	100	&	1	&	70	&	350	&	3	&	120	&	600	&	6	\\
21	&	105	&	1	&	71	&	355	&	3	&	121	&	605	&	6	\\
22	&	110	&	1	&	72	&	360	&	3	&	122	&	610	&	6	\\
23	&	115	&	1	&	73	&	365	&	3	&	123	&	615	&	6	\\
24	&	120	&	1	&	74	&	370	&	3	&	124	&	620	&	6	\\
25	&	125	&	1	&	75	&	375	&	3	&	125	&	625	&	6	\\
26	&	130	&	1	&	76	&	380	&	3	&	126	&	630	&	6	\\
27	&	135	&	1	&	77	&	385	&	3	&	127	&	635	&	6	\\
28	&	140	&	1	&	78	&	390	&	3	&	128	&	640	&	6	\\
29	&	145	&	1	&	79	&	395	&	3	&	129	&	645	&	6	\\
30	&	150	&	1	&	80	&	400	&	4	&	130	&	650	&	6	\\
31	&	155	&	1	&	81	&	405	&	4	&	131	&	655	&	6	\\
32	&	160	&	1	&	82	&	410	&	4	&	132	&	660	&	6	\\
33	&	165	&	1	&	83	&	415	&	4	&	133	&	665	&	6	\\
34	&	170	&	1	&	84	&	420	&	4	&	134	&	670	&	6	\\
35	&	175	&	1	&	85	&	425	&	4	&	135	&	675	&	6	\\
36	&	180	&	1	&	86	&	430	&	4	&	136	&	680	&	6	\\
37	&	185	&	1	&	87	&	435	&	4	&	137	&	685	&	6	\\
38	&	190	&	1	&	88	&	440	&	4	&	138	&	690	&	6	\\
39	&	195	&	1	&	89	&	445	&	4	&	139	&	695	&	6	\\
40	&	200	&	2	&	90	&	450	&	4	&	140	&	700	&	7	\\
41	&	205	&	2	&	91	&	455	&	4	&	141	&	705	&	7	\\
42	&	210	&	2	&	92	&	460	&	4	&	142	&	710	&	7	\\
43	&	215	&	2	&	93	&	465	&	4	&	143	&	715	&	7	\\
44	&	220	&	2	&	94	&	470	&	4	&	144	&	720	&	7	\\
45	&	225	&	2	&	95	&	475	&	4	&	145	&	725	&	7	\\
46	&	230	&	2	&	96	&	480	&	4	&	146	&	730	&	7	\\
47	&	235	&	2	&	97	&	485	&	4	&	147	&	735	&	7	\\
48	&	240	&	2	&	98	&	490	&	4	&	148	&	740	&	7	\\
49	&	245	&	2	&	99	&	495	&	4	&	149	&	745	&	7	\\
50	&	250	&	2	&	100	&	500	&	5	&	150	&	750	&	7	\\ \\ \\
151	&	755	&	7	&	168	&	840	&	8	&	185	&	925	&	9	\\
152	&	760	&	7	&	169	&	845	&	8	&	186	&	930	&	9	\\
153	&	765	&	7	&	170	&	850	&	8	&	187	&	935	&	9	\\
154	&	770	&	7	&	171	&	855	&	8	&	188	&	940	&	9	\\
155	&	775	&	7	&	172	&	860	&	8	&	189	&	945	&	9	\\
156	&	780	&	7	&	173	&	865	&	8	&	190	&	950	&	9	\\
157	&	785	&	7	&	174	&	870	&	8	&	191	&	955	&	9	\\
158	&	790	&	7	&	175	&	875	&	8	&	192	&	960	&	9	\\
159	&	795	&	7	&	176	&	880	&	8	&	193	&	965	&	9	\\
160	&	800	&	8	&	177	&	885	&	8	&	194	&	970	&	9	\\
161	&	805	&	8	&	178	&	890	&	8	&	195	&	975	&	9	\\
162	&	810	&	8	&	179	&	895	&	8	&	196	&	980	&	9	\\
163	&	815	&	8	&	180	&	900	&	9	&	197	&	985	&	9	\\
164	&	820	&	8	&	181	&	905	&	9	&	198	&	990	&	9	\\
165	&	825	&	8	&	182	&	910	&	9	&	199	&	995	&	9	\\
166	&	830	&	8	&	183	&	915	&	9	&	200	&	1000	&		\\
167	&	835	&	8	&	184	&	920	&	9	&		&		&		\\
\bottomrule \\
\end{longtable}

\begin{table}[hbt]
\caption{Analysis of $V(t)$ values first digits}
\label{P_k_Linear}
\centering
\begin{tabular}{cccc}
\toprule
Digit & Count & Frequency &  $P_k$ values \\  
      &       & (Count / Sum) & for $y = mx$ \\ 
\midrule
1	&	22	&	 1/9	&	 1/9	\\
2	&	22	&	 1/9	&	 1/9	\\
3	&	22	&	 1/9	&	 1/9	\\
4	&	22	&	 1/9	&	 1/9	\\
5	&	22	&	 1/9	&	 1/9	\\
6	&	22	&	 1/9	&	 1/9	\\
7	&	22	&	 1/9	&	 1/9	\\
8	&	22	&	 1/9	&	 1/9	\\
9	&	22	&	 1/9	&	 1/9	\\
\bottomrule
Sum	&	198	&	1.00000000	&	1.00000000	\\
\bottomrule
\end{tabular}
\end{table}

\FloatBarrier
As we can see, the calculated water volume $fdp$ values are the same as the corresponding $fdp$ values of the liner function. Moreover, all the water volume $fdp$ values are equal to each other, as it should be.

\newpage
\subsection{$P_k$ values of the root function $y = \sqrt[a]{x}$}
\subsubsection{Deriving the $P_k$ formula}

Let us talk about the function $y = \sqrt[a]{x}$ where $a$ is a natural number greater than $1$. Again, we will look at the function on the range $[1, \infty)$. Moreover, we will pick the same interval $[1, \infty)$ for the domain, so that they will both match each other.

\begin{figure}[h]
\centering 
\includegraphics[width=1\columnwidth]{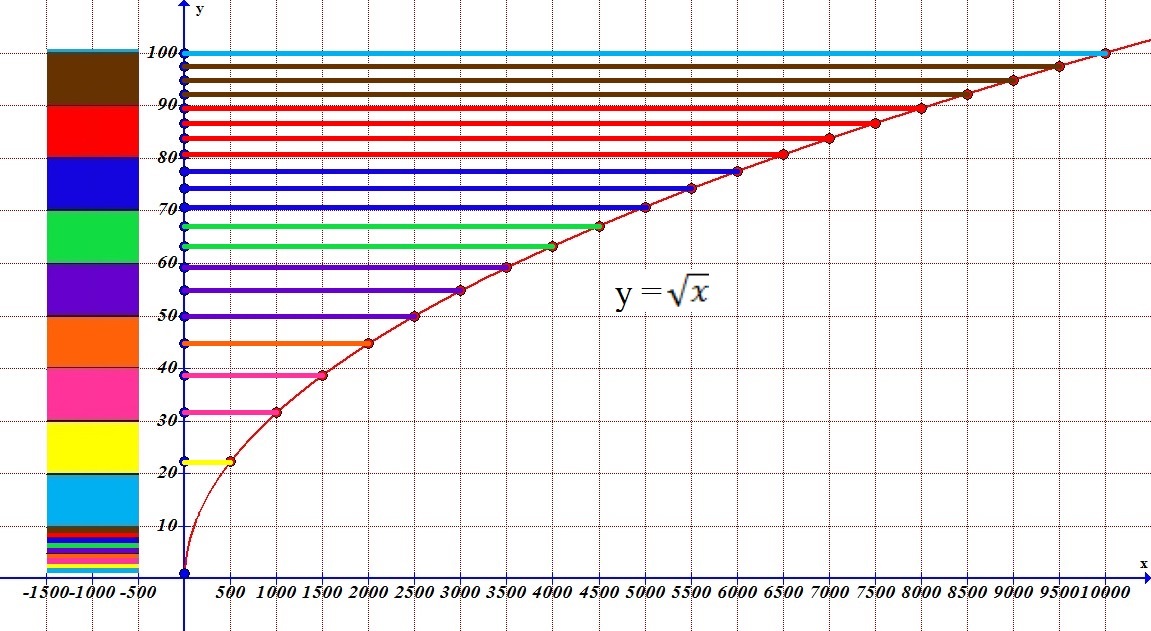} 
\caption{The function $y = \sqrt{x}$}   
\label{Square_Root}
\end{figure}

First, let $a = 2$. The graph of the function $y = \sqrt{x}$ on the limited range $[1, 100)$ is shown in Figure~\vref{Square_Root}. The same ruler is used to measure the probabilities. Looking at the graph, we can see that bigger digits have higher probabilities. However, more work needs to be done to prove it. 

\FloatBarrier
Let us look now at the function $y = \sqrt[3]{x}$. Figure~\vref{Cube_Root} shows the function on the range $[1, 100)$. Again, bigger digits have higher probabilities and we will prove it.

\begin{figure}[h]
\centering 
\includegraphics[width=1\columnwidth]{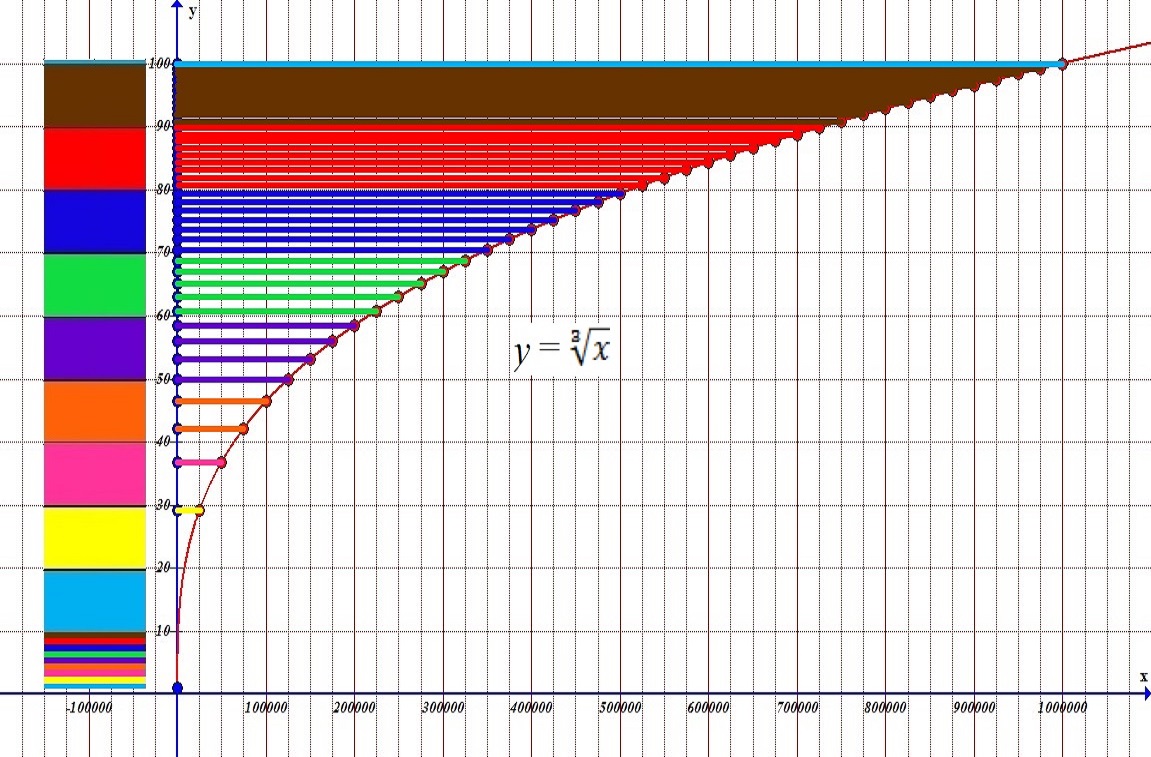} 
\caption{The function $y = \sqrt[3]{x}$}
\label{Cube_Root}
\end{figure}

\FloatBarrier
We had a chance to see the functions $y = \sqrt[a]{x}$ with different indexes $a$. Let us derive the formula for $P_k$ of the function.

Formula~\vref{s_formula} should be used. In addition, we know that $1 \leq x$ < $\infty$ and $a$ is a natural number greater than $1$.

Let $y = \sqrt[a]{x}$. This function is a one-to-one function on any part of its domain regardless of the index $a$. That is why the function has an inverse, which is  $f^{-1}(x) = x^a$. Next, 
$$
P_k 
= \frac{((k+1)\cdot 10^{n-1})^a - (k \cdot 10^{n-1})^a}{(10^n)^a - (10^{n-1})^a}$$
$$= \frac{(10^{n-1})^a((k+1)^a - k^a)}{(10^{n-1})^a(10^a-1)} 
= \frac{(k+1)^a - k^a}{10^a-1};$$

$$ $$
Consequently,
\begin{equation}
P_k = \frac{(k+1)^a - k^a}{10^a-1}
\label{root_formula}
\end{equation}
$$ $$

The final formula does not contain the exponent $n$ as well as any other references to the function’s range. It means that the function will have the same $fdp$ on any of its range. However, the formula contains the index $a$, which looks like an exponent now. Thus, we can conclude that the $fdp$ of the square root function $y = \sqrt{x}$ is different from the $fdp$ of the cube root function $y = \sqrt[3]{x}$ or any other function of that nature that has a different index.

Let us prove that the sum of all the $P_k$ of the function equals to $1$.

\begin{proof}
$\Sigma^9_{k=1}P_k = \Sigma^9_{k=1}\frac{(k+1)^a - k^a}{10^a-1}$\\
\\$ = \frac{(1+1)^a - 1^a}{10^a-1} + \frac{(2+1)^a - 2^a}{10^a-1}+ \frac{(3+1)^a - 3^a}{10^a-1} + \frac{(4+1)^a - 4^a}{10^a-1} + \frac{(5+1)^a - 5^a}{10^a - 1}$\\ 
\\$+ \frac{(6+1)^a - 6^a}{10^a-1} + \frac{(7+1)^a - 7^a}{10^a-1} + \frac{(8+1)^a - 8^a}{10^a-1} + \frac{(9+1)^a - 9^a}{10^a-1}$\\
\\$= \frac{2^a - 1^a + 3^a - 2^a+ 4^a - 3^a + 5^a - 4^a + 6^a - 5^a + 7^a - 6^a + 8^a - 7^a + 9^a - 8^a + 10^a - 9^a}{10^a-1}$\\
\\$= \frac{10^a - 1^a}{10^a-1} = \frac{10^a - 1}{10^a-1} = 1;$\\
\\
Therefore, $\Sigma^9_{k=1}P_k = 1$, which satisfies one of the rules of probability.
\end{proof}

Let us now derive the $P_k$ formulas for the square root and the cube root functions separately.

First, let $y = \sqrt{x}$. Then $a = 2$ and

$$P_k = \frac{(k + 1)^2 - k^2}{10^2 - 1} = 
\frac{k^2 + 2k + 1 - k^2}{100 - 1} = \frac{2k + 1}{99};$$

Next, let $y = \sqrt[3]{x}$. Then $a = 3$ and

$$P_k = \frac{(k + 1)^3 - k^3}{10^3 - 1} = 
\frac{k^3 + 3k^2 +3k + 1 - k^3}{1000 - 1} = \frac{3k^2 +3k + 1}{999};$$

\begin{figure}[h]
\centering 
\includegraphics[width=1\columnwidth]{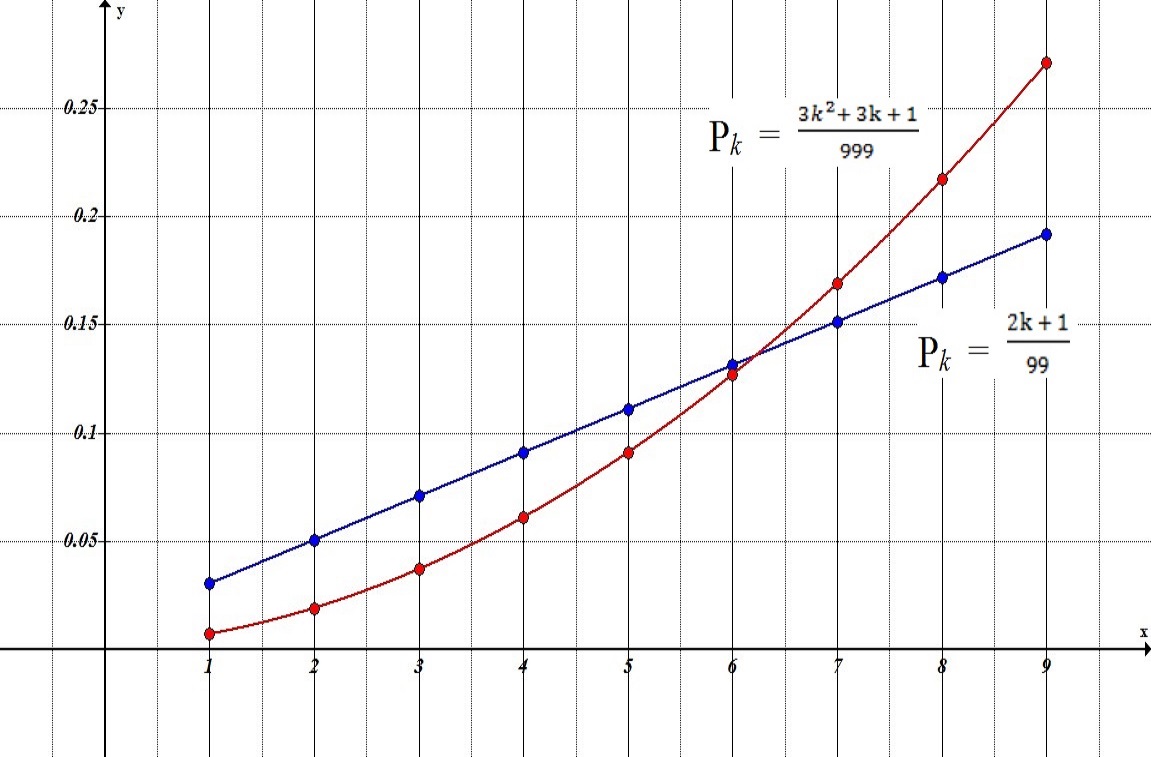} 
\caption{Functions $P_k = \frac{2k + 1}{99}$ and $P_k = \frac{3k^2 +3k + 1}{999}$}   
\label{Square_Root_Probability_1}
\end{figure}

Let us recall that $k$ is a whole number from $1$ to $9$. We should graph both functions derived above as discrete ones on the domain $[1, 9]$. The graphs of both of them are shown in Figure~\vref{Square_Root_Probability_1}.

Both graphs show that the bigger the first digit is, the higher its probability and that all the probabilities are positive numbers less than $1$. The last fact satisfies the probability rules. 

Since we already checked that the sum of all nine probabilities totals to $1$ using the formulas, let us prove the same fact using graphs this time.

$P_k$ formula depends on the value of $a$, so we should sketch one graph proving that $\Sigma^9_{k=1}P_k = 1$ for the function $y = \sqrt{x}$ and another graph proving the same thing for the function $y = \sqrt[3]{x}$.

\begin{figure}[h]
\centering 
\includegraphics[width=1\columnwidth]{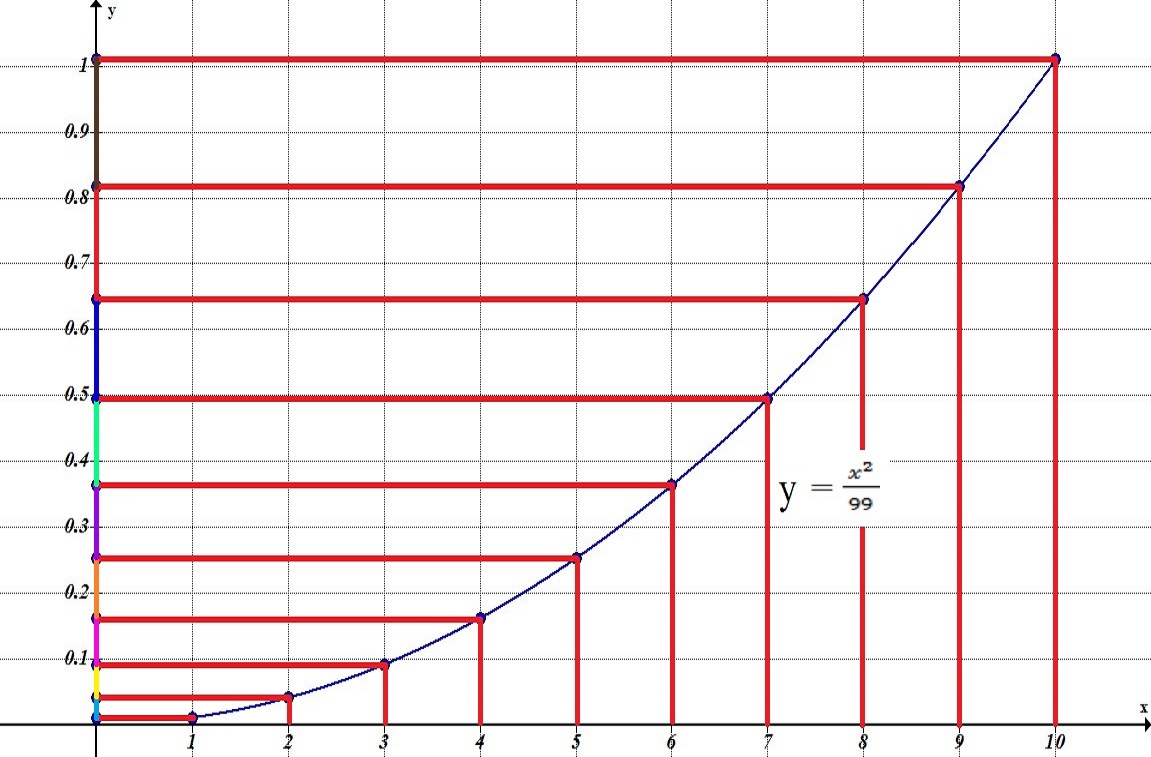} 
\caption{The function $f(x) = \frac{x^2}{99}$ on the domain $[1, 10]$}
\label{Square_Root_Probability_2}
\end{figure}

First, let $a = 2$. The original formula $P_k = \frac{(k + 1)^2 - k^2}{99}$ should be used. Let us graph the function $f(x) = \frac{x^2}{99}$ on the domain $[1, 10]$. The graph is shown in Figure~\vref{Square_Root_Probability_2}.

This graph can be used to find the values of $\frac{(k + 1)^2}{99}$ and $\frac{k^2}{99}$ where $k$ is a whole number from $1$ to $9$. Both of the values can be found on the $y$-axis. After picking a particular $k$, subtracting $\frac{k^2}{99}$ from $\frac{(k + 1)^2}{99}$ will give us the corresponding $P_k$ value. 

The graph in Figure~\vref{Square_Root_Probability_2} shows the differences between the values $\frac{(k + 1)^2}{99}$ and $\frac{k^2}{99}$ for all nine values of $k$. All the differences are equal to the corresponding $P_k$ values. They are labeled with the corresponding colors that are taken from the ruler described above.
\\

The graph clearly shows that the sum of all nine $P_k$ values of the function equals to the difference 

$$f(10) - f(1)  = \frac{100}{99} - \frac{1}{99} = \frac{99}{99} = 1$$

We have proved once again that $\Sigma^9_{k=1}P_k = 1$ for the function $y = \sqrt{x}$. 

\begin{figure}[h]
\centering 
\includegraphics[width=1\columnwidth]{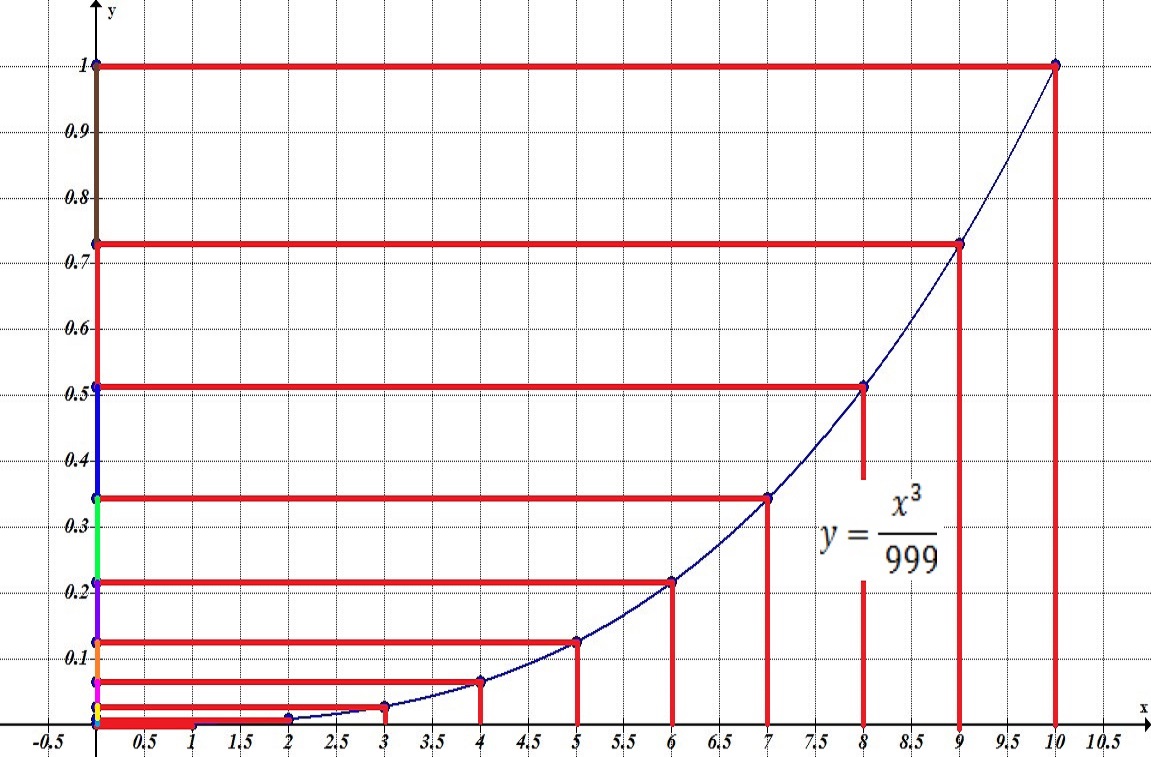} 
\caption{The function $f(x) = \frac{x^3}{999}$ on the domain $[1, 10]$}
\label{Cube_Root_Probability_2}
\end{figure}

Let us now do the same for the function $y = \sqrt[3]{x}$. Let $a = 3$. The initial formula is $P_k = \frac{(k + 1)^3 - k^3}{999}$. Let us graph the function $f(x) = \frac {x^3}{999}$ on the domain $[1, 10]$. The graph is shown in Figure~\vref{Cube_Root_Probability_2}.

Using the same technique that was used for the square root function, we can see that the sum of all nine $P_k$ values of the cubic function equals to the difference 

$$f(10) - f(1)  = \frac{1000}{999} - \frac{1}{999} = \frac{999}{999} = 1$$

It proves once again that $\Sigma^9_{k=1}P_k = 1$ for the function $y = \sqrt[3]{x}$. 

Let us now look at the $P_k$ values of the square root and the cube root functions. They are listed in Table~\vref{P_k_Root}. As we can see, bigger $k$ values have higher $P_k$ values in both cases. In addition, all the $P_k$ values for both functions are positive numbers less than $1$. If we total the exact $P_k$ values for each function (not the approximated ones with three decimal places), the totals will be equal to the number $1$. 

We also can see that for the digits from one to six, the fdp values of the square root function are higher than those of the cube root function. However, the fdp values of the digits from seven to nine are higher for the cube root function than those of the square root function.

\begin{table}[hbt]
\caption{$P_k$ values for the square root and cube root functions}
\label{P_k_Root}
\centering
\begin{tabular}{ccc}
\toprule
Digit & $P_k$ for $y = \sqrt{x}$ & $P_k$ for $y = \sqrt[3]{x}$ \\
\midrule
1	&	0.03030303	&	0.00700701	\\
2	&	0.05050505	&	0.01901902	\\
3	&	0.07070707	&	0.03703704	\\
4	&	0.09090909	&	0.06106106	\\
5	&	0.11111111	&	0.09109109	\\
6	&	0.13131313	&	0.12712713	\\
7	&	0.15151515	&	0.16916917	\\
8	&	0.17171717	&	0.21721722	\\
9	&	0.19191919	&	0.27127127	\\
\bottomrule
Sum & 1.00000000 & 1.00000000  \\
\bottomrule
\end{tabular}
\end{table}

Prior to talking about how the function $y = \sqrt[a]{x}$ can be used to model a real-life situation, let us check how a coefficient in front of the variable $x$ and another one in front of the whole function will affect the $P_k$ values of the function. 

Let us consider two cases similar to those that were done before: $y = \sqrt[a]{mx}$ and $y = h\sqrt[a]{x}$ where $h$ and $m$ are positive constants.  

First, let $y = \sqrt[a]{mx}$. Then, $y = \sqrt[a]{m} \cdot \sqrt[a]{x}$. Obviously, $\sqrt[a]{m}$ is a coefficient in front of the whole function and it falls into the  $y = h\sqrt[a]{x}$ category. Therefore, we do not need to analyze both functions separately. Let us just talk about the function 
 $y = h\sqrt[a]{x}$ because it looks easier than $y = \sqrt[a]{mx}$. Thus, $y = h\sqrt[a]{x}$. 

\FloatBarrier
Next, we will switch $x$ and $y$:

$$x = h\sqrt[a]{y}; \: \frac{x}{h} = \sqrt[a]{y}; \: y = (\frac{x}{h})^a; \: f(x)^{-1} =  (\frac{x}{h})^a$$

$$P_k = \frac{f^{-1}((k + 1) \cdot 10^{n-1}) - f^{-1}(k \cdot 10^{n-1})}{f^{-1}(10 \cdot 10^{n-1}) - f^{-1}(1 \cdot 10^{n-1})}= \frac{(\frac{(k + 1) \cdot 10^{n-1}}{h})^a - (\frac{k \cdot 10^{n-1}}{h})^a}{(\frac{10 \cdot 10^{n-1}}{h})^a - (\frac {10^{n-1}}{h})^a}$$ 
$$= \frac{(\frac{10^{n-1}}{h})^a [(k + 1)^a - k^a]}{(\frac{10^{n-1}}{h})^a (10^a - 1)}
= \frac{(k + 1)^a - k^a}{10^a - 1};$$

Thus, we came back to Formula~\vref{s_formula}. It means that the $P_k$ values of the function $y = h\sqrt[a]{x}$ are the same as the $P_k$ values of the function $y = \sqrt[a]{x}$.

\subsubsection{A real-life example}

Let us suppose that we are blowing up a spherical balloon at a rate of $4$ cubic inches per second. The volume of the balloon can be found by the formula 

\begin{equation}
V = {\frac{4}{3}}\pi R^3
\label{Cube_Root}
\end{equation}

\noindent
where $\pi \approx 3.14$ and $R$ is the balloon’s radius. We also understand that $V(t) = 4t$, where $t$ is time in seconds. 

If we isolate the radius, we get 

$$R = \sqrt[3]{\frac{3V}{4\pi}} = \sqrt[3]{\frac{12t}{4\pi}} = \sqrt[3]{\frac{3t}{\pi}}$$ 

Thus, we will use the formula 

\begin{equation}
R = \sqrt[3]{\frac{3t}{\pi}}
\label{Cube_Root_Example}
\end{equation}

\begin{figure}[h]
\centering 
\includegraphics[width=1\columnwidth]{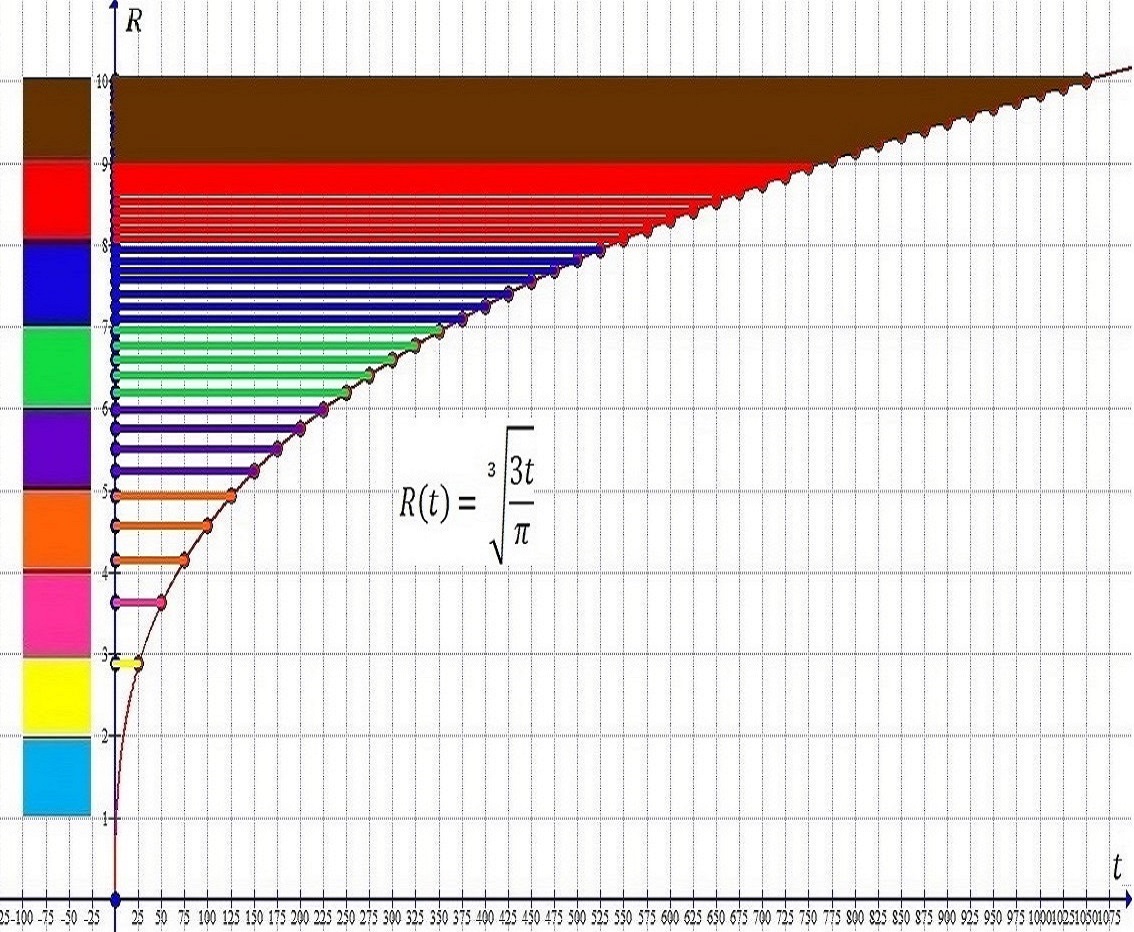} 
\caption{The function $R(t) = \sqrt[3]{\frac{3t}{\pi}}$}
\label{Radius_Growth}
\end{figure}

\FloatBarrier
Let us graph the function and see how the radius changes over time. The graph describing the radius's growth during the first $1050$ seconds is shown in Figure~\vref{Radius_Growth}.

As we can see from the graph, it takes about $1050$ seconds to make the radius equal to $10$ inches. 

Since the function grows relatively slow, we should project it on the range $[1, 10)$ only. If we calculate the radius every $5$ seconds, it will give us $210$ time intervals which should be enough for our table. 

Let us look at Table~\vref{R(t)}. The first column shows how many seconds have passed after the balloon started being blown up. The second column shows the radius of the balloon in inches after the corresponding number of seconds have passed. The third column shows the first digit of the number representing the radius. 

We will use the table to analyze the first digit probabilities of the numbers representing the radius of the balloon. As it was done before, we will exclude those values which do not satisfy our requirements.

Table~\vref{P_k_Cube_Root} shows the calculations and the cube root function probabilities for making a comparison. We can see that the probabilities from Table~\vref{R(t)} are very close to those of the cube root function. We can also check that the sum of each column with probabilities is equal to $1$. 

\newpage
\begin{longtable}[c]{lllllllll} 
 \caption{Calculated $R(t)$ values and their first digits}
 \label{R(t)}\\
\toprule
$t$- & $R$- & First & $t$- & $R$- & First & $t$- & $R$- & First \\
values & values & digit & values & values & digit & values & values & digit \\ 
\midrule
\endfirsthead 

\toprule
$t$- & $R$- & First & $t$- & $R$- & First & $t$- & $R$- & First \\
values & values & digit & values & values & digit & values & values & digit \\ 
\midrule
\endhead 

5	&	1.68	&	1	&	255	&	6.24	&	6	&	505	&	7.84	&	7	\\
10	&	2.12	&	2	&	260	&	6.29	&	6	&	510	&	7.87	&	7	\\
15	&	2.43	&	2	&	265	&	6.33	&	6	&	515	&	7.89	&	7	\\
20	&	2.67	&	2	&	270	&	6.36	&	6	&	520	&	7.92	&	7	\\
25	&	2.88	&	2	&	275	&	6.40	&	6	&	525	&	7.94	&	7	\\
30	&	3.06	&	3	&	280	&	6.44	&	6	&	530	&	7.97	&	7	\\
35	&	3.22	&	3	&	285	&	6.48	&	6	&	535	&	7.99	&	7	\\
40	&	3.37	&	3	&	290	&	6.52	&	6	&	540	&	8.02	&	8	\\
45	&	3.50	&	3	&	295	&	6.56	&	6	&	545	&	8.04	&	8	\\
50	&	3.63	&	3	&	300	&	6.59	&	6	&	550	&	8.07	&	8	\\
55	&	3.74	&	3	&	305	&	6.63	&	6	&	555	&	8.09	&	8	\\
60	&	3.86	&	3	&	310	&	6.66	&	6	&	560	&	8.12	&	8	\\
65	&	3.96	&	3	&	315	&	6.70	&	6	&	565	&	8.14	&	8	\\
70	&	4.06	&	4	&	320	&	6.74	&	6	&	570	&	8.16	&	8	\\
75	&	4.15	&	4	&	325	&	6.77	&	6	&	575	&	8.19	&	8	\\
80	&	4.24	&	4	&	330	&	6.81	&	6	&	580	&	8.21	&	8	\\
85	&	4.33	&	4	&	335	&	6.84	&	6	&	585	&	8.24	&	8	\\
90	&	4.41	&	4	&	340	&	6.87	&	6	&	590	&	8.26	&	8	\\
95	&	4.49	&	4	&	345	&	6.91	&	6	&	595	&	8.28	&	8	\\
100	&	4.57	&	4	&	350	&	6.94	&	6	&	600	&	8.31	&	8	\\
105	&	4.65	&	4	&	355	&	6.97	&	6	&	605	&	8.33	&	8	\\
110	&	4.72	&	4	&	360	&	7.01	&	7	&	610	&	8.35	&	8	\\
115	&	4.79	&	4	&	365	&	7.04	&	7	&	615	&	8.37	&	8	\\
120	&	4.86	&	4	&	370	&	7.07	&	7	&	620	&	8.40	&	8	\\
125	&	4.92	&	4	&	375	&	7.10	&	7	&	625	&	8.42	&	8	\\
130	&	4.99	&	4	&	380	&	7.13	&	7	&	630	&	8.44	&	8	\\
135	&	5.05	&	5	&	385	&	7.16	&	7	&	635	&	8.46	&	8	\\
140	&	5.11	&	5	&	390	&	7.19	&	7	&	640	&	8.49	&	8	\\
145	&	5.17	&	5	&	395	&	7.23	&	7	&	645	&	8.51	&	8	\\
150	&	5.23	&	5	&	400	&	7.26	&	7	&	650	&	8.53	&	8	\\
155	&	5.29	&	5	&	405	&	7.29	&	7	&	655	&	8.55	&	8	\\
160	&	5.35	&	5	&	410	&	7.32	&	7	&	660	&	8.57	&	8	\\
165	&	5.40	&	5	&	415	&	7.35	&	7	&	665	&	8.60	&	8	\\
170	&	5.46	&	5	&	420	&	7.37	&	7	&	670	&	8.62	&	8	\\
175	&	5.51	&	5	&	425	&	7.40	&	7	&	675	&	8.64	&	8	\\
180	&	5.56	&	5	&	430	&	7.43	&	7	&	680	&	8.66	&	8	\\
185	&	5.61	&	5	&	435	&	7.46	&	7	&	685	&	8.68	&	8	\\
190	&	5.66	&	5	&	440	&	7.49	&	7	&	690	&	8.70	&	8	\\
195	&	5.71	&	5	&	445	&	7.52	&	7	&	695	&	8.72	&	8	\\
200	&	5.76	&	5	&	450	&	7.55	&	7	&	700	&	8.74	&	8	\\
205	&	5.81	&	5	&	455	&	7.57	&	7	&	705	&	8.76	&	8	\\
210	&	5.85	&	5	&	460	&	7.60	&	7	&	710	&	8.79	&	8	\\
215	&	5.90	&	5	&	465	&	7.63	&	7	&	715	&	8.81	&	8	\\
220	&	5.94	&	5	&	470	&	7.66	&	7	&	720	&	8.83	&	8	\\
225	&	5.99	&	5	&	475	&	7.68	&	7	&	725	&	8.85	&	8	\\
230	&	6.03	&	6	&	480	&	7.71	&	7	&	730	&	8.87	&	8	\\
235	&	6.08	&	6	&	485	&	7.74	&	7	&	735	&	8.89	&	8	\\
240	&	6.12	&	6	&	490	&	7.76	&	7	&	740	&	8.91	&	8	\\
245	&	6.16	&	6	&	495	&	7.79	&	7	&	745	&	8.93	&	8	\\
250	&	6.20	&	6	&	500	&	7.82	&	7	&	750	&	8.95	&	8	\\ \\  \\
755	&	8.97	&	8	&	855	&	9.35	&	9	&	955	&	9.70	&	9	\\
760	&	8.99	&	8	&	860	&	9.36	&	9	&	960	&	9.71	&	9	\\
765	&	9.01	&	9	&	865	&	9.38	&	9	&	965	&	9.73	&	9	\\
770	&	9.03	&	9	&	870	&	9.40	&	9	&	970	&	9.75	&	9	\\
775	&	9.05	&	9	&	875	&	9.42	&	9	&	975	&	9.76	&	9	\\
780	&	9.06	&	9	&	880	&	9.44	&	9	&	980	&	9.78	&	9	\\
785	&	9.08	&	9	&	885	&	9.45	&	9	&	985	&	9.80	&	9	\\
790	&	9.10	&	9	&	890	&	9.47	&	9	&	990	&	9.81	&	9	\\
795	&	9.12	&	9	&	895	&	9.49	&	9	&	995	&	9.83	&	9	\\
800	&	9.14	&	9	&	900	&	9.51	&	9	&	1000	&	9.85	&	9	\\
805	&	9.16	&	9	&	905	&	9.53	&	9	&	1005	&	9.86	&	9	\\
810	&	9.18	&	9	&	910	&	9.54	&	9	&	1010	&	9.88	&	9	\\
815	&	9.20	&	9	&	915	&	9.56	&	9	&	1015	&	9.90	&	9	\\
820	&	9.22	&	9	&	920	&	9.58	&	9	&	1020	&	9.91	&	9	\\
825	&	9.24	&	9	&	925	&	9.59	&	9	&	1025	&	9.93	&	9	\\
830	&	9.25	&	9	&	930	&	9.61	&	9	&	1030	&	9.94	&	9	\\
835	&	9.27	&	9	&	935	&	9.63	&	9	&	1035	&	9.96	&	9	\\
840	&	9.29	&	9	&	940	&	9.65	&	9	&	1040	&	9.98	&	9	\\
845	&	9.31	&	9	&	945	&	9.66	&	9	&	1045	&	9.99	&	9	\\
850	&	9.33	&	9	&	950	&	9.68	&	9	&	1050	&	10.01	&	1	\\
\bottomrule \\
\end{longtable}

\FloatBarrier

\begin{table}[hbt]
\caption{Analysis of $R(t)$ values first digits}
\label{P_k_Cube_Root}
\centering
\begin{tabular}{cccc}
\toprule
Digit & Count & Frequency & $P_k$ values \\  
      &       & (Count / Sum) & for $y = \sqrt[3]{x}$ \\ 
\midrule
1	&	2	&	0.00952381	&	0.00700701	\\
2	&	4	&	0.019047619	&	0.01901902	\\
3	&	8	&	0.038095238	&	0.03703704	\\
4	&	13	&	0.061904762	&	0.06106106	\\
5	&	19	&	0.09047619	&	0.09109109	\\
6	&	26	&	0.123809524	&	0.12712713	\\
7	&	36	&	0.171428571	&	0.16916917	\\
8	&	45	&	0.214285714	&	0.21721722	\\
9	&	57	&	0.271428571	&	0.27127127	\\
\bottomrule
Sum	&	210	&	1.00000000	&	1.00000000	\\
\bottomrule

\end{tabular}
\end{table}

\newpage
\subsection{$P_k$ values of the logarithmic function $y = \log_a(x)$}
\subsubsection{Deriving the $P_k$ formula}

The function $y = \log_a(x)$ is the most challenging one in this chapter. 

We will talk about this function having base $a$ as a positive number greater than $1$. As it was done before with other functions, we will look at the function $y = \log_a(x)$ on the range $[1, \infty)$. In addition, we will pick an interval $[a, \infty)$ for the domain, so they both will match each other.

To be able to graph the function, we have to assign a particular value to $a$. First, let $a = 2$; then we will assign $a = 3$. 

Due to that fact that the logarithmic function with any value of $a$ grows much slower than all previous ones that we discussed before, we will graph both functions $y = \log_2⁡(x)$ and $y = \log_3(x)$ on the limited range $[1, 10)$. However, the range $[1, \infty)$ will be used for calculating the $P_k$ values initially.

\begin{figure}[h]
\centering 
\includegraphics[width=1\columnwidth]{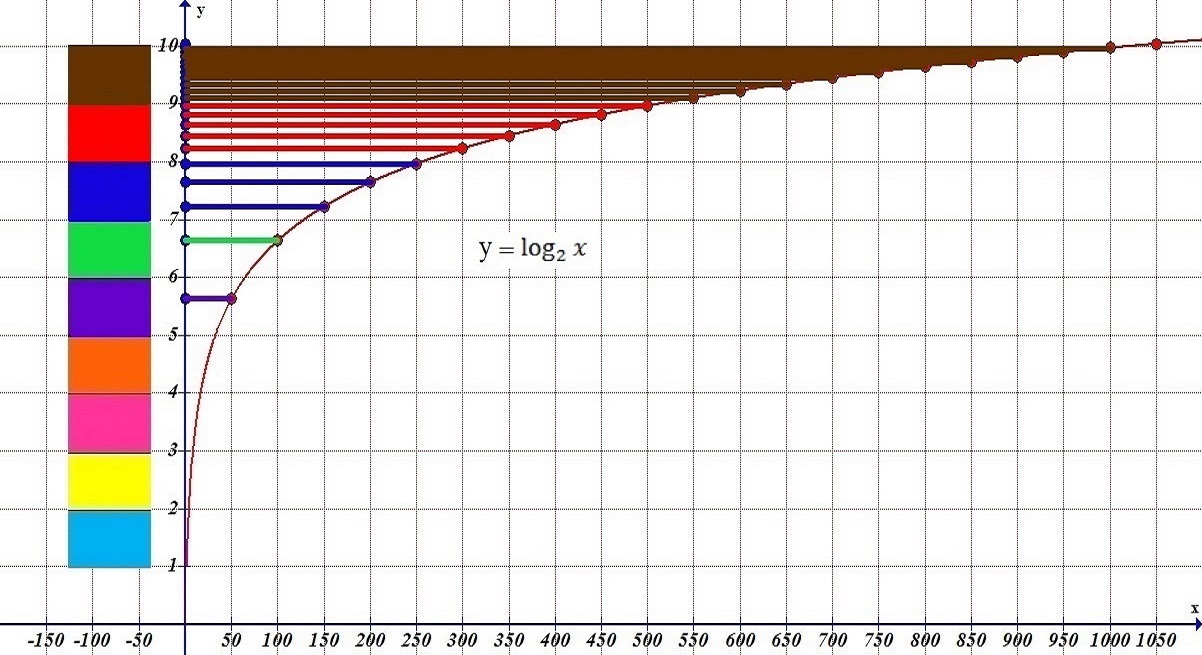} 
\caption{The function $y = \log_2(x)$}
\label{Log_Of_2}
\end{figure}

The graph of the function $y = \log_2(⁡x)$ is shown in Figure~\vref{Log_Of_2}. It is obvious that bigger digits have higher probabilities. We will do more work later to prove it. 

\begin{figure}[h]
\centering 
\includegraphics[width=1\columnwidth]{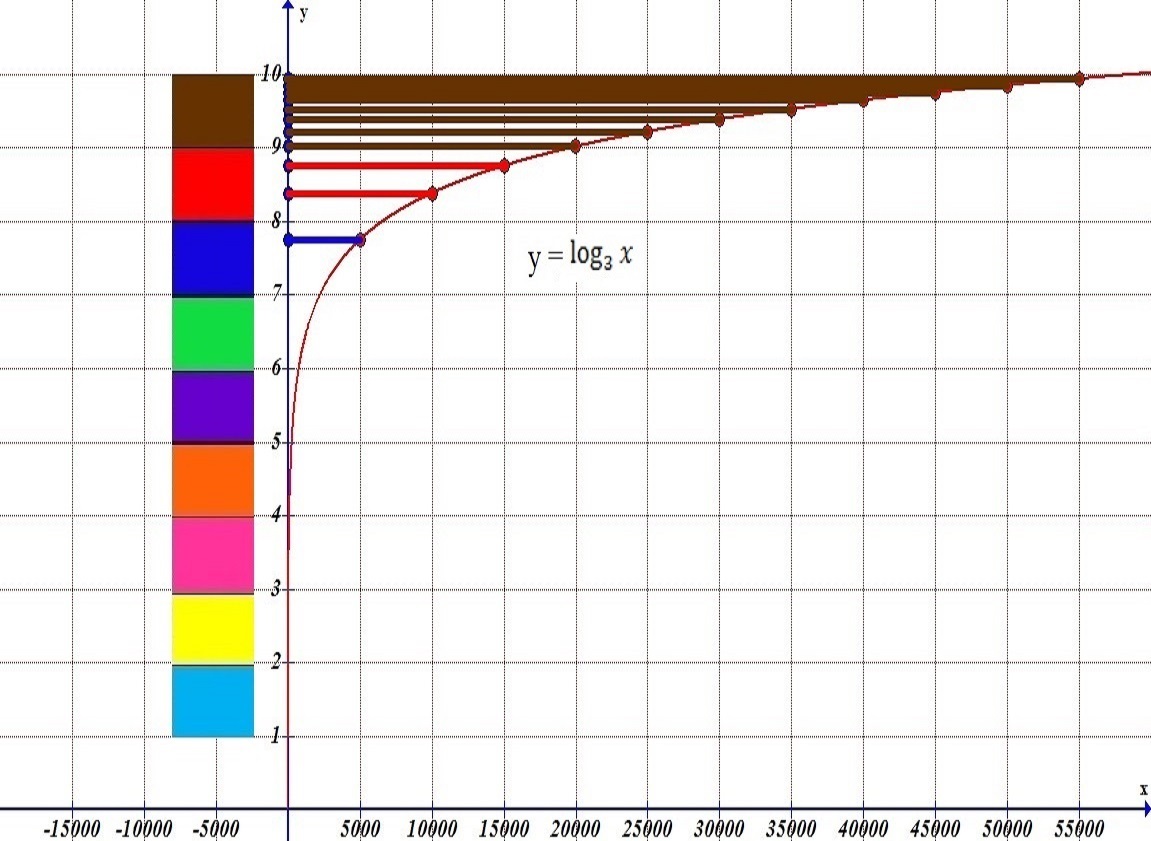} 
\caption{The function $y = \log_3(x)$}
\label{Log_Of_3}
\end{figure}

\FloatBarrier
Let us look at the function $y = \log_3⁡(x)$. Figure~\vref{Log_Of_3} shows the function and yet again, bigger digits have higher probabilities. 

Thus, let us prove it. Let $f(x) = \log_a(x)$. According to the definition of the logarithmic function, $f^{-1}(x) = a^x$. Next, 

$$P_k = \frac{a^{(k+1)\cdot10^{n-1}} - a^{k\cdot10^{n-1}}}{a^{10\cdot10^{n-1}} - a^{10^{n-1}}} = 
\frac{(a^{10^{n-1}})^{k+1} - (a^{10^{n-1}})^k}{(a^{10^{n-1}})^{10} - a^{10^{n-1}}};$$

Let us use a substitution $c = a^{10^{n-1}}$. Then,

$$\frac{(a^{10^{n-1}})^{k+1} - (a^{10^{n-1}})^k}{(a^{10^{n-1}})^{10} - a^{10^{n-1}}} = \frac{c^{k+1} - c^k}{c^{10} - c} = \frac{c^k(c - 1)}{c(c^9 - 1)};$$

Next, let us modify the denominator.

$$c^9 - 1 = (c^3)^3 - (1^3)^3 = (c^3 - 1^3 )(c^6 + c^3 + 1)=
(c - 1)(c^2 + c + 1)(c^6 + c^3 + 1)$$
$$= (c - 1)(c^8 + c^7 + c^6 + c^5 + c^4 + c^3 + c^2 + c + 1);$$

Then, we will bring it back.

$$\frac{c^k(c - 1)}{c(c^9 - 1)} = 
\frac{c^k(c - 1)}{c(c - 1)(c^8 + c^7 + c^6 + c^5 + c^4 + c^3 + c^2 + c + 1)}$$
$$=\frac{c^k}{c(c^8 + c^7 + c^6 + c^5 + c^4 + c^3 + c^2 + c + 1)}$$ 
$$=\frac{c^k}{c^9 + c^8 + c^7 + c^6 + c^5 + c^4 + c^3 + c^2 + c} =
\frac{c^k}{\sum^9_{i=1}{c^i}};$$

Let us bring our variables $a$ and $n$ back.

$$\frac{c^k}{\sum^9_{i=1}{c^i}} = 
\frac{(a^{10^{n-1}})^k}{\sum^9_{i=1}{(a^{10^{n-1}})^i}} = 
\frac{a^{k \cdot 10^{n-1}}}{\sum^9_{i=1}{a^{i \cdot 10^{n-1}}}};$$

Finaly,
\begin{equation}
P_k = \frac{a^{k \cdot 10^{n-1}}}{\sum^9_{i=1}{a^{i \cdot 10^{n-1}}}}\label{log_formula}
\end{equation}

Let us analyze the result. The numerator of the fraction is one of the sum’s components, located in the fraction’s denominator. Thus, the fraction’s value is always less than $1$. Let us check if the sum of all the $P_k$ values equals to $1$.

$$\sum^9_{k=1}P_k = 
\sum^9_{k=1}\frac{a^{k \cdot 10^{n-1}}}{\sum^9_{i=1}{a^{i \cdot 10^{n-1}}}} =
\frac{\sum^9_{k=1}a^{k \cdot 10^{n-1}}}{\sum^9_{i=1}{a^{i \cdot 10^{n-1}}}} = 1;$$

Since the top and the bottom sums are equal to each other, the fraction equals to $1$.

Let us look at the Formula~\vref{log_formula}. It contains the exponent $n$, which is the reference to the function’s range. It means that the function’s $fdp$ depends on the range used for calculations. In addition, the formula contains the base $a$. Thus, we can conclude that the $fdp$ of the logarithmic function depends on its base.

Let us start from figuring out what is going on in the range $[1, 10)$ first. It means that we will assign $n = 1$ for now. In addition, since the variable $a$ is a part of the formula, we will first assign $a = 2$, and then $a = 3$. 
\\

If $n = 1$, then

$$P_k = \frac{a^{k \cdot 10^{1-1}}}{\sum^9_{i=1}{a^{i \cdot 10^{1-1}}}} =
\frac{a^{k \cdot 10^0}}{\sum^9_{i=1}{a^{i \cdot 10^0}}} =
\frac{a^k}{\sum^9_{i=1}{a^i}};$$

If $a = 2$, then

$$P_k = \frac{a^k}{\sum^9_{i=1}{a^i}} = 
\frac{2^k}{\sum^9_{i=1}{2^i}}  = \frac{2^k}{1022}  = \frac{2^{k-1}}{511};$$

If $a = 3$, then

$$P_k = \frac{a^k}{\sum^9_{i=1}{a^i}} = 
\frac{3^k}{\sum^9_{i=1}{3^i}}  = \frac{3^k}{29523}  = \frac{3^{k-1}}{9841};$$
\\ \\
Let us look at a table with $P_k$ values for the functions $y = \log_2(⁡x)$ and $y = \log_3⁡(x)$ on the range $[1, 10)$. Table~\vref{P_k_Log_1_10} shows the data.

\begin{table}[hbt]
\caption{$P_k$ values for $y = \log_2⁡(x)$ and $y = \log_3(⁡x)$ on the range $[1, 10)$}
\label{P_k_Log_1_10}
\centering
\begin{tabular}{ccc}
\toprule
Digit & $P_k$ for $y = \log_2⁡(x)$ & $P_k$ for $y = \log_3⁡(x)$ \\
\midrule
1	&	0.00195695	&	0.00010162	\\
2	&	0.00391389	&	0.00030485	\\
3	&	0.00782779	&	0.00091454	\\
4	&	0.01565558	&	0.00274362	\\
5	&	0.03131115	&	0.00823087	\\
6	&	0.06262231	&	0.02469261	\\
7	&	0.12524462	&	0.07407784	\\
8	&	0.25048924	&	0.22223351	\\
9	&	0.50097847	&	0.66670054	\\
\bottomrule
Sum & 1.00000000 & 1.00000000  \\
\bottomrule
\end{tabular}
\end{table}

While looking at the $P_k$ formulas and graphs above, we could already see that bigger digits have higher probabilities. This conclusion is proven again by our table for $n = 1$. In addition, it looks like it will be true for higher $n$-values as well. 

As we can see, each $P_k$ value is bigger than the previous one exactly $a$ times. Moreover, the difference between first and last $fdp$ values for both $a$ is much higher than in all previous cases of this work. Furthermore, it looks like this difference is even higher for functions with higher bases.

\begin{figure}[h]
\centering 
\includegraphics[width=1\columnwidth]{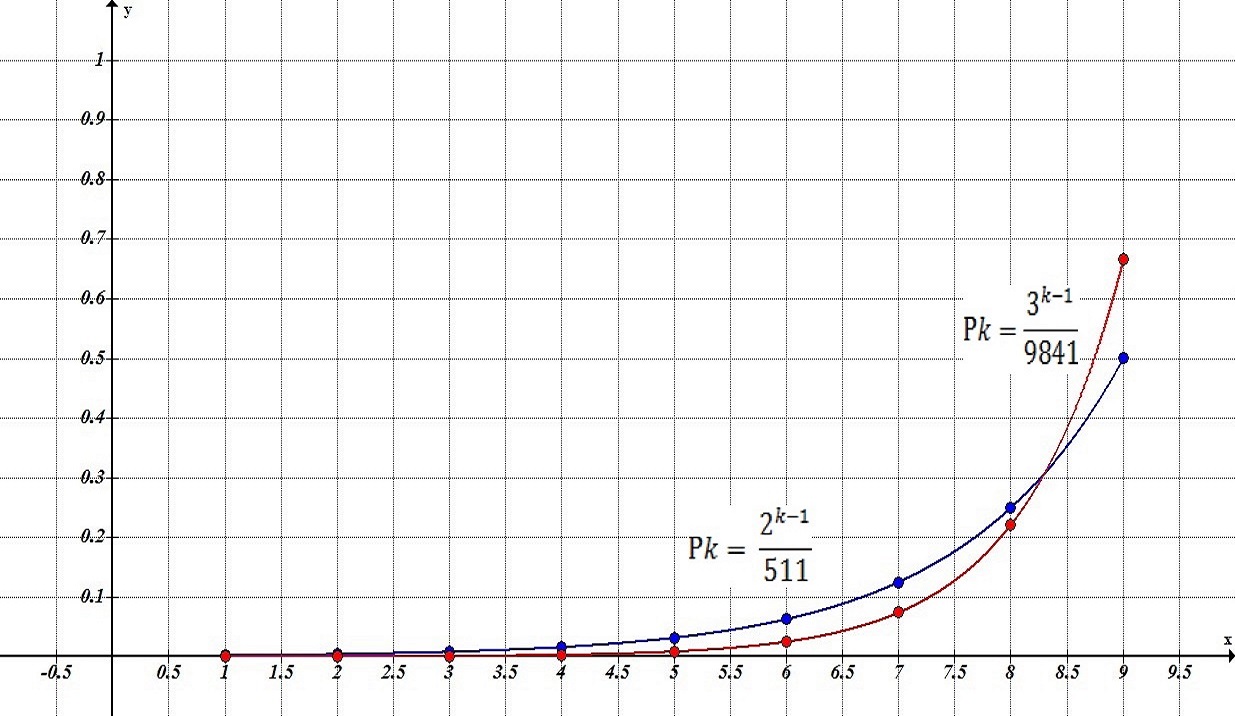} 
\caption{The functions $P_k=  \frac{2^{k-1}}{511}$ and $P_k=\frac{3^{k-1}}{9841}$}
\label{Log_Prob_1_10}
\end{figure}

After the formulas above were defined, let us see graphs of the functions for $a = 2$ and $a = 3$ on the range $[1, 10)$. The functions are $P_k=  \frac{2^{k-1}}{511}$ and $P_k=\frac{3^{k-1}}{9841}$. Both graphs are shown in Figure~\vref{Log_Prob_1_10}.
\\

Let us see now how the formulas will look like on the range $[10, 100)$. Thus, we will assign $n = 2$. Then,

$$P_k = \frac{a^{k \cdot 10^{2-1}}}{\sum^9_{i=1}{a^{i \cdot 10^{2-1}}}} =
\frac{a^{k \cdot 10^1}}{\sum^9_{i=1}{a^{i \cdot 10^1}}} =
\frac{a^{10k}}{\sum^9_{i=1}{a^{10i}}};$$

If $a = 2$, then

$$P_k = \frac{a^{10k}}{\sum^9_{i=1}{a^{10i}}} = 
\frac{2^{10k}}{\sum^9_{i=1}{2^{10i}}} =
\frac{(2^{10})^{k}}{\sum^9_{i=1}{(2^{10})^i}} =
\frac{1024^k}{\sum^9_{i=1}{1024^i}};$$

If $a = 3$, then

$$P_k = \frac{a^{10k}}{\sum^9_{i=1}{a^{10i}}} = 
\frac{3^{10k}}{\sum^9_{i=1}{3^{10i}}} =
\frac{(3^{10})^{k}}{\sum^9_{i=1}{(3^{10})^i}} =
\frac{59049^k}{\sum^9_{i=1}{59049^i}};$$
\\ \\
We should now look at a table with $P_k$ values for the functions $y = \log_2⁡(x)$ and 
$y = \log_3⁡(x)$ on the range $[10, 100)$. Table~\vref{P_k_Log_10_100} shows the data. As we can see, both $P_9$ values almost equal to $1$. The rest of $P_k$ values are very small. 

\begin{table}[hbt]
\caption{$P_k$ values for $y = \log_2⁡(x)$ and $y = \log_3(⁡x)$ on the range $[10, 100)$}
\label{P_k_Log_10_100}
\centering
\begin{tabular}{cll}
\toprule
Digit & $P_k$ for $y = \log_2⁡(x)$ & $P_k$ for $y = \log_3⁡(x)$ \\
\midrule
1	&	8.263728 $\cdot$ 10$^{-25}$	&	6.765381 $\cdot$ 10$^{-39}$	\\
2	&	8.462058 $\cdot$ 10$^{-22}$	&	3.994890 $\cdot$ 10$^{-34}$	\\
3	&	8.665147 $\cdot$ 10$^{-19}$	&	2.358943 $\cdot$ 10$^{-29}$	\\
4	&	8.873111 $\cdot$ 10$^{-16}$	&	1.392932 $\cdot$ 10$^{-24}$	\\
5	&	9.086065 $\cdot$ 10$^{-13}$	&	8.225124 $\cdot$ 10$^{-20}$	\\
6	&	9.304131 $\cdot$ 10$^{-10}$	&	4.856853 $\cdot$ 10$^{-15}$	\\
7	&	9.527430 $\cdot$ 10$^{-7}$		&	2.867923 $\cdot$ 10$^{-10}$	\\
8	&	9.756088 $\cdot$ 10$^{-4}$		&	1.693480 $\cdot$ 10$^{-5}$		\\
9	&	9.990234 $\cdot$ 10$^{-1}$		&	9.999831 $\cdot$ 10$^{-1}$ \\
\bottomrule
Sum & 1.00000000 & 1.00000000  \\
\bottomrule
\end{tabular}
\end{table}

\begin{figure}[h]
\centering 
\includegraphics[width=1\columnwidth]{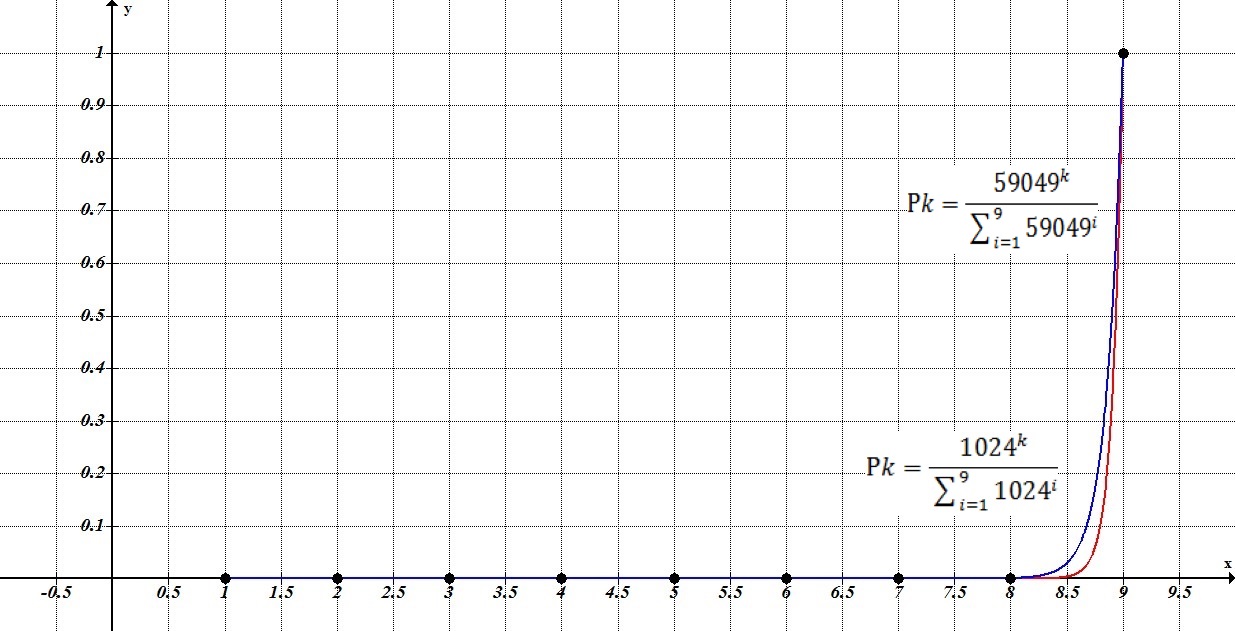} 
\caption{The functions $P_k = \frac{1024^k}{\sum^9_{i=1}{1024^i}}$ and $P_k = \frac{59049^k}{\sum^9_{i=1}{59049^i}}$}
\label{Log_Prob_10_100}
\end{figure}

Let us also compare the functions $P_k = \frac{1024^k}{\sum^9_{i=1}{1024^i}}$ and $P_k = \frac{59049^k}{\sum^9_{i=1}{59049^i}}$ while looking at their graphs. Both graphs corresponding to the range $[10, 100)$ are presented in Figure~\vref{Log_Prob_10_100}. As we can see, they are located really close to each other. However, the curve connecting the points of the $P_k$ values of  $y = \log_3⁡(x)$ function is located at the right of the curve connecting the points of the $P_k$ values of  $y = \log_2⁡(x)$ function in between $k$ values of $8$ and $9$. It is the only way to differentiate them in the graph. As it was mentioned above, all the $P_k$ values except $P_9$ of both functions are so small, that they almost equal to zero.

Let us think about a proper example that can be used in this section in such a way that it describes a process that has logarithmic growth and has $P_k$ values equal or approximately equal to ours. 

First, let us analyze what will happen if we stretch or squeeze our initial $y = \log_a⁡(x)$ function. In particular, we need to know if it affects its $P_k$ values. Only if the $P_k$ values stay unchanged, we can try to find a proper example. To make our steps a little bit easier later, let us divide the whole function by a constant $h$ instead of multiplying by it as we did before. Depending on value of $h$, we either stretch or squeeze our function.  

Let $y = \frac{\log_a⁡(x)}{h}$. Then, $hy = \log_a⁡(x)$. Let us switch $x$ and $y$: 

$$hx = \log_a(y); \: y = a^{hx}; \: f(x)^{-1} = a^{hx}$$

$$P_k = \frac{a^{h(k+1)\cdot10^{n-1}} - a^{h\cdot k\cdot10^{n-1}}}{a^{h\cdot10\cdot10^{n-1}} - a^{h\cdot10^{n-1}}} = 
 \frac{(a^h)^{(k+1)\cdot10^{n-1}} - (a^h)^{k\cdot10^{n-1}}}{(a^h)^{10\cdot10^{n-1}} - (a^h)^{10^{n-1}}};$$
 
Let us use a substitution $w = a^h$. Since $a$ and $h$ are constants, $w$ is a constant as well. Then:

$$P_k = \frac{w^{(k+1)\cdot10^{n-1}} - w^{k\cdot10^{n-1}}}{w^{10\cdot10^{n-1}} - w^{10^{n-1}}};$$
\\ \\
As we can see, we came to the very first step of getting the $P_k$ values for a logarithmic function. The only difference is that the formula has a base $w$ now, not a base $a$ as it was before. We know that logarithmic functions with different bases have different $P_k$ values. In other words, stretching or squeezing of the logarithmic function will change its $P_k$ values, which tells us that we will have to change the base of our $\log$ function in a certain way if $h \neq 1$. In particular, the base $a$ will be changed for the base $a^h$.

Next, let us check what will happen if we multiply or divide the variable $x$ by a constant directly. Let us have $y = \log_a⁡(\frac{x}{m})$, where $m$ is a constant. Then we have:

$$y = \log_a⁡(\frac{x}{m}) = \log_a⁡(x) - \log_a⁡(m);$$
\\
Since both $a$ and $m$ are constants, $\log_a(⁡m)$ is a constant as well. As we can see, dividing or multiplying the variable $x$ by a constant will result in a vertical shift of the initial function. Let us assign $q = \log_a(⁡m)$ and check how a vertical shift effects the $P_k$ values. 

Let $y= \log_a⁡(x) - q$. Then $y + q = \log_a⁡(x)$. After we switch $x$ and $y$, we get 

$$\log_a(⁡y) = x + q; \: y= a^{x+q}; \: f(x)^{-1} = a^{x+q}$$

$$P_k = \frac{a^{(k+1)\cdot10^{n-1}+q} - a^{k\cdot10^{n-1}+q}}{a^{10\cdot10^{n-1}+q} - a^{10^{n-1}+q}} =
\frac{a^{(k+1)\cdot10^{n-1}}\cdot a^q - a^{k\cdot10^{n-1}}\cdot a^q}{a^{10\cdot10^{n-1}}\cdot a^q - a^{10^{n-1}}\cdot a^q}$$
$$= \frac{a^q(a^{(k+1)\cdot10^{n-1}} - a^{k\cdot10^{n-1}})}{a^q(a^{10\cdot10^{n-1}} - a^{10^{n-1}})} = 
\frac{a^{(k+1)\cdot10^{n-1}} - a^{k\cdot10^{n-1}}}{a^{10\cdot10^{n-1}} - a^{10^{n-1}}};$$

Obviously, we came back to the first step of getting the $P_k$ values for our initial function.  The last fraction above is identical to the one with no vertical shift. It means that a vertical shift of a logarithmic function does not change the function’s initial $P_k$ values. Thus, we do not have to do anything about multiplying or dividing the variable $x$ by a constant if we face a case like that.

It looks like we have enough information by now to start thinking about our example. We had a chance to talk about bacterial growth before. Let us talk about population growth now. 

\subsubsection{A real-life example}

We will be using the population growth formula 

\begin{equation}
P = 1600e^{0.025T}
\label{Population_Growth}
\end{equation}

\noindent
where $P$ is the number of people at time $T$, $1600$ is the number of people when $T = 0$, $0.025$ is the population growth rate, and $T$ is time in years. Let us isolate $T$. It will allow us to get a logarithmic function.

$$\frac{P}{1600} = e^{0.025T}; \: 0.025T = \ln(\frac{P}{1600}); \: \frac{1}{40}T = \ln(\frac{P}{1600}); \: T = 40\ln(\frac{P}{1600});$$

Thus, we will be using the following function in our example.

\begin{equation}
T(P) = 40\ln(\frac{P}{1600})
\label{Population_Growth_Example}
\end{equation}

This equation allows us to calculate how many years are needed to get a certain population. Let us recall that a $\log$ function changes its $P_k$ values depending on its range. To make our work meaningful, let us choose the range $[1, 10)$. All other ranges are much harder to analyze. Let us calculate the initial and final values of $P$ in such a way that they will be starting and ending points of our domain. 

First, let $T = 1$. Then,

$$1 = 40\ln(\frac{P}{1600}); \: \frac{1}{40} = \ln(\frac{P}{1600}); \: \ln(\frac{P}{1600}) = 0.025; \: \frac{P}{1600} = e^{0.025};$$

$$P = 1600e^{0.025} \approx 1600 \cdot 1.0253 = 1640.48;$$

We will approximate the initial $P$ value to $1640$, so it would bring us to the initial $T$ value that is a little bit less than $1$.

Next, let $T = 10$. Then,

$$10 = 40\ln(\frac{P}{1600}); \: \frac{1}{4} = \ln(\frac{P}{1600}); \: \ln(\frac{P}{1600}) = 0.25; \: \frac{P}{1600} = e^{0.25};$$

$$P = 1600e^{0.25} \approx 1600 \cdot 1.2840 = 2054.4;$$

We will approximate the final $P$ value to $2055$, so it would give us the final $T$ value that is a little bit greater than $10$. The difference $2055 - 1640 = 415$ will give us enough steps for our calculations. 

We will calculate a new $T$ value as soon as the number of people gets bigger by one. If calculated $T$ values exceed our range $[1, 10)$, we will not use them for our $fdp$ calculations. Table~\vref{T(P)} shows the results.

\newpage
\begin{longtable}[c]{lllllllll} 
 \caption{Calculated $T(P)$ values and their first digits}
 \label{T(P)}\\
\toprule
$P$- & $T$- & First & $P$- & $T$- & First & $P$- & $T$- & First \\
values & values & digit & values & values & digit & values & values & digit \\ 
\midrule
\endfirsthead 

\toprule
$P$- & $T$- & First & $P$- & $T$- & First & $P$- & $T$- & First \\
values & values & digit & values & values & digit & values & values & digit \\ 
\midrule
\endhead 
1640	&	0.9877	&		&	1690	&	2.1890	&	2	&	1740	&	3.3553	&	3	\\
1641	&	1.0121	&	1	&	1691	&	2.2127	&	2	&	1741	&	3.3782	&	3	\\
1642	&	1.0365	&	1	&	1692	&	2.2363	&	2	&	1742	&	3.4012	&	3	\\
1643	&	1.0608	&	1	&	1693	&	2.2599	&	2	&	1743	&	3.4242	&	3	\\
1644	&	1.0851	&	1	&	1694	&	2.2836	&	2	&	1744	&	3.4471	&	3	\\
1645	&	1.1095	&	1	&	1695	&	2.3072	&	2	&	1745	&	3.4700	&	3	\\
1646	&	1.1338	&	1	&	1696	&	2.3308	&	2	&	1746	&	3.4930	&	3	\\
1647	&	1.1581	&	1	&	1697	&	2.3543	&	2	&	1747	&	3.5159	&	3	\\
1648	&	1.1824	&	1	&	1698	&	2.3779	&	2	&	1748	&	3.5387	&	3	\\
1649	&	1.2066	&	1	&	1699	&	2.4014	&	2	&	1749	&	3.5616	&	3	\\
1650	&	1.2309	&	1	&	1700	&	2.4250	&	2	&	1750	&	3.5845	&	3	\\
1651	&	1.2551	&	1	&	1701	&	2.4485	&	2	&	1751	&	3.6073	&	3	\\
1652	&	1.2793	&	1	&	1702	&	2.4720	&	2	&	1752	&	3.6302	&	3	\\
1653	&	1.3035	&	1	&	1703	&	2.4955	&	2	&	1753	&	3.6530	&	3	\\
1654	&	1.3277	&	1	&	1704	&	2.5190	&	2	&	1754	&	3.6758	&	3	\\
1655	&	1.3519	&	1	&	1705	&	2.5425	&	2	&	1755	&	3.6986	&	3	\\
1656	&	1.3761	&	1	&	1706	&	2.5659	&	2	&	1756	&	3.7214	&	3	\\
1657	&	1.4002	&	1	&	1707	&	2.5894	&	2	&	1757	&	3.7442	&	3	\\
1658	&	1.4243	&	1	&	1708	&	2.6128	&	2	&	1758	&	3.7669	&	3	\\
1659	&	1.4485	&	1	&	1709	&	2.6362	&	2	&	1759	&	3.7897	&	3	\\
1660	&	1.4726	&	1	&	1710	&	2.6596	&	2	&	1760	&	3.8124	&	3	\\
1661	&	1.4966	&	1	&	1711	&	2.6830	&	2	&	1761	&	3.8351	&	3	\\
1662	&	1.5207	&	1	&	1712	&	2.7063	&	2	&	1762	&	3.8578	&	3	\\
1663	&	1.5448	&	1	&	1713	&	2.7297	&	2	&	1763	&	3.8805	&	3	\\
1664	&	1.5688	&	1	&	1714	&	2.7530	&	2	&	1764	&	3.9032	&	3	\\
1665	&	1.5929	&	1	&	1715	&	2.7764	&	2	&	1765	&	3.9259	&	3	\\
1666	&	1.6169	&	1	&	1716	&	2.7997	&	2	&	1766	&	3.9485	&	3	\\
1667	&	1.6409	&	1	&	1717	&	2.8230	&	2	&	1767	&	3.9712	&	3	\\
1668	&	1.6649	&	1	&	1718	&	2.8463	&	2	&	1768	&	3.9938	&	3	\\
1669	&	1.6888	&	1	&	1719	&	2.8696	&	2	&	1769	&	4.0164	&	4	\\
1670	&	1.7128	&	1	&	1720	&	2.8928	&	2	&	1770	&	4.0390	&	4	\\
1671	&	1.7367	&	1	&	1721	&	2.9161	&	2	&	1771	&	4.0616	&	4	\\
1672	&	1.7607	&	1	&	1722	&	2.9393	&	2	&	1772	&	4.0842	&	4	\\
1673	&	1.7846	&	1	&	1723	&	2.9625	&	2	&	1773	&	4.1068	&	4	\\
1674	&	1.8085	&	1	&	1724	&	2.9857	&	2	&	1774	&	4.1293	&	4	\\
1675	&	1.8324	&	1	&	1725	&	3.0089	&	3	&	1775	&	4.1519	&	4	\\
1676	&	1.8563	&	1	&	1726	&	3.0321	&	3	&	1776	&	4.1744	&	4	\\
1677	&	1.8801	&	1	&	1727	&	3.0553	&	3	&	1777	&	4.1969	&	4	\\
1678	&	1.9040	&	1	&	1728	&	3.0784	&	3	&	1778	&	4.2194	&	4	\\
1679	&	1.9278	&	1	&	1729	&	3.1016	&	3	&	1779	&	4.2419	&	4	\\
1680	&	1.9516	&	1	&	1730	&	3.1247	&	3	&	1780	&	4.2644	&	4	\\
1681	&	1.9754	&	1	&	1731	&	3.1478	&	3	&	1781	&	4.2869	&	4	\\
1682	&	1.9992	&	1	&	1732	&	3.1709	&	3	&	1782	&	4.3093	&	4	\\
1683	&	2.0230	&	2	&	1733	&	3.1940	&	3	&	1783	&	4.3317	&	4	\\
1684	&	2.0467	&	2	&	1734	&	3.2171	&	3	&	1784	&	4.3542	&	4	\\
1685	&	2.0705	&	2	&	1735	&	3.2402	&	3	&	1785	&	4.3766	&	4	\\
1686	&	2.0942	&	2	&	1736	&	3.2632	&	3	&	1786	&	4.3990	&	4	\\
1687	&	2.1179	&	2	&	1737	&	3.2862	&	3	&	1787	&	4.4214	&	4	\\
1688	&	2.1416	&	2	&	1738	&	3.3093	&	3	&	1788	&	4.4438	&	4	\\
1689	&	2.1653	&	2	&	1739	&	3.3323	&	3	&	1789	&	4.4661	&	4	\\ \\ \\
1790	&	4.4885	&	4	&	1840	&	5.5905	&	5	&	1890	&	6.6629	&	6	\\
1791	&	4.5108	&	4	&	1841	&	5.6122	&	5	&	1891	&	6.6841	&	6	\\
1792	&	4.5331	&	4	&	1842	&	5.6339	&	5	&	1892	&	6.7052	&	6	\\
1793	&	4.5555	&	4	&	1843	&	5.6556	&	5	&	1893	&	6.7264	&	6	\\
1794	&	4.5778	&	4	&	1844	&	5.6773	&	5	&	1894	&	6.7475	&	6	\\
1795	&	4.6001	&	4	&	1845	&	5.6990	&	5	&	1895	&	6.7686	&	6	\\
1796	&	4.6223	&	4	&	1846	&	5.7207	&	5	&	1896	&	6.7897	&	6	\\
1797	&	4.6446	&	4	&	1847	&	5.7424	&	5	&	1897	&	6.8108	&	6	\\
1798	&	4.6669	&	4	&	1848	&	5.7640	&	5	&	1898	&	6.8319	&	6	\\
1799	&	4.6891	&	4	&	1849	&	5.7857	&	5	&	1899	&	6.8530	&	6	\\
1800	&	4.7113	&	4	&	1850	&	5.8073	&	5	&	1900	&	6.8740	&	6	\\
1801	&	4.7335	&	4	&	1851	&	5.8289	&	5	&	1901	&	6.8951	&	6	\\
1802	&	4.7557	&	4	&	1852	&	5.8505	&	5	&	1902	&	6.9161	&	6	\\
1803	&	4.7779	&	4	&	1853	&	5.8721	&	5	&	1903	&	6.9371	&	6	\\
1804	&	4.8001	&	4	&	1854	&	5.8937	&	5	&	1904	&	6.9581	&	6	\\
1805	&	4.8223	&	4	&	1855	&	5.9152	&	5	&	1905	&	6.9791	&	6	\\
1806	&	4.8444	&	4	&	1856	&	5.9368	&	5	&	1906	&	7.0001	&	7	\\
1807	&	4.8666	&	4	&	1857	&	5.9583	&	5	&	1907	&	7.0211	&	7	\\
1808	&	4.8887	&	4	&	1858	&	5.9799	&	5	&	1908	&	7.0421	&	7	\\
1809	&	4.9108	&	4	&	1859	&	6.0014	&	6	&	1909	&	7.0630	&	7	\\
1810	&	4.9329	&	4	&	1860	&	6.0229	&	6	&	1910	&	7.0840	&	7	\\
1811	&	4.9550	&	4	&	1861	&	6.0444	&	6	&	1911	&	7.1049	&	7	\\
1812	&	4.9771	&	4	&	1862	&	6.0659	&	6	&	1912	&	7.1258	&	7	\\
1813	&	4.9992	&	4	&	1863	&	6.0874	&	6	&	1913	&	7.1468	&	7	\\
1814	&	5.0212	&	5	&	1864	&	6.1088	&	6	&	1914	&	7.1677	&	7	\\
1815	&	5.0433	&	5	&	1865	&	6.1303	&	6	&	1915	&	7.1886	&	7	\\
1816	&	5.0653	&	5	&	1866	&	6.1517	&	6	&	1916	&	7.2094	&	7	\\
1817	&	5.0873	&	5	&	1867	&	6.1732	&	6	&	1917	&	7.2303	&	7	\\
1818	&	5.1093	&	5	&	1868	&	6.1946	&	6	&	1918	&	7.2512	&	7	\\
1819	&	5.1313	&	5	&	1869	&	6.2160	&	6	&	1919	&	7.2720	&	7	\\
1820	&	5.1533	&	5	&	1870	&	6.2374	&	6	&	1920	&	7.2929	&	7	\\
1821	&	5.1753	&	5	&	1871	&	6.2588	&	6	&	1921	&	7.3137	&	7	\\
1822	&	5.1972	&	5	&	1872	&	6.2801	&	6	&	1922	&	7.3345	&	7	\\
1823	&	5.2192	&	5	&	1873	&	6.3015	&	6	&	1923	&	7.3553	&	7	\\
1824	&	5.2411	&	5	&	1874	&	6.3229	&	6	&	1924	&	7.3761	&	7	\\
1825	&	5.2631	&	5	&	1875	&	6.3442	&	6	&	1925	&	7.3969	&	7	\\
1826	&	5.2850	&	5	&	1876	&	6.3655	&	6	&	1926	&	7.4177	&	7	\\
1827	&	5.3069	&	5	&	1877	&	6.3868	&	6	&	1927	&	7.4384	&	7	\\
1828	&	5.3288	&	5	&	1878	&	6.4082	&	6	&	1928	&	7.4592	&	7	\\
1829	&	5.3506	&	5	&	1879	&	6.4294	&	6	&	1929	&	7.4799	&	7	\\
1830	&	5.3725	&	5	&	1880	&	6.4507	&	6	&	1930	&	7.5007	&	7	\\
1831	&	5.3943	&	5	&	1881	&	6.4720	&	6	&	1931	&	7.5214	&	7	\\
1832	&	5.4162	&	5	&	1882	&	6.4933	&	6	&	1932	&	7.5421	&	7	\\
1833	&	5.4380	&	5	&	1883	&	6.5145	&	6	&	1933	&	7.5628	&	7	\\
1834	&	5.4598	&	5	&	1884	&	6.5357	&	6	&	1934	&	7.5835	&	7	\\
1835	&	5.4816	&	5	&	1885	&	6.5570	&	6	&	1935	&	7.6041	&	7	\\
1836	&	5.5034	&	5	&	1886	&	6.5782	&	6	&	1936	&	7.6248	&	7	\\
1837	&	5.5252	&	5	&	1887	&	6.5994	&	6	&	1937	&	7.6455	&	7	\\
1838	&	5.5470	&	5	&	1888	&	6.6206	&	6	&	1938	&	7.6661	&	7	\\
1839	&	5.5687	&	5	&	1889	&	6.6418	&	6	&	1939	&	7.6867	&	7	\\ \\ \\ \\
1940	&	7.7074	&	7	&	1979	&	8.5035	&	8	&	2018	&	9.2841	&	9	\\
1941	&	7.7280	&	7	&	1980	&	8.5237	&	8	&	2019	&	9.3039	&	9	\\
1942	&	7.7486	&	7	&	1981	&	8.5439	&	8	&	2020	&	9.3238	&	9	\\
1943	&	7.7692	&	7	&	1982	&	8.5641	&	8	&	2021	&	9.3436	&	9	\\
1944	&	7.7898	&	7	&	1983	&	8.5843	&	8	&	2022	&	9.3633	&	9	\\
1945	&	7.8103	&	7	&	1984	&	8.6045	&	8	&	2023	&	9.3831	&	9	\\
1946	&	7.8309	&	7	&	1985	&	8.6246	&	8	&	2024	&	9.4029	&	9	\\
1947	&	7.8514	&	7	&	1986	&	8.6448	&	8	&	2025	&	9.4226	&	9	\\
1948	&	7.8720	&	7	&	1987	&	8.6649	&	8	&	2026	&	9.4424	&	9	\\
1949	&	7.8925	&	7	&	1988	&	8.6850	&	8	&	2027	&	9.4621	&	9	\\
1950	&	7.9130	&	7	&	1989	&	8.7051	&	8	&	2028	&	9.4819	&	9	\\
1951	&	7.9335	&	7	&	1990	&	8.7252	&	8	&	2029	&	9.5016	&	9	\\
1952	&	7.9540	&	7	&	1991	&	8.7453	&	8	&	2030	&	9.5213	&	9	\\
1953	&	7.9745	&	7	&	1992	&	8.7654	&	8	&	2031	&	9.5410	&	9	\\
1954	&	7.9950	&	7	&	1993	&	8.7855	&	8	&	2032	&	9.5607	&	9	\\
1955	&	8.0155	&	8	&	1994	&	8.8056	&	8	&	2033	&	9.5804	&	9	\\
1956	&	8.0359	&	8	&	1995	&	8.8256	&	8	&	2034	&	9.6000	&	9	\\
1957	&	8.0564	&	8	&	1996	&	8.8457	&	8	&	2035	&	9.6197	&	9	\\
1958	&	8.0768	&	8	&	1997	&	8.8657	&	8	&	2036	&	9.6393	&	9	\\
1959	&	8.0972	&	8	&	1998	&	8.8857	&	8	&	2037	&	9.6590	&	9	\\
1960	&	8.1176	&	8	&	1999	&	8.9057	&	8	&	2038	&	9.6786	&	9	\\
1961	&	8.1380	&	8	&	2000	&	8.9257	&	8	&	2039	&	9.6982	&	9	\\
1962	&	8.1584	&	8	&	2001	&	8.9457	&	8	&	2040	&	9.7178	&	9	\\
1963	&	8.1788	&	8	&	2002	&	8.9657	&	8	&	2041	&	9.7375	&	9	\\
1964	&	8.1992	&	8	&	2003	&	8.9857	&	8	&	2042	&	9.7570	&	9	\\
1965	&	8.2195	&	8	&	2004	&	9.0057	&	9	&	2043	&	9.7766	&	9	\\
1966	&	8.2399	&	8	&	2005	&	9.0256	&	9	&	2044	&	9.7962	&	9	\\
1967	&	8.2602	&	8	&	2006	&	9.0456	&	9	&	2045	&	9.8158	&	9	\\
1968	&	8.2806	&	8	&	2007	&	9.0655	&	9	&	2046	&	9.8353	&	9	\\
1969	&	8.3009	&	8	&	2008	&	9.0854	&	9	&	2047	&	9.8549	&	9	\\
1970	&	8.3212	&	8	&	2009	&	9.1053	&	9	&	2048	&	9.8744	&	9	\\
1971	&	8.3415	&	8	&	2010	&	9.1252	&	9	&	2049	&	9.8939	&	9	\\
1972	&	8.3618	&	8	&	2011	&	9.1451	&	9	&	2050	&	9.9134	&	9	\\
1973	&	8.3821	&	8	&	2012	&	9.1650	&	9	&	2051	&	9.9330	&	9	\\
1974	&	8.4023	&	8	&	2013	&	9.1849	&	9	&	2052	&	9.9525	&	9	\\
1975	&	8.4226	&	8	&	2014	&	9.2048	&	9	&	2053	&	9.9719	&	9	\\
1976	&	8.4428	&	8	&	2015	&	9.2246	&	9	&	2054	&	9.9914	&	9	\\
1977	&	8.4631	&	8	&	2016	&	9.2445	&	9	&	2055	&	10.0109	&		\\
1978	&	8.4833	&	8	&	2017	&	9.2643	&	9	&			&			&		\\	
\bottomrule \\
\end{longtable}

Before we summarize the results from Table~\vref{T(P)}, let us prepare the data that will be compared to the mentioned above results. We already discovered that the $P_k$ values of a function $y = \frac{1}{h}\log_a⁡(\frac{x}{m})$ are the same as the $P_k$ values of a function $y = \log_{(a^h)}⁡(x)$ and the constant $m$ has no effect on it.

In the example above, $a = e$ because we are dealing with natural $\log$ there and $h = 0.025$. Thus, $a^h = e^{0.025} \approx 1.0253$. It means that we will be calculating the $P_k$ values of the function $y = \log_{1.0253}(x)$. Moreover, we are expecting that the values will be at least approximately equal to those from Table~\vref{T(P)}.

\begin{table}[hbt]
\caption{Analysis of $T(p)$ values first digits}
\label{T(p)_analysis}
\centering
\begin{tabular}{cccc}
\toprule
Digit & Count & Frequency       & $P_k$ values \\  
       &           & (Count / Sum) &  for $y = \log_{1.0253}⁡(x)$ \\ 
\midrule
1	&	42	&	0.10144928	&	0.10033451	\\
2	&	42	&	0.10144928	&	0.10287297	\\
3	&	44	&	0.10628019	&	0.10547566	\\
4	&	45	&	0.10869565	&	0.10814419	\\
5	&	45	&	0.10869565	&	0.11088024	\\
6	&	47	&	0.11352657	&	0.11368551	\\
7	&	49	&	0.11835749	&	0.11656175	\\
8	&	49	&	0.11835749	&	0.11951077	\\
9	&	51	&	0.12318841	&	0.12253439	\\
\bottomrule
Sum	&	217	&	1.00000000	&	1.00000000	\\
\bottomrule
\end{tabular}
\end{table}

\FloatBarrier
Both groups of the $P_k$ values are presented in Table~\vref{T(p)_analysis}. The middle column shows the result of the calculations based on the formula $y= 40 \ln(\frac{x}{1600})$ and Table~\vref{T(P)}. The last column shows the results based on the general $P_k$ formula that was defined for any $\log$ function and  the formula $y = \log_{1.0253}⁡(x)$. The range $[1, 10)$ was used in both cases.

As we can see, the corresponding $P_k$ values are really close to each other. The last column shows numbers that are more accurate because their calculations were derived directly from the formula. The middle column has numbers that are less accurate because the $\log$ function used in the example was approximated to the discrete one with steps of one unit. Having more steps would definitely help to get numbers that are more accurate.

We also can see that the difference between the first and last $P_k$ values is not as big as in all our previous calculations. We already discovered that a bigger $a$ value causes bigger gaps between $P_1$ and $P_9$. However, the base of the function $\log_{1.0253⁡}(x)$ is very small. It is the reason why the differences between $P_k$ values in Table~\vref{T(p)_analysis} are not that big.

\begin{figure}[h]
\centering 
\includegraphics[width=1\columnwidth]{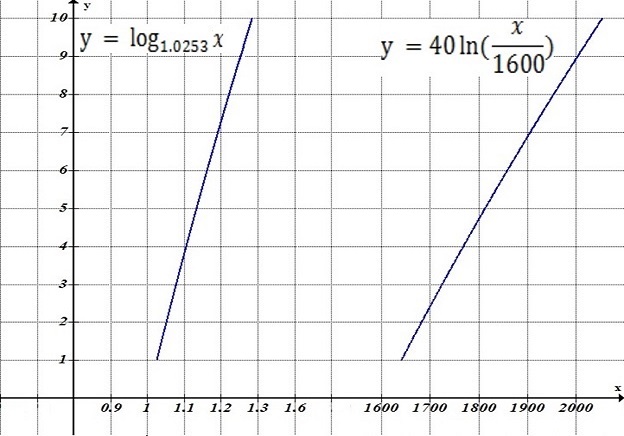} 
\caption{The functions $y = \log_{1.0253}(⁡x)$ and $y = 40\ln(\frac{x}{1600})$}
\label{Log_From_The_Formula_And_Table}
\end{figure}

To finish our comparison, let us graph both functions. Unfortunately, we will not be able to graph them in the same $xy$-plane without doing extra tricks due to a big gap in between their domains that we have to pick. The function $y = 40\ln(\frac{x}{1600})$ is presented above on the domain $[1640, 2055]$. 

Let us calculate the domain of the function $y = \log_{1.0253}(⁡x)$ in such a way that it will match to the range $[1, 10)$. 

$$1.0253^1=1.0253; \: 1.0253^{10} \approx 1.2838;$$  

Thus, the domain is $[1.0253, 1.2838)$. If we skip certain portions of the $x$-axis, we will be able to locate both graphs together. 

We are ready to look at the graphs now. Both of them are shown in Figure~\vref{Log_From_The_Formula_And_Table}. Please, notice that their actual domains are really far from each other. 
As we can see, both $\log$ functions are stretched so much that their graphed portions look almost like straight lines.  That is why their $P_k$ values are close to those of straight lines.

\newpage
\subsection{$P_k$ values of the reciprocal function $y = \frac{a}{x}$}
\subsubsection{Deriving the $P_k$ formula}

The reciprocal function $y = \frac{a}{x}$, where $a$ is a constant, is the last one in this chapter. Since we are graphing all our functions in the first quadrant only, the constant $a$ for this function will be always positive. 

We will consider three cases: $a$ > $1$, $a = 1$, and $0$ < $a$ < $1$. We also should keep in mind that our range is always restricted to $[1, \infty)$, that is why our domain for this function is the interval $(0, a]$.

All three functions are graphed on the range $[1, 10)$. Let us look at each case separately. 

\begin{figure}[h]
\centering 
\includegraphics[width=1\columnwidth]{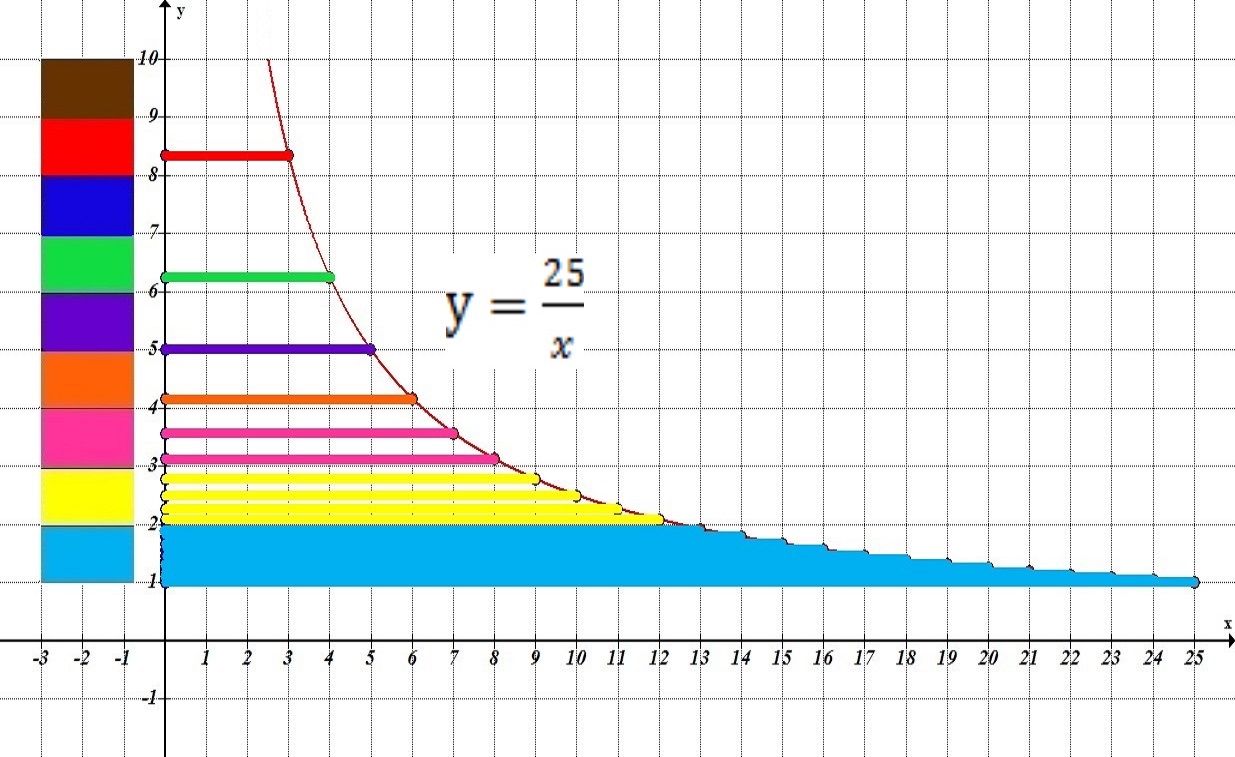} 
\caption{The function $y = \frac{25}{x}$}
\label{Reciprocal_25}
\end{figure}

\begin{figure}[h]
\centering 
\includegraphics[width=1\columnwidth]{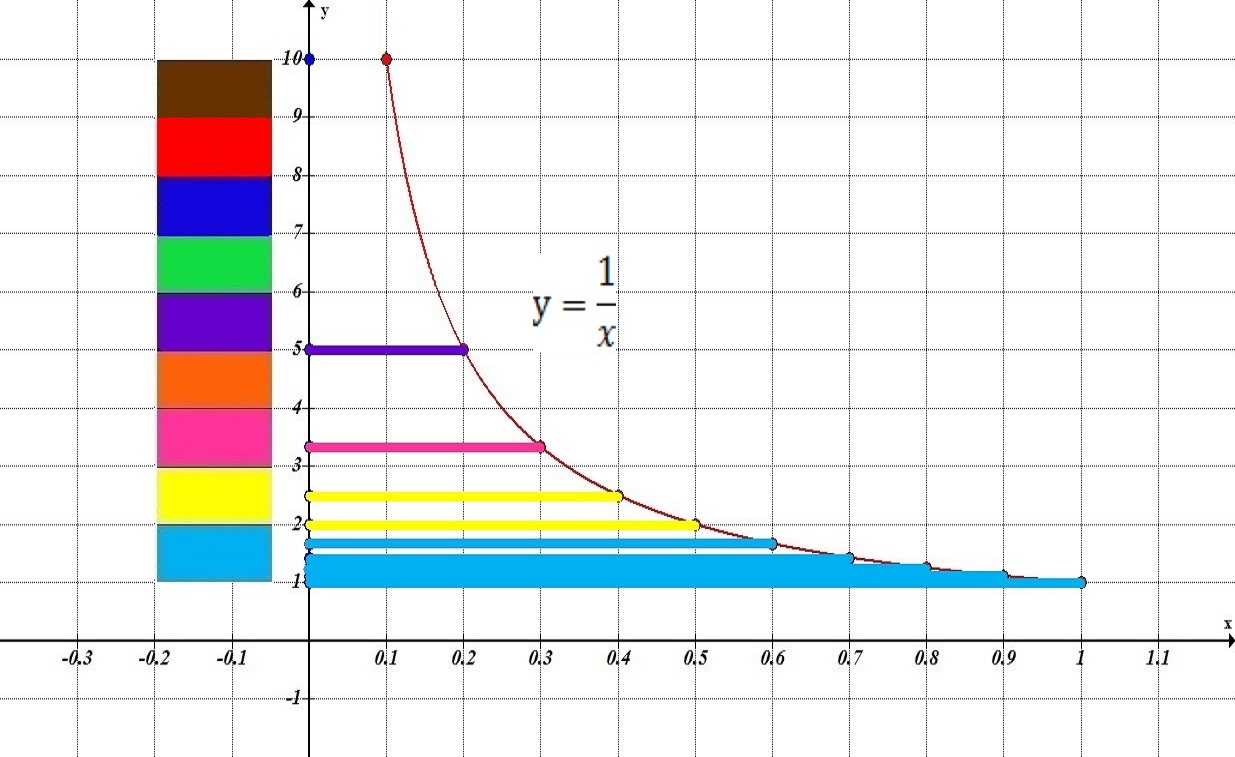} 
\caption{The function $y = \frac{1}{x}$}
\label{Reciprocal_1}
\end{figure}

\begin{figure}[h]
\centering 
\includegraphics[width=1\columnwidth]{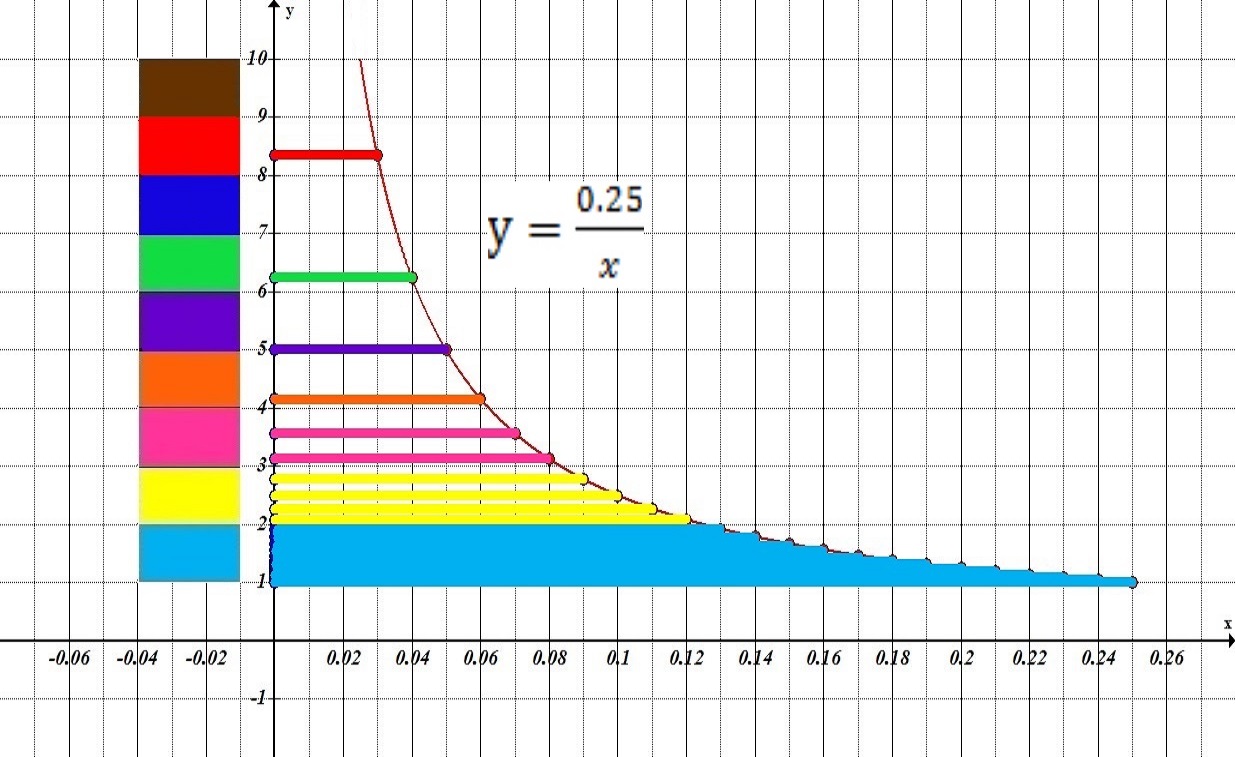} 
\caption{The function $y = \frac{0.25}{x}$}
\label{Reciprocal_Quarter}
\end{figure}

The function $y = \frac{25}{x}$ is shown in Figure~\vref{Reciprocal_25}, the function $y = \frac{1}{x}$ is shown in Figure~\vref{Reciprocal_1}, and the function $y = \frac{0.25}{x}$ is shown in Figure~\vref{Reciprocal_Quarter}.

\FloatBarrier
As we can see, two of the graphs are stretched horizontally. Thus, all three functions look alike now. We can also see that in all three cases smaller digits have higher $fdp$ values. We already know that this feature depends on the function's shape. Let us calculate the values.

\FloatBarrier
The general Formula~\vref{l_formula} will be used. First, let $y = \frac{a}{x}$. Then, if we switch $x$ and $y$, $x = \frac{a}{y}$. Thus, $y = \frac{a}{x}$ again and $f^{-1}(x) = \frac{a}{x}$. 

$$P_k = \frac{\frac{a}{(k+1) \cdot 10^{n-1}} - \frac{a}{k \cdot 10^{n-1}}}{\frac{a}{10 \cdot 10^{n-1}} - \frac{a}{10^{n-1}}} = 
\frac{\frac{ak}{k(k+1) \cdot 10^{n-1}} - \frac{a(k+1)}{k(k+1) \cdot 10^{n-1}}}{\frac{a}{10 \cdot 10^{n-1}} - \frac{10a}{10 \cdot 10^{n-1}}} =
\frac{\frac{ak - a(k+1)}{k(k+1) \cdot 10^{n-1}}}{\frac{a - 10a}{10 \cdot 10^{n-1}}}$$

$$= \frac{\frac{ak - ak - a}{k(k+1) \cdot 10^{n-1}}}{\frac{-9a}{10 \cdot 10^{n-1}}}
=  \frac{\frac{-a}{k(k+1) \cdot 10^{n-1}}}{\frac{-9a}{10 \cdot 10^{n-1}}}
= \frac{-10 \cdot 10^{n-1}a}{-9ak(k+1) \cdot 10^{n-1}}
= \frac{10}{9k(k+1)};$$

$$ $$
Thus,
\begin{equation}
P_k = \frac{10}{9k(k+1)}
\label{rec_formula}
\end{equation}
$$ $$

\begin{figure}[h]
\centering 
\includegraphics[width=1\columnwidth]{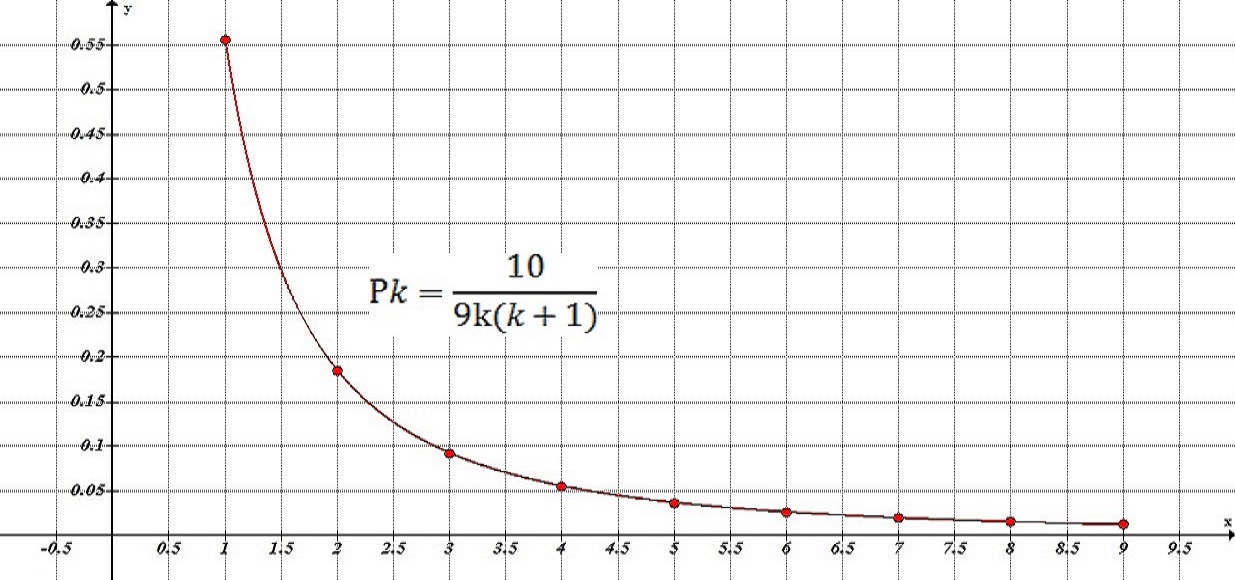} 
\caption{The function $y = \frac{10}{9k(k+1)}$}
\label{Reciprocal_Probability}
\end{figure}

We can see that the $P_k$ values of the function will be the same on any range and the constant $a$ does not change them as well. The graph of the $P_k$ function is shown in Figure~\vref{Reciprocal_Probability}.

It is obvious now that smaller first digits have higher probabilities. Besides, all the probabilities are positive numbers less than $1$. Let us prove that the sum of all nine probabilities totals to $1$. 

\FloatBarrier
\begin{proof}
$$\frac{1}{k(k+1)} = \frac{A}{k} + \frac{B}{k+1} = \frac{A(k+1)}{k(k+1)} + \frac{Bk}{k(k+1)}$$

$$= \frac{A(k+1) + Bk}{k(k+1)} = \frac{Ak+A+Bk}{k(k+1)} = \frac{k(A+B)+A}{k(k+1)};$$

$$k(A+B)+ A=1; \: A+B=0; \: A=1; \: B=-1; \: \frac{1}{k(k+1)} = \frac{1}{k} - \frac{1}{(k+1)}$$

Thus,

$$P_k = \frac{10}{9k(k+1)} = \frac{10}{9} \cdot \frac{1}{k(k+1)} = \frac{10}{9}(\frac{1}{k} - \frac{1}{k+1});$$

Then, 

$$\sum^9_{k=1}P_k = \sum^9_{k=1}\frac{10}{9}(\frac{1}{k} - \frac{1}{k+1}) =
\frac{10}{9}\sum^9_{k=1}(\frac{1}{k} - \frac{1}{k+1})=
\frac{10}{9}(\frac{1}{1} - \frac{1}{1+1} + \frac{1}{2} - \frac{1}{2+1}$$

$$ + \frac{1}{3} - \frac{1}{3+1} + \frac{1}{4} - \frac{1}{4+1} + \frac{1}{5} - 
\frac{1}{5+1} + \frac{1}{6} - \frac{1}{6+1} + \frac{1}{7} - \frac{1}{7+1} + 
\frac{1}{8} - \frac{1}{8+1} + \frac{1}{9} - \frac{1}{9+1})$$ 

$$=\frac{10}{9}(1 - \frac{1}{2} + \frac{1}{2} - \frac{1}{3}+ \frac{1}{3} - \frac{1}{4} + \frac{1}{4} - \frac{1}{5} + \frac{1}{5} - 
\frac{1}{6} + \frac{1}{6} - \frac{1}{7} + \frac{1}{7} - \frac{1}{8} + 
\frac{1}{8} - \frac{1}{9} + \frac{1}{9} - \frac{1}{10})$$

$$=\frac{10}{9}(1 - \frac{1}{10}) = \frac{10}{9} \cdot \frac{9}{10} = 1;$$

Thus, $\sum^9_{k=1}P_k = 1$, so we have obtained one of the proofs. 
\end{proof}

\FloatBarrier
Let us do one more proof that will be done by graphing after the following step is performed. 
$$P_k = \frac{10}{9}(\frac{1}{k} - \frac{1}{k+1}) 
= \frac{10}{9k} - \frac{10}{9(k+1)};$$

\begin{figure}[h]
\centering 
\includegraphics[width=1\columnwidth]{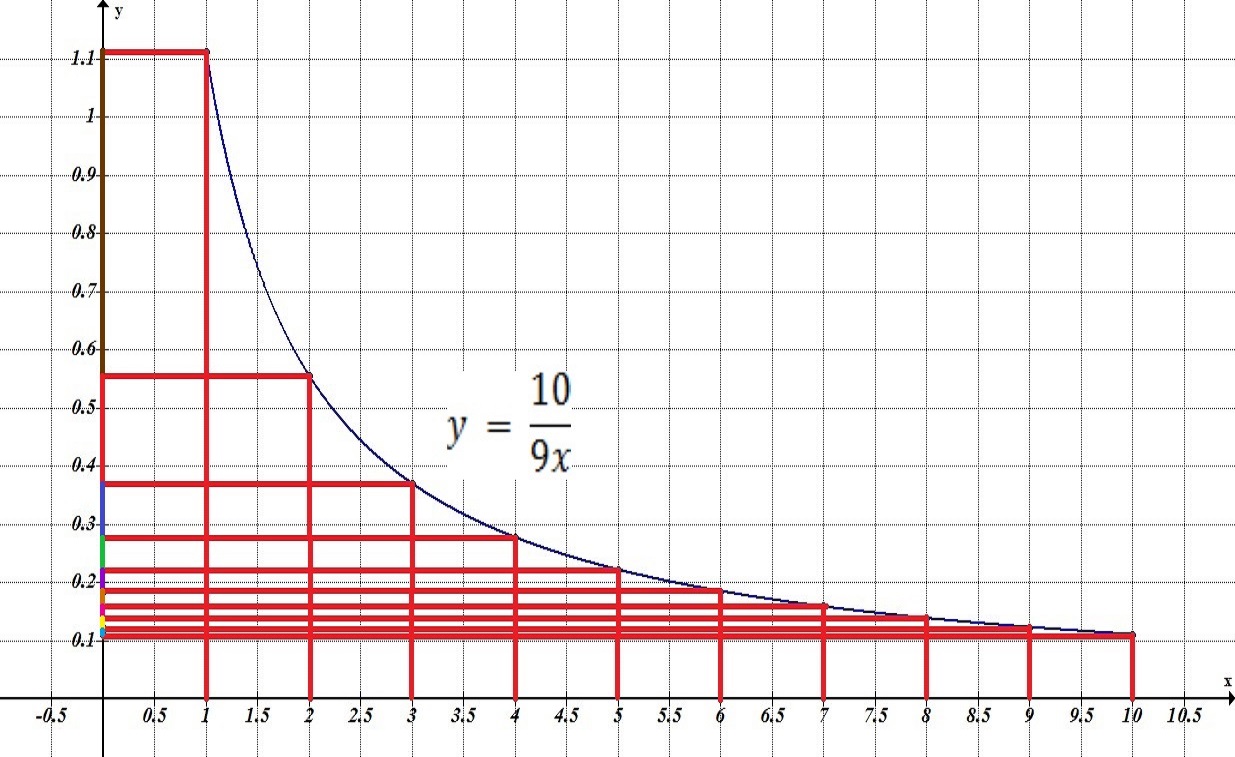} 
\caption{The function $f(x) = \frac{10}{9x}$ on the domain [1, 10]}
\label{Reciprocal_Sum}
\end{figure}

Then, we will graph the function $f(x) =  \frac{10}{9x}$ on the domain from $1$ to $10$. Figure~\vref{Reciprocal_Sum} shows the graph.

We will use the graph for locating the values $\frac{10}{9k}$ and $\frac{10}{9(k+1)}$ where $k$ is a whole number from $1$ to $9$. Both values can be found on the $y$-axis. After picking a particular $k$, subtracting $\frac{10}{9(k+1)}$ from $\frac{10}{9k}$ will give us the corresponding $P_k$ value. The graph shows the differences between the values $\frac{10}{9k}$ and $\frac{10}{9(k+1)}$ for all nine values of $k$. All the differences are equal to the corresponding $P_k$ values. They are labeled with the corresponding colors that are taken from the ruler described above.

The graph clearly shows that the sum of all nine $P_k$ values equals to the difference between $f(1)$ and $f(10)$. 
Thus, $$f(1) - f(10) = \frac{10}{9 \cdot 1} - \frac{10}{9 \cdot 10} = \frac{10}{9} - \frac{10}{90} = \frac{100}{90} - \frac{10}{90} =  \frac{90}{90} = 1;$$
It is another proof of the statement $\sum^9_{k=1}P_k = 1$.

\FloatBarrier
Let us calculate the $P_k$ values of our function and list them in Table~\vref{P_k_Reciprocal}. We should recall that the values do not depend on the constant $a$ and on the chosen range. It means that numbers presented in the table will stay the same as long as we are dealing with any reciprocal function $y = \frac{a}{x}$.

\begin{table}[hbt]
\caption{$P_k$ values of the reciprocal function}
\label{P_k_Reciprocal}
\centering
\begin{tabular}{cc}
\toprule
Digit & $P_k$ for $y = \frac{a}{x}$ \\
\midrule
1	&	0.55555556	\\
2	&	0.18518519	\\
3	&	0.09259259	\\
4	&	0.05555556	\\
5	&	0.03703704	\\
6	&	0.02645503	\\
7	&	0.01984127	\\
8	&	0.01543210  \\
9	&	0.01234568	\\
\bottomrule
Sum & 1.00000000  \\
\bottomrule
\end{tabular}
\end{table}

Before we talk about how the reciprocal function can be used to model a real-life situation, we should check how a horizontal shift would affect the $P_k$ values of the function. Our intuition should tell us that a horizontal shift of any function should not change its $P_k$ values because a horizontally shifted function will be reflected on the $y$-axis in exactly the same way as it was reflected initially. However, we just like proofs.

Let us consider a case $y = \frac{a}{x - h}$ where $h$ is a constant which shows a horizontal shift. First, we will switch $x$ and $y$. 

$$x = \frac{a}{y - h}; \: y - h = \frac{a}{x}; \: y = \frac{a}{x} + h; \: f^{-1}(x) = \frac{a}{x} + h$$

$$P_k = \frac{f^{-1}((k+1)\cdot 10^{n-1}) - f^{-1}(k \cdot 10^{n-1})}{f^{-1}(10^n) - f^{-1}(10^{n-1})}
= \frac{(\frac{a}{(k+1)\cdot 10^{n-1}} + h) - (\frac{a}{k \cdot 10^{n-1}} + h)}{(\frac{a}{10^n} + h) - (\frac{a}{10^{n-1}} + h)}$$

$$ = \frac{\frac{a}{(k+1)\cdot 10^{n-1}} + h - \frac{a}{k \cdot 10^{n-1}} - h}{\frac{a}{10^n} + h - \frac{a}{10^{n-1}} - h} = \frac{\frac{a}{(k+1)\cdot 10^{n-1}} - \frac{a}{k \cdot 10^{n-1}}}{\frac{a}{10^n} - \frac{a}{10^{n-1}}};$$

Thus, we came back to the initial $P_k$ formula. It means that the $P_k$ values of our function will not be affected by a horizontal shift. Let us talk about the actual application now.

\subsubsection{A real-life example}

We will be solving a problem about a scuba diver. The formula 

\begin{equation}
T =  \frac{525}{D - 10}
\label{r_formula}
\end{equation}

\noindent
is used to calculate the greatest amount of time in minutes that a scuba diver can take to rise toward the water surface without stopping for decompression, where $D$ is the depth, in meters, of the diver. As we can see, it is a reciprocal function with a horizontal shift. In addition, since the formula represents a real-life process, its $T$-value cannot be negative. Moreover, the denominator cannot be equal to zero. It means that the domain should be initially restricted to $(10, \infty)$.

Let us calculate the $fdp$ of the $T$-values of the function on the range $[1, 100)$. We also should find the domain that will match to this range.

$$1 = \frac {525}{D - 10}; \: D - 10 = 525; \: D = 535;$$

$$100 = \frac {525}{D - 10}; \: D - 10 = \frac {525}{100} = 5.25; \: D = 15.25;$$

We will pick a domain interval $[15, 535]$. However, if some calculated $T$-values exceed the range $[1, 100)$, we will not include them in our further calculations. We also should check that selected by us domain is located inside of the initial domain $(10, \infty)$, which is the case. 

\begin{figure}[h]
\centering 
\includegraphics[width=1\columnwidth]{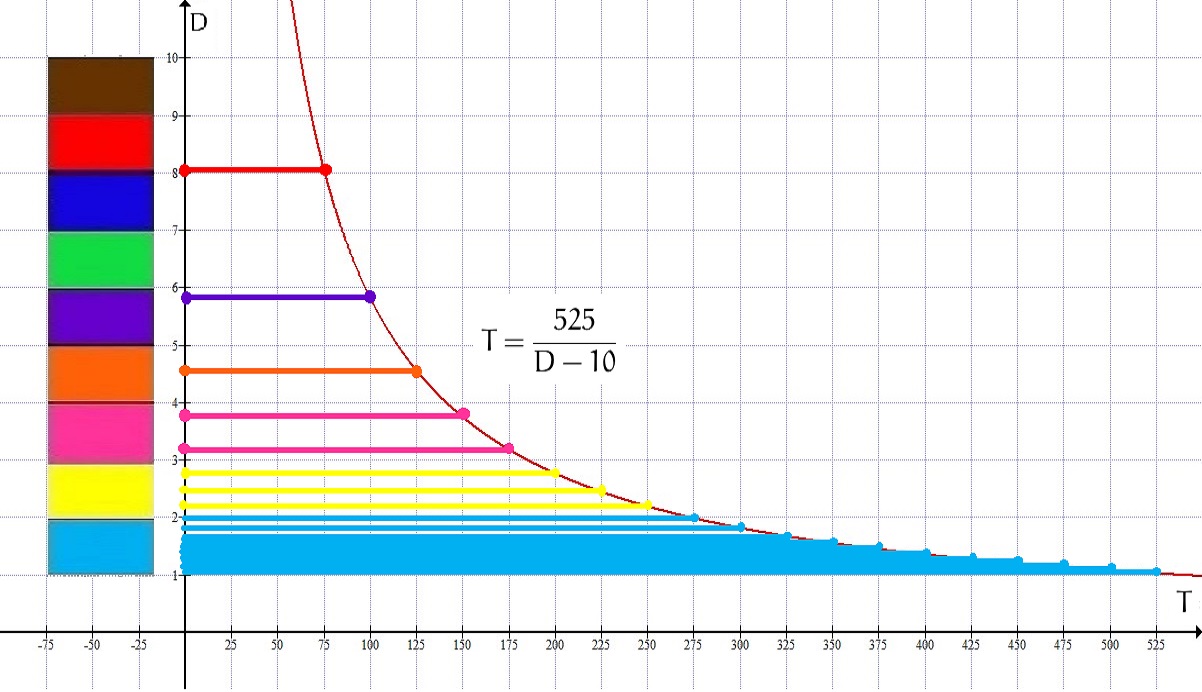} 
\caption{The function $T =  \frac{525}{D - 10}$}
\label{Reciprocal_Example}
\end{figure}

We will plug every integer from $15$ to $535$ into the variable $D$ in Formula~\vref{r_formula}, calculate the corresponding $T$-value, retrieve the $T$-value's first digit, and list all the results in Table~\vref{T(D)}. However, since x-values are changing very slow on the y-interval [10, 100), we will look at the graph of the function on the range [1, 10) only. The graph is shown in Figure~\vref{Reciprocal_Example}.

\newpage
\begin{longtable}[c]{lllllllll} 
 \caption{Calculated $T(D)$ values and their first digits}
 \label{T(D)}\\
\toprule
$D$- & $T$- & First & $D$- & $T$- & First & $D$- & $T$- & First \\
values & values & digit & values & values & digit & values & values & digit \\ 
\midrule
\endfirsthead 

\toprule
$D$- & $T$- & First & $D$- & $T$- & First & $D$- & $T$- & First \\
values & values & digit & values & values & digit & values & values & digit \\ 
\midrule
\endhead 
15	&	105	&		&	65	&	9.545455	&	9	&	115	&	5	&	5	\\		
16	&	87.5	&	8	&	66	&	9.375	&	9	&	116	&	4.95283	&	4	\\		
17	&	75	&	7	&	67	&	9.210526	&	9	&	117	&	4.906542	&	4	\\		
18	&	65.625	&	6	&	68	&	9.051724	&	9	&	118	&	4.861111	&	4	\\		
19	&	58.33333	&	5	&	69	&	8.898305	&	8	&	119	&	4.816514	&	4	\\		
20	&	52.5	&	5	&	70	&	8.75	&	8	&	120	&	4.772727	&	4	\\		
21	&	47.72727	&	4	&	71	&	8.606557	&	8	&	121	&	4.72973	&	4	\\		
22	&	43.75	&	4	&	72	&	8.467742	&	8	&	122	&	4.6875	&	4	\\		
23	&	40.38462	&	4	&	73	&	8.333333	&	8	&	123	&	4.646018	&	4	\\		
24	&	37.5	&	3	&	74	&	8.203125	&	8	&	124	&	4.605263	&	4	\\		
25	&	35	&	3	&	75	&	8.076923	&	8	&	125	&	4.565217	&	4	\\		
26	&	32.8125	&	3	&	76	&	7.954545	&	7	&	126	&	4.525862	&	4	\\		
27	&	30.88235	&	3	&	77	&	7.835821	&	7	&	127	&	4.487179	&	4	\\		
28	&	29.16667	&	2	&	78	&	7.720588	&	7	&	128	&	4.449153	&	4	\\		
29	&	27.63158	&	2	&	79	&	7.608696	&	7	&	129	&	4.411765	&	4	\\		
30	&	26.25	&	2	&	80	&	7.5	&	7	&	130	&	4.375	&	4	\\		
31	&	25	&	2	&	81	&	7.394366	&	7	&	131	&	4.338843	&	4	\\		
32	&	23.86364	&	2	&	82	&	7.291667	&	7	&	132	&	4.303279	&	4	\\		
33	&	22.82609	&	2	&	83	&	7.191781	&	7	&	133	&	4.268293	&	4	\\		
34	&	21.875	&	2	&	84	&	7.094595	&	7	&	134	&	4.233871	&	4	\\		
35	&	21	&	2	&	85	&	7	&	7	&	135	&	4.2	&	4	\\		
36	&	20.19231	&	2	&	86	&	6.907895	&	6	&	136	&	4.166667	&	4	\\		
37	&	19.44444	&	1	&	87	&	6.818182	&	6	&	137	&	4.133858	&	4	\\		
38	&	18.75	&	1	&	88	&	6.730769	&	6	&	138	&	4.101563	&	4	\\		
39	&	18.10345	&	1	&	89	&	6.64557	&	6	&	139	&	4.069767	&	4	\\		
40	&	17.5	&	1	&	90	&	6.5625	&	6	&	140	&	4.038462	&	4	\\		
41	&	16.93548	&	1	&	91	&	6.481481	&	6	&	141	&	4.007634	&	4	\\		
42	&	16.40625	&	1	&	92	&	6.402439	&	6	&	142	&	3.977273	&	3	\\		
43	&	15.90909	&	1	&	93	&	6.325301	&	6	&	143	&	3.947368	&	3	\\		
44	&	15.44118	&	1	&	94	&	6.25	&	6	&	144	&	3.91791	&	3	\\		
45	&	15	&	1	&	95	&	6.176471	&	6	&	145	&	3.888889	&	3	\\		
46	&	14.58333	&	1	&	96	&	6.104651	&	6	&	146	&	3.860294	&	3	\\		
47	&	14.18919	&	1	&	97	&	6.034483	&	6	&	147	&	3.832117	&	3	\\		
48	&	13.81579	&	1	&	98	&	5.965909	&	5	&	148	&	3.804348	&	3	\\		
49	&	13.46154	&	1	&	99	&	5.898876	&	5	&	149	&	3.776978	&	3	\\		
50	&	13.125	&	1	&	100	&	5.833333	&	5	&	150	&	3.75	&	3	\\		
51	&	12.80488	&	1	&	101	&	5.769231	&	5	&	151	&	3.723404	&	3	\\		
52	&	12.5	&	1	&	102	&	5.706522	&	5	&	152	&	3.697183	&	3	\\		
53	&	12.2093	&	1	&	103	&	5.645161	&	5	&	153	&	3.671329	&	3	\\		
54	&	11.93182	&	1	&	104	&	5.585106	&	5	&	154	&	3.645833	&	3	\\		
55	&	11.66667	&	1	&	105	&	5.526316	&	5	&	155	&	3.62069	&	3	\\		
56	&	11.41304	&	1	&	106	&	5.46875	&	5	&	156	&	3.59589	&	3	\\		
57	&	11.17021	&	1	&	107	&	5.412371	&	5	&	157	&	3.571429	&	3	\\		
58	&	10.9375	&	1	&	108	&	5.357143	&	5	&	158	&	3.547297	&	3	\\		
59	&	10.71429	&	1	&	109	&	5.30303	&	5	&	159	&	3.52349	&	3	\\		
60	&	10.5	&	1	&	110	&	5.25	&	5	&	160	&	3.5	&	3	\\		
61	&	10.29412	&	1	&	111	&	5.19802	&	5	&	161	&	3.476821	&	3	\\		
62	&	10.09615	&	1	&	112	&	5.147059	&	5	&	162	&	3.453947	&	3	\\		
63	&	9.90566	&	9	&	113	&	5.097087	&	5	&	163	&	3.431373	&	3	\\		
64	&	9.722222	&	9	&	114	&	5.048077	&	5	&	164	&	3.409091	&	3	\\	\\	\\
165	&	3.387097	&	3	&	215	&	2.560976	&	2	&	265	&	2.058824	&	2	\\		
166	&	3.365385	&	3	&	216	&	2.548544	&	2	&	266	&	2.050781	&	2	\\		
167	&	3.343949	&	3	&	217	&	2.536232	&	2	&	267	&	2.042802	&	2	\\		
168	&	3.322785	&	3	&	218	&	2.524038	&	2	&	268	&	2.034884	&	2	\\		
169	&	3.301887	&	3	&	219	&	2.511962	&	2	&	269	&	2.027027	&	2	\\		
170	&	3.28125	&	3	&	220	&	2.5	&	2	&	270	&	2.019231	&	2	\\		
171	&	3.26087	&	3	&	221	&	2.488152	&	2	&	271	&	2.011494	&	2	\\		
172	&	3.240741	&	3	&	222	&	2.476415	&	2	&	272	&	2.003817	&	2	\\		
173	&	3.220859	&	3	&	223	&	2.464789	&	2	&	273	&	1.996198	&	1	\\		
174	&	3.20122	&	3	&	224	&	2.453271	&	2	&	274	&	1.988636	&	1	\\		
175	&	3.181818	&	3	&	225	&	2.44186	&	2	&	275	&	1.981132	&	1	\\		
176	&	3.162651	&	3	&	226	&	2.430556	&	2	&	276	&	1.973684	&	1	\\		
177	&	3.143713	&	3	&	227	&	2.419355	&	2	&	277	&	1.966292	&	1	\\		
178	&	3.125	&	3	&	228	&	2.408257	&	2	&	278	&	1.958955	&	1	\\		
179	&	3.106509	&	3	&	229	&	2.39726	&	2	&	279	&	1.951673	&	1	\\		
180	&	3.088235	&	3	&	230	&	2.386364	&	2	&	280	&	1.944444	&	1	\\		
181	&	3.070175	&	3	&	231	&	2.375566	&	2	&	281	&	1.937269	&	1	\\		
182	&	3.052326	&	3	&	232	&	2.364865	&	2	&	282	&	1.930147	&	1	\\		
183	&	3.034682	&	3	&	233	&	2.35426	&	2	&	283	&	1.923077	&	1	\\		
184	&	3.017241	&	3	&	234	&	2.34375	&	2	&	284	&	1.916058	&	1	\\		
185	&	3	&	3	&	235	&	2.333333	&	2	&	285	&	1.909091	&	1	\\		
186	&	2.982955	&	2	&	236	&	2.323009	&	2	&	286	&	1.902174	&	1	\\		
187	&	2.966102	&	2	&	237	&	2.312775	&	2	&	287	&	1.895307	&	1	\\		
188	&	2.949438	&	2	&	238	&	2.302632	&	2	&	288	&	1.888489	&	1	\\		
189	&	2.932961	&	2	&	239	&	2.292576	&	2	&	289	&	1.88172	&	1	\\		
190	&	2.916667	&	2	&	240	&	2.282609	&	2	&	290	&	1.875	&	1	\\		
191	&	2.900552	&	2	&	241	&	2.272727	&	2	&	291	&	1.868327	&	1	\\		
192	&	2.884615	&	2	&	242	&	2.262931	&	2	&	292	&	1.861702	&	1	\\		
193	&	2.868852	&	2	&	243	&	2.253219	&	2	&	293	&	1.855124	&	1	\\		
194	&	2.853261	&	2	&	244	&	2.24359	&	2	&	294	&	1.848592	&	1	\\		
195	&	2.837838	&	2	&	245	&	2.234043	&	2	&	295	&	1.842105	&	1	\\		
196	&	2.822581	&	2	&	246	&	2.224576	&	2	&	296	&	1.835664	&	1	\\		
197	&	2.807487	&	2	&	247	&	2.21519	&	2	&	297	&	1.829268	&	1	\\		
198	&	2.792553	&	2	&	248	&	2.205882	&	2	&	298	&	1.822917	&	1	\\		
199	&	2.777778	&	2	&	249	&	2.196653	&	2	&	299	&	1.816609	&	1	\\		
200	&	2.763158	&	2	&	250	&	2.1875	&	2	&	300	&	1.810345	&	1	\\		
201	&	2.748691	&	2	&	251	&	2.178423	&	2	&	301	&	1.804124	&	1	\\		
202	&	2.734375	&	2	&	252	&	2.169421	&	2	&	302	&	1.797945	&	1	\\		
203	&	2.720207	&	2	&	253	&	2.160494	&	2	&	303	&	1.791809	&	1	\\		
204	&	2.706186	&	2	&	254	&	2.151639	&	2	&	304	&	1.785714	&	1	\\		
205	&	2.692308	&	2	&	255	&	2.142857	&	2	&	305	&	1.779661	&	1	\\		
206	&	2.678571	&	2	&	256	&	2.134146	&	2	&	306	&	1.773649	&	1	\\		
207	&	2.664975	&	2	&	257	&	2.125506	&	2	&	307	&	1.767677	&	1	\\		
208	&	2.651515	&	2	&	258	&	2.116935	&	2	&	308	&	1.761745	&	1	\\		
209	&	2.638191	&	2	&	259	&	2.108434	&	2	&	309	&	1.755853	&	1	\\		
210	&	2.625	&	2	&	260	&	2.1	&	2	&	310	&	1.75	&	1	\\		
211	&	2.61194	&	2	&	261	&	2.091633	&	2	&	311	&	1.744186	&	1	\\		
212	&	2.59901	&	2	&	262	&	2.083333	&	2	&	312	&	1.738411	&	1	\\		
213	&	2.586207	&	2	&	263	&	2.075099	&	2	&	313	&	1.732673	&	1	\\		
214	&	2.573529	&	2	&	264	&	2.066929	&	2	&	314	&	1.726974	&	1	\\	\\	\\	\\
315	&	1.721311	&	1	&	365	&	1.478873	&	1	&	415	&	1.296296	&	1	\\		
316	&	1.715686	&	1	&	366	&	1.474719	&	1	&	416	&	1.293103	&	1	\\		
317	&	1.710098	&	1	&	367	&	1.470588	&	1	&	417	&	1.289926	&	1	\\		
318	&	1.704545	&	1	&	368	&	1.46648	&	1	&	418	&	1.286765	&	1	\\		
319	&	1.699029	&	1	&	369	&	1.462396	&	1	&	419	&	1.283619	&	1	\\		
320	&	1.693548	&	1	&	370	&	1.458333	&	1	&	420	&	1.280488	&	1	\\		
321	&	1.688103	&	1	&	371	&	1.454294	&	1	&	421	&	1.277372	&	1	\\		
322	&	1.682692	&	1	&	372	&	1.450276	&	1	&	422	&	1.274272	&	1	\\		
323	&	1.677316	&	1	&	373	&	1.446281	&	1	&	423	&	1.271186	&	1	\\		
324	&	1.671975	&	1	&	374	&	1.442308	&	1	&	424	&	1.268116	&	1	\\		
325	&	1.666667	&	1	&	375	&	1.438356	&	1	&	425	&	1.26506	&	1	\\		
326	&	1.661392	&	1	&	376	&	1.434426	&	1	&	426	&	1.262019	&	1	\\		
327	&	1.656151	&	1	&	377	&	1.430518	&	1	&	427	&	1.258993	&	1	\\		
328	&	1.650943	&	1	&	378	&	1.42663	&	1	&	428	&	1.255981	&	1	\\		
329	&	1.645768	&	1	&	379	&	1.422764	&	1	&	429	&	1.252983	&	1	\\		
330	&	1.640625	&	1	&	380	&	1.418919	&	1	&	430	&	1.25	&	1	\\		
331	&	1.635514	&	1	&	381	&	1.415094	&	1	&	431	&	1.247031	&	1	\\		
332	&	1.630435	&	1	&	382	&	1.41129	&	1	&	432	&	1.244076	&	1	\\		
333	&	1.625387	&	1	&	383	&	1.407507	&	1	&	433	&	1.241135	&	1	\\		
334	&	1.62037	&	1	&	384	&	1.403743	&	1	&	434	&	1.238208	&	1	\\		
335	&	1.615385	&	1	&	385	&	1.4	&	1	&	435	&	1.235294	&	1	\\		
336	&	1.610429	&	1	&	386	&	1.396277	&	1	&	436	&	1.232394	&	1	\\		
337	&	1.605505	&	1	&	387	&	1.392573	&	1	&	437	&	1.229508	&	1	\\		
338	&	1.60061	&	1	&	388	&	1.388889	&	1	&	438	&	1.226636	&	1	\\		
339	&	1.595745	&	1	&	389	&	1.385224	&	1	&	439	&	1.223776	&	1	\\		
340	&	1.590909	&	1	&	390	&	1.381579	&	1	&	440	&	1.22093	&	1	\\		
341	&	1.586103	&	1	&	391	&	1.377953	&	1	&	441	&	1.218097	&	1	\\		
342	&	1.581325	&	1	&	392	&	1.374346	&	1	&	442	&	1.215278	&	1	\\		
343	&	1.576577	&	1	&	393	&	1.370757	&	1	&	443	&	1.212471	&	1	\\		
344	&	1.571856	&	1	&	394	&	1.367188	&	1	&	444	&	1.209677	&	1	\\		
345	&	1.567164	&	1	&	395	&	1.363636	&	1	&	445	&	1.206897	&	1	\\		
346	&	1.5625	&	1	&	396	&	1.360104	&	1	&	446	&	1.204128	&	1	\\		
347	&	1.557864	&	1	&	397	&	1.356589	&	1	&	447	&	1.201373	&	1	\\		
348	&	1.553254	&	1	&	398	&	1.353093	&	1	&	448	&	1.19863	&	1	\\		
349	&	1.548673	&	1	&	399	&	1.349614	&	1	&	449	&	1.1959	&	1	\\		
350	&	1.544118	&	1	&	400	&	1.346154	&	1	&	450	&	1.193182	&	1	\\		
351	&	1.539589	&	1	&	401	&	1.342711	&	1	&	451	&	1.190476	&	1	\\		
352	&	1.535088	&	1	&	402	&	1.339286	&	1	&	452	&	1.187783	&	1	\\		
353	&	1.530612	&	1	&	403	&	1.335878	&	1	&	453	&	1.185102	&	1	\\		
354	&	1.526163	&	1	&	404	&	1.332487	&	1	&	454	&	1.182432	&	1	\\		
355	&	1.521739	&	1	&	405	&	1.329114	&	1	&	455	&	1.179775	&	1	\\		
356	&	1.517341	&	1	&	406	&	1.325758	&	1	&	456	&	1.17713	&	1	\\		
357	&	1.512968	&	1	&	407	&	1.322418	&	1	&	457	&	1.174497	&	1	\\		
358	&	1.508621	&	1	&	408	&	1.319095	&	1	&	458	&	1.171875	&	1	\\		
359	&	1.504298	&	1	&	409	&	1.315789	&	1	&	459	&	1.169265	&	1	\\		
360	&	1.5	&	1	&	410	&	1.3125	&	1	&	460	&	1.166667	&	1	\\		
361	&	1.495726	&	1	&	411	&	1.309227	&	1	&	461	&	1.16408	&	1	\\		
362	&	1.491477	&	1	&	412	&	1.30597	&	1	&	462	&	1.161504	&	1	\\		
363	&	1.487252	&	1	&	413	&	1.30273	&	1	&	463	&	1.15894	&	1	\\		
364	&	1.483051	&	1	&	414	&	1.299505	&	1	&	464	&	1.156388	&	1	\\	\\	\\	\\
465	&	1.153846	&	1	&	489	&	1.096033	&	1	&	513	&	1.043738	&	1	\\		
466	&	1.151316	&	1	&	490	&	1.09375	&	1	&	514	&	1.041667	&	1	\\		
467	&	1.148796	&	1	&	491	&	1.091476	&	1	&	515	&	1.039604	&	1	\\		
468	&	1.146288	&	1	&	492	&	1.089212	&	1	&	516	&	1.037549	&	1	\\		
469	&	1.143791	&	1	&	493	&	1.086957	&	1	&	517	&	1.035503	&	1	\\		
470	&	1.141304	&	1	&	494	&	1.084711	&	1	&	518	&	1.033465	&	1	\\		
471	&	1.138829	&	1	&	495	&	1.082474	&	1	&	519	&	1.031434	&	1	\\		
472	&	1.136364	&	1	&	496	&	1.080247	&	1	&	520	&	1.029412	&	1	\\		
473	&	1.133909	&	1	&	497	&	1.078029	&	1	&	521	&	1.027397	&	1	\\		
474	&	1.131466	&	1	&	498	&	1.07582	&	1	&	522	&	1.025391	&	1	\\		
475	&	1.129032	&	1	&	499	&	1.07362	&	1	&	523	&	1.023392	&	1	\\		
476	&	1.126609	&	1	&	500	&	1.071429	&	1	&	524	&	1.021401	&	1	\\		
477	&	1.124197	&	1	&	501	&	1.069246	&	1	&	525	&	1.019417	&	1	\\		
478	&	1.121795	&	1	&	502	&	1.067073	&	1	&	526	&	1.017442	&	1	\\		
479	&	1.119403	&	1	&	503	&	1.064909	&	1	&	527	&	1.015474	&	1	\\		
480	&	1.117021	&	1	&	504	&	1.062753	&	1	&	528	&	1.013514	&	1	\\		
481	&	1.11465	&	1	&	505	&	1.060606	&	1	&	529	&	1.011561	&	1	\\		
482	&	1.112288	&	1	&	506	&	1.058468	&	1	&	530	&	1.009615	&	1	\\		
483	&	1.109937	&	1	&	507	&	1.056338	&	1	&	531	&	1.007678	&	1	\\		
484	&	1.107595	&	1	&	508	&	1.054217	&	1	&	532	&	1.005747	&	1	\\		
485	&	1.105263	&	1	&	509	&	1.052104	&	1	&	533	&	1.003824	&	1	\\		
486	&	1.102941	&	1	&	510	&	1.05	&	1	&	534	&	1.001908	&	1	\\		
487	&	1.100629	&	1	&	511	&	1.047904	&	1	&	535	&	1	&	1	\\		
488	&	1.098326	&	1	&	512	&	1.045817	&	1	&						\\		
\bottomrule \\
\end{longtable}

We have all the information we need to perform our final calculations and comparison. Let us summarize all the numbers using Table~\vref{T(D)_analysis}. 

\begin{table}[hbt]
\caption{Analysis of $T(D)$ values first digits}
\label{T(D)_analysis}
\centering
\begin{tabular}{cccc}
\toprule
Digit & Count & Frequency       & $P_k$ values    \\  
       &           & (Count / Sum) & for $y = a / x$  \\ 
\midrule
1	&	289	&	0.55576923	&	0.55555556	\\
2	&	96	&	0.18461538	&	0.18518519	\\
3	&	48	&	0.09230769	&	0.09259259	\\
4	&	29	&	0.05576923	&	0.05555556	\\
5	&	20	&	0.03846154	&	0.03703704	\\
6	&	13	&	0.02500000	&	0.02645503	\\
7	&	11	&	0.02115385	&	0.01984127	\\
8	&	8	&	0.01538462	&	0.01543210	\\
9	&	6	&	0.01153846	&	0.01234568	\\
\bottomrule
Sum	&	520	&	1.00000000	&	1.00000000	\\
\bottomrule
\end{tabular}
\end{table}

As we can see, the corresponding numbers are really close to each other. We got such a nice result due to that fact that our application problem calculation had $521$ steps.

\clearpage
\subsection{Summary of $P_k$ values of continuous functions}

We know by now that all functions that we have discussed have certain rules for their $P_k$ values. The values are either increasing, decreasing, or they stay constant. Let us analyze how all of our functions are broken into these three groups.

Let us summarize everything that we have discovered in this work about continuous functions. Our intuition tells us that behavior of $P_k$ numbers has something to do with the function's shape. In particular, we should analyze how the slope of a tangent to the function line affects its $P_k$ numbers. We will do this analysis for all our continuous functions by creating one picture per group. 

\begin{figure}[h]
\centering 
\includegraphics[width=1\columnwidth]{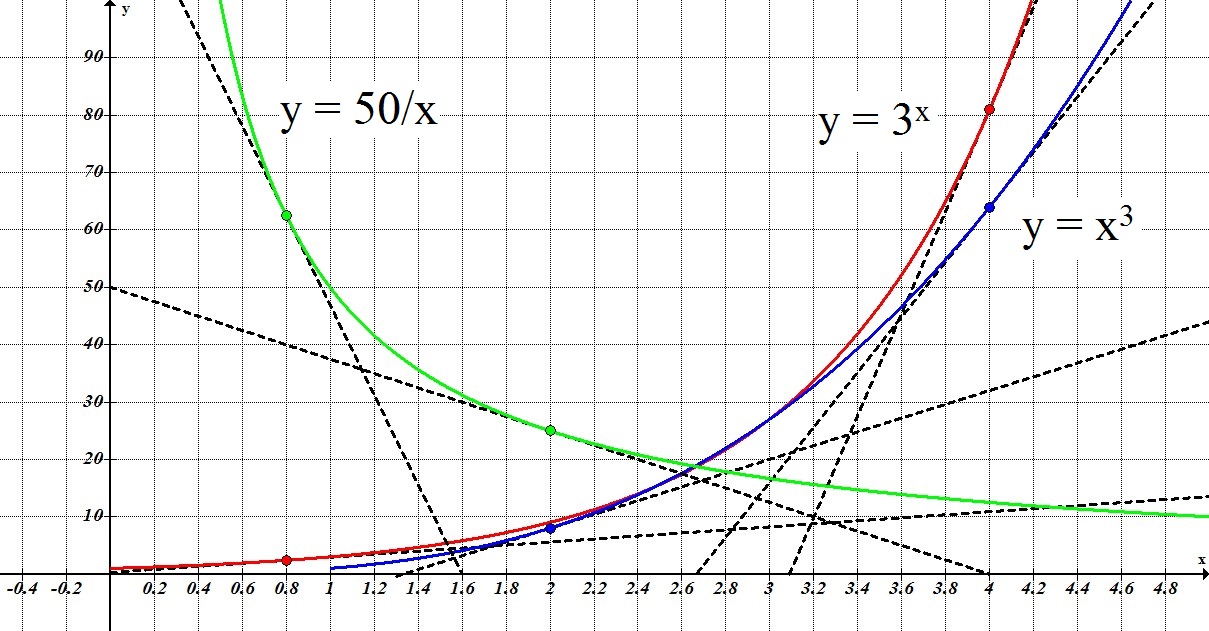} 
\caption{Graphs of the functions $y = 3^x$, $y = x^3$, and $y = \frac{50}{x}$}
\label{Lines_Increasing}
\end{figure}

Figure~\vref{Lines_Increasing} shows us graphs of the functions with decreasing $P_k$ numbers. All of them represent their sections and one representative per section should be enough. They are: 

$$y = 3^x; \: y = x^3; \: y = \frac{50}{x}$$

Since the function $y = \frac{a}{x}$ is decreasing on its entire domain if $a$ is positive, it is decreasing on the interval $(0, a]$. However, we are moving from $a$ to $0$ along the $x$-axis while analyzing the $y$-values of this function. We are doing it because our goal is to see how the first digits of $y$-values change while the $y$-values are increasing. Thus, we have to move in the same direction while analyzing the tangent lines of the function. As we can see, tangent lines of all the functions in Figure~\vref{Lines_Increasing} increase their slopes while $y$-coordinates of their functions are increasing.

\begin{figure}[h]
\centering 
\includegraphics[width=1\columnwidth]{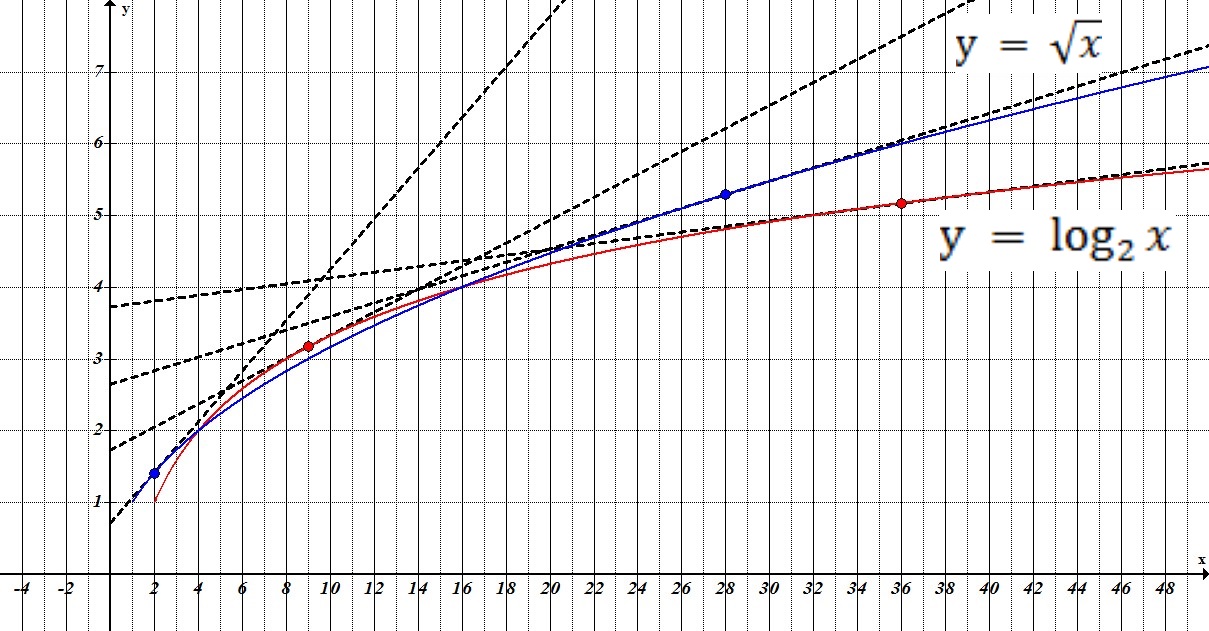} 
\caption{Graphs of the functions $y = \sqrt x$ and $y = \log_2⁡(x)$}
\label{Lines_Decreasing}
\end{figure}

Let us look at another picture. The functions with increasing $P_k$ numbers are shown in Figure~\vref{Lines_Decreasing}. They are: 

$$y = \sqrt x; \: y = \log_2⁡(x)$$

Again, one representative per section is enough to see the whole picture. It is obvious that the slopes of the  tangent lines are decreasing while the functions $y$-coordinates are increasing.

We should also cover the case when our function is linear. There is no need to graph the function $y = ax$ because it is obvious that all tangent lines will have the same slopes as the function does due to the function's nature. 

Thus, we are ready to draw certain conclusions. Functions with decreasing $P_k$ numbers have tangent lines with increasing slopes and functions with increasing $P_k$ numbers have tangent lines with decreasing slopes. Moreover, functions with constant $P_k$ numbers have tangent lines with constant slopes. We know that changes in directions of tangent lines are directly relevant to first and second derivatives of the initial functions. Let us find them using our knowledge of Calculus. 

First, we will find the first and second derivatives of our continuous functions and list the results. Then, we will analyze whether the first derivatives are increasing or decreasing and whether the second derivatives are positive or negative. These analyses will be done on the chosen domain. The conclusions are directly related to each other and to our $P_k$ numbers.
\\

$f(x) = a^x$, $a$ > $1$, $x \geq 0$  

$f'(x) = a^x \cdot \ln (a)$ - increasing

$f''(x) = \ln (⁡     a) \cdot a^x \cdot \ln (⁡     a) = (\ln (⁡     a))^2 \cdot a^x$ > $0$ 
\\

$f(x) = x^a$, $a$ > $1$, $x \geq 1$
 
$f'(x) = ax^{a - 1}$ - increasing

$f''(x) = a(a - 1)x^{a - 1}$ > $0$ 
\\

$f(x) = mx$, $m$ > $0$, $x \geq  \frac{1}{m}$
 
$f'(x) = m$ - constant

$f''(x) = 0$
\\
 
$f(x) = \sqrt [a]{x} = x^{\frac{1}{a}}$, $a$ > $1$, $x \geq 1$
 
$f'(x) = \frac{1}{a} \cdot x^{\frac{1}{a} - 1} = 
\frac{1}{a} \cdot x^{\frac{1}{a}} \cdot x^{-1} =  
\frac {\sqrt[a]{x}}{ax}$ - decreasing

$f''(x) = (\frac{1}{a} \cdot x^{\frac{1-a}{a}})' = 
\frac{1}{a} \cdot \frac{1 - a}{a} \cdot x^{\frac{1 - a}{a}-1} = 
\frac{1 - a}{a^2} \cdot x^{\frac{1}{a}-2} = 
\frac{\sqrt [a]{x} \cdot (1 - a)}{a^2x^2}$ < $0$
\\

$f(x) = \log_a ⁡(x)$, $a$ > $1$, $x \geq a$ 

$f'(x) = \frac{1}{\ln (a) \cdot x}$ - decreasing

$f''(x) = (\frac{1}{\ln (a)} \cdot x^{-1})' = 
(-1) \cdot \frac{1}{\ln (a)} \cdot x^{-2} =  
- \frac{1}{\ln (a) \cdot x^2}$ < $0$ 
\\

$f(x) = \frac{a}{x} = ax^{-1}$, $a$ > $0$, $0$ < $x \leq a$
 
$f'(x) = (-1)ax^{-2} = - \frac{a}{x^2}$ - decreasing

$f''(x) = ((-1)ax^{-2})' = (-1)(-2)ax^{-3} = \frac{2a}{x^3} $ > $0$
\\
 
It is obvious that all functions with decreasing $P_k$ numbers have positive second derivatives and those with increasing $P_k$ numbers have negative second derivatives. In addition, the second derivative of the function with constant $P_k$ numbers equals to zero.

\begin{figure}[h]
\centering 
\includegraphics[width=1\columnwidth]{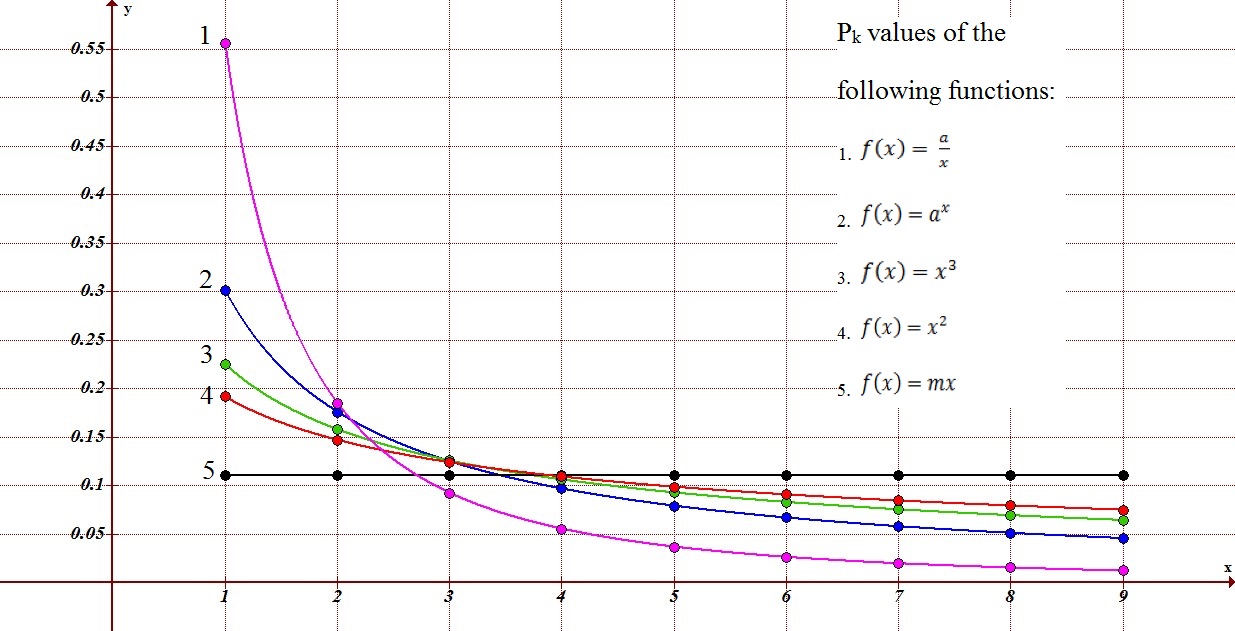} 
\caption{$P_k$ values of a linear function and functions with positive second derivatives}
\label{Combined_Probability_1}
\end{figure}

\FloatBarrier
Let us graph all the $P_k$ values of the functions with positive second derivatives and the $P_k$ numbers of a linear function on the same $xy$-plane for a comparison. Figure~\vref{Combined_Probability_1} shows the graphs. As we can see, all the decreasing $P_k$ functions change their values differently. Let us analyze these differences. 

We will use our knowledge of Calculus again to find the limits of the ratio $\frac {f(x+1)}{f(x)}$ for all the functions whose $P_k$ numbers were just graphed. All the $x$-variables of our limits will approach to positive infinity except the $x$-variable of the limit of the reciprocal function. Since the $x$-values of the function are decreasing while the $y$-values are increasing in our domain, the $x$-variable of this function's ratio's limit will approach to zero. For the same reason, we will find the limit of the ratio $\frac {f(x)}{f(x+1)}$  for the reciprocal function instead of the one mentioned above. 
\\ \\
\noindent
$\lim_{x \to 0}\frac{\frac{a}{x}}{\frac{a}{x+1}} = \lim_{x \to 0}\frac{x+1}{x} = 
\lim_{x \to 0}(1 + \frac{1}{x}) = \infty;$
\\ \\
$\lim_{x \to \infty}\frac{a^{x+1}}{a^x} = \lim_{x \to \infty}a = a, a $ > $ 1;$
\\ \\
$\lim_{x \to \infty}\frac{(x+1)^a}{x^a} = \lim_{x \to \infty}(\frac{x+1}{x})^a = 
\lim_{x \to \infty}(1 + \frac{1}{x})^a = 1^a = 1;$
\\ \\
$\lim_{x \to \infty}\frac{m(x+1)}{mx} = \lim_{x \to \infty}\frac{x+1}{x} = 
\lim_{x \to \infty}(1 + \frac{1}{x}) = 1;$ \\

Let us summarize the results. It is ovbious that all exponential functions have the same probabilities and all their ratios are equal to their bases, which are greater than $1$. $P_k$ numbers of the reciprocal function start with bigger values than those of the exponential functions and their ratios approach to infinity. $P_k$ numbers of the power functions start with smaller values than those of the exponential ones and their ratios approach to $1$. Thus, we can see a pattern that higher limits belong to functions with higher $P_1$ values.

Moreover, we might notice that $P_k$ numbers of power functions with higher exponents start with bigger values than those with lower exponents. However, we do not know how high the $P_k$ numbers will start if the exponent gets infinitely big. Let us take another limit to answer this question. We will check what will happen to the ratio $\frac{f(x+1)}{f(x)}$, where $f(x) = x^a$ if both $x$ and $a$ approach to infinity.

$$\lim_{a,x \to \infty}\frac{(x+1)^a}{x^a} = \lim_{a,x \to \infty}(\frac{x+1}{x})^a = 
\lim_{a,x \to \infty}(1 + \frac{1}{x})^a = e;$$

The last step is valid due to the fact that \cite{r_limit}

\begin{equation}
\lim_{x \to \infty}(1 + \frac{1}{x})^x = e
\end{equation}

According to the result, we can make a hypothesis that as the exponent of a power function gets bigger, its $P_1$ approaches to $P_1$ of an exponential function, but never gets equal or bigger than it. 

\begin{figure}[h]
\centering 
\includegraphics[width=1\columnwidth]{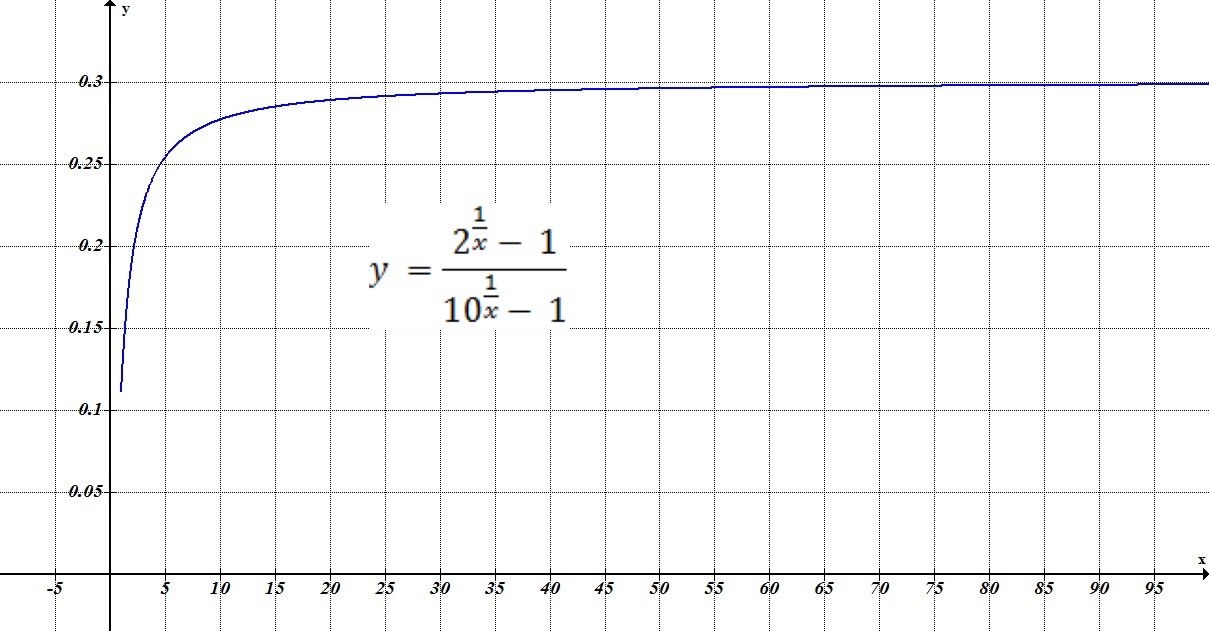} 
\caption{The function $y = \frac{2^\frac{1}{x} - 1}{10^\frac{1}{x} - 1}$ on the domain $[1, \infty)$}
\label{Power_Approaches_To_Exponential}
\end{figure}

Based on Formula~\vref{pow_formula}, 

$$P_1 = \frac{2^{1/a} - 1^{1/a}}{10^{1/a} - 1} = \frac{2^{1/a} - 1}{10^{1/a} - 1},$$  

\noindent
where $a$ is the initial function's exponent. Let us graph the function $y = \frac{2^{1/x} - 1}{10^{1/x} - 1}$ on the domain $[1, \infty)$ and see where it approaches to while $x$ approaches to infinity. Figure~\vref{Power_Approaches_To_Exponential} shows the graph.

As we can see, the function mentioned above approaches to a number which is close to the number $0.3$ and never gets above it. Our intuition tells us that this number is $P_1$ of the exponential function, which is 

$$\log_{10}(\frac{1+1}{1}) = \log_{10}(2) = \frac{\ln (a)}{\ln (10)} \approx 0.301;$$

Let us prove it algebraically. First, let $x = \frac{1}{c}$ . Then, $c \to 0$ as $x \to \infty$.

$$\lim_{x \to \infty}\frac{2^\frac{1}{x} - 1}{10^\frac{1}{x} - 1} = 
\lim_{c \to 0}\frac{2^c - 1}{10^c - 1} = \lim_{c \to 0}\frac{(2^c - 1)'}{(10^c - 1)'} = 
\lim_{c \to 0}\frac{2^c \cdot \ln (a)}{10^c \cdot \ln (10)} = \frac{\ln (2)}{\ln (10)};$$

Thus, our assumption was correct. In addition, we can make a conclusion that $P_1$ values of the reciprocal, exponential, and power functions have their boundaries and these boundaries do not overlap each other.

\begin{figure}[h]
\centering 
\includegraphics[width=1\columnwidth]{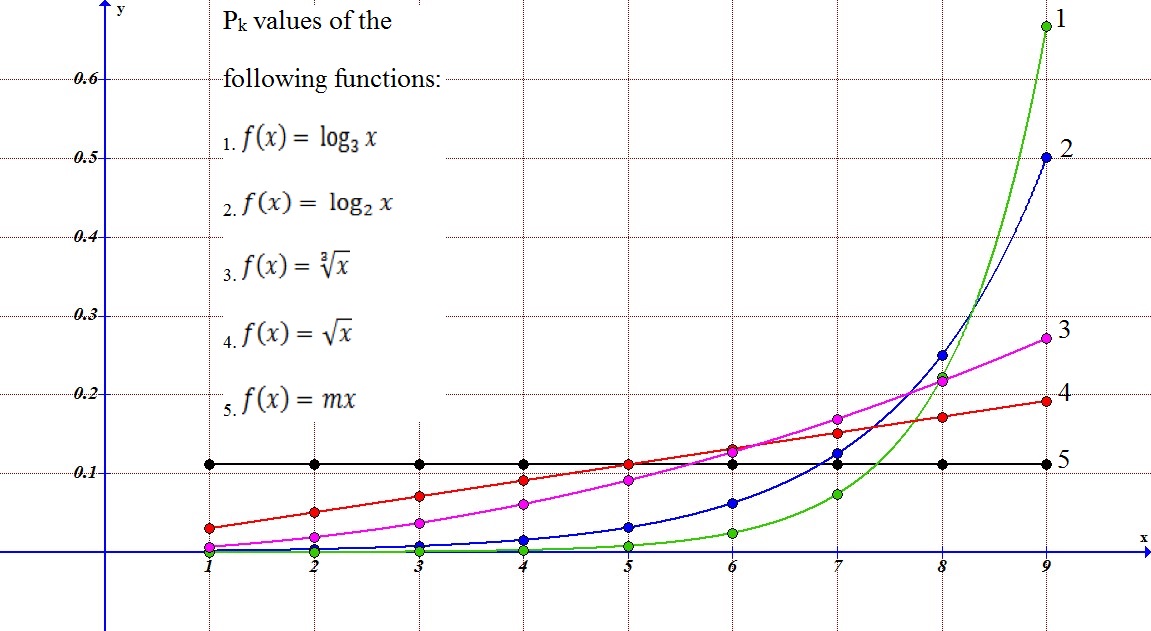} 
\caption{$P_k$ values of a linear function and functions with negative second derivatives}
\label{Combined_Probability_2}
\end{figure}

Let us look at the $P_k$ graphs of functions with increasing $P_k$ numbers. Those functions have negative second derivatives. The $P_k$ graph of a linear function was added to the picture as well for a comparison. The graphs are shown in Figure~\vref{Combined_Probability_2}. 

The $P_k$ values of both logarithmic functions are shown for the range $[1, 10)$. It was proven before that these $P_k$ numbers depend on the logarithmic function's range. It  was also shown earlier that $P_9$ values of the function are almost equal to $1$ and the rest of the $P_k$ values are almost equal to $0$ on all other ranges. Thus, we are using the range mentioned above.

According to the graph, the $P_1$ value of a root function with index $a$ is greater than the $P_1$ value of a logarithmic function with base $a$. In addition, $P_1$ of a root function with higher index will be smaller than the one with lower index. Moreover, $P_1$ of a logarithmic function with higher base will be smaller than the one with lower base.

Let us look at limits of the ratios of the functions.

$$\lim_{x \to \infty} \frac{\sqrt[a]{x+1}}{\sqrt[a]x} = 
\lim_{x \to \infty} \sqrt[a] {\frac{x+1}{x}} =
\lim_{x \to \infty} \sqrt[a] {1 + \frac{1}{x}} =
\lim_{x \to \infty} \sqrt[a] {1 + \frac{1}{x}} =
\sqrt[a] {1 + 0} = 1;$$

$$\lim_{x \to \infty} \frac{\log_a(x + 1)}{\log_a⁡(x)} = 
\lim_{x \to \infty} \frac{\frac{\ln(x + 1)}{\ln (a)}}{\frac{\ln (x)}{\ln (a)}} =
\lim_{x \to \infty} \frac{\ln (x + 1)}{\ln (x)} = 
\lim_{x \to \infty} \frac{(\ln (x + 1))'}{(\ln (x))'}$$
$$= \lim_{x \to \infty} \frac{\frac{1}{x + 1}}{\frac{1}{x}} =
\lim_{x \to \infty} \frac{x}{x + 1} =
\lim_{x \to \infty} \frac{\frac{x}{x}}{\frac{x}{x} + \frac{1}{x}} =
\lim_{x \to \infty} \frac{1}{1 + \frac{1}{x}} =
\frac{1}{1 + 0} = 1;$$

Thus, the ratios of both functions with increasing $P_k$ numbers approach to $1$. We can just assume that the reason for it is the fact that $P_1$ of those functions is always smaller than $P_1$ of the exponential one.

We had a chance to discuss how $P_k$ numbers of the continuous functions covered in this chapter relate to each other. Summarizing all the mentioned above relationships, we can conclude that if we have a real-life data set, there is a possibility that we can find a continuous function that is a good fit for the data and has similar $fdp$ values.

Nevertheless, if we have a relatively big set of random numbers, their $fdp$ values will be very close to each other. Moreover, the larger the set is, the closer to each other the $fdp$ values will be. Clearly, if the numbers are collected with no apparent rule, in other words randomly, there is no rule which will make their $fdp$ values much different from each other.

\newpage
\section{First digits probabilities of basic discrete functions}

\subsection{$P_k$ values of the Fibonacci function}

No formulas that were used by us to find the $fdp$ values of a continuous function can be used for a discrete function due to the way discrete functions are defined. However, we still can use a technique that we applied for every example done previously. Thus, we will build a table for $x$- and $y$-values and count the number of each $y$-value's first digit for every discrete function that is mentioned in this chapter.

Let us recall that the Fibonacci sequence is achieved by adding the last two consecutive numbers. The first two terms of the sequence are both $1$. In other words, 

\begin{equation}
F_n = F_{n-1} + F_{n-2}; \: n \geq 2; \: F_0 = 1; \: F_1 = 1\label{F_n_First}
\end{equation}

\begin{longtable}[c]{lllllllll} 
\caption{Fibonacci numbers and their first digits}
\label{Fibonacci}\\
\toprule
$n$- & $F_n$ & First & $n$- & $F_n$ & First & $n$- & $F_n$ & First \\
value & value & digit & value & value & digit & value & value & digit \\ 
\midrule
\endfirsthead 

\toprule
$n$- & $F_n$ & First & $n$- & $F_n$ & First & $n$- & $F_n$ & First \\
value & value & digit & value & value & digit & value & value & digit \\ 
\midrule
\endhead 
1	&	1	&	1	&	35	&	9227465	&	9	&	467	&	1.77E+97	&	1	\\
2	&	1	&	1	&	36	&	14930352	&	1	&	468	&	2.86E+97	&	2	\\
3	&	2	&	2	&	37	&	24157817	&	2	&	469	&	4.63E+97	&	4	\\
4	&	3	&	3	&	38	&	39088169	&	3	&	470	&	7.49E+97	&	7	\\
5	&	5	&	5	&	39	&	63245986	&	6	&	471	&	1.21E+98	&	1	\\
6	&	8	&	8	&	40	&	102334155	&	1	&	472	&	1.96E+98	&	1	\\
7	&	13	&	1	&	41	&	165580141	&	1	&	473	&	3.17E+98	&	3	\\
8	&	21	&	2	&	42	&	267914296	&	2	&	474	&	5.14E+98	&	5	\\
9	&	34	&	3	&	43	&	433494437	&	4	&	475	&	8.31E+98	&	8	\\
10	&	55	&	5	&	44	&	701408733	&	7	&	476	&	1.34E+99	&	1	\\
11	&	89	&	8	&	45	&	1134903170	&	1	&	477	&	2.18E+99	&	2	\\
12	&	144	&	1	&	46	&	1836311903	&	1	&	478	&	3.52E+99	&	3	\\
13	&	233	&	2	&	47	&	2971215073	&	2	&	479	&	5.70E+99	&	5	\\
14	&	377	&	3	&	48	&	4807526976	&	4	&	480	&	9.22E+99	&	9	\\
15	&	610	&	6	&	49	&	7778742049	&	7	&	481	&	1.49E+100	&	1	\\
16	&	987	&	9	&	50	&	12586269025	&	1	&	482	&	2.41E+100	&	2	\\
17	&	1597	&	1	&	...	&	...	&	...	&	483	&	3.90E+100	&	3	\\
18	&	2584	&	2	&	...	&	...	&	...	&	484	&	6.32E+100	&	6	\\
19	&	4181	&	4	&	451	&	8.02E+93	&	8	&	485	&	1.02E+101	&	1	\\
20	&	6765	&	6	&	452	&	1.30E+94	&	1	&	486	&	1.65E+101	&	1	\\
21	&	10946	&	1	&	453	&	2.10E+94	&	2	&	487	&	2.68E+101	&	2	\\
22	&	17711	&	1	&	454	&	3.40E+94	&	3	&	488	&	4.33E+101	&	4	\\
23	&	28657	&	2	&	455	&	5.49E+94	&	5	&	489	&	7.01E+101	&	7	\\
24	&	46368	&	4	&	456	&	8.89E+94	&	8	&	490	&	1.13E+102	&	1	\\
25	&	75025	&	7	&	457	&	1.44E+95	&	1	&	491	&	1.83E+102	&	1	\\
26	&	121393	&	1	&	458	&	2.33E+95	&	2	&	492	&	2.97E+102	&	2	\\
27	&	196418	&	1	&	459	&	3.77E+95	&	3	&	493	&	4.80E+102	&	4	\\
28	&	317811	&	3	&	460	&	6.09E+95	&	6	&	494	&	7.77E+102	&	7	\\
29	&	514229	&	5	&	461	&	9.86E+95	&	9	&	495	&	1.26E+103	&	1	\\
30	&	832040	&	8	&	462	&	1.60E+96	&	1	&	496	&	2.03E+103	&	2	\\
31	&	1346269	&	1	&	463	&	2.58E+96	&	2	&	497	&	3.29E+103	&	3	\\
32	&	2178309	&	2	&	464	&	4.18E+96	&	4	&	498	&	5.33E+103	&	5	\\
33	&	3524578	&	3	&	465	&	6.76E+96	&	6	&	499	&	8.62E+103	&	8	\\
34	&	5702887	&	5	&	466	&	1.09E+97	&	1	&	500	&	1.39E+104	&	1	\\
\bottomrule \\
\end{longtable}

There is one more formula for defining each term of the sequence. Since we will use it later, let us state it. 

\begin{equation}
F_n = \frac{\Phi^n - (-\Phi)^{-n}}{\sqrt 5}; \: \Phi = \frac{1 + \sqrt 5}{2} \approx 1.618
\label{F_n_Second}
\end{equation}

A list of the first $500$ Fibonacci numbers was created by using Formula~\vref{F_n_First}. 
The beginning and ending portions of the list are shown in Table~\vref{Fibonacci}. 

\begin{figure}[h]
\centering 
\includegraphics[width=1\columnwidth]{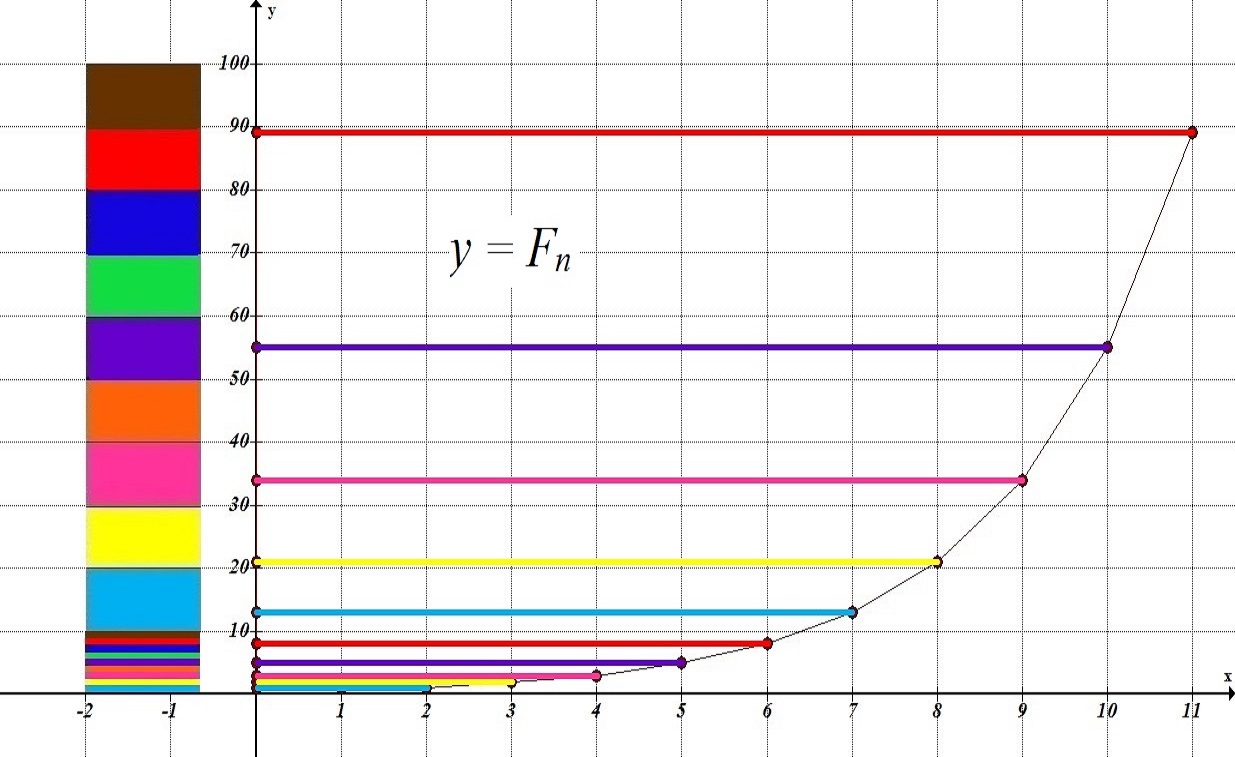} 
\caption{Fibonacci numbers in the domain $[1, 11]$}
\label{Fibonacci_Small}
\end{figure}

\begin{figure}[h]
\centering 
\includegraphics[width=1\columnwidth]{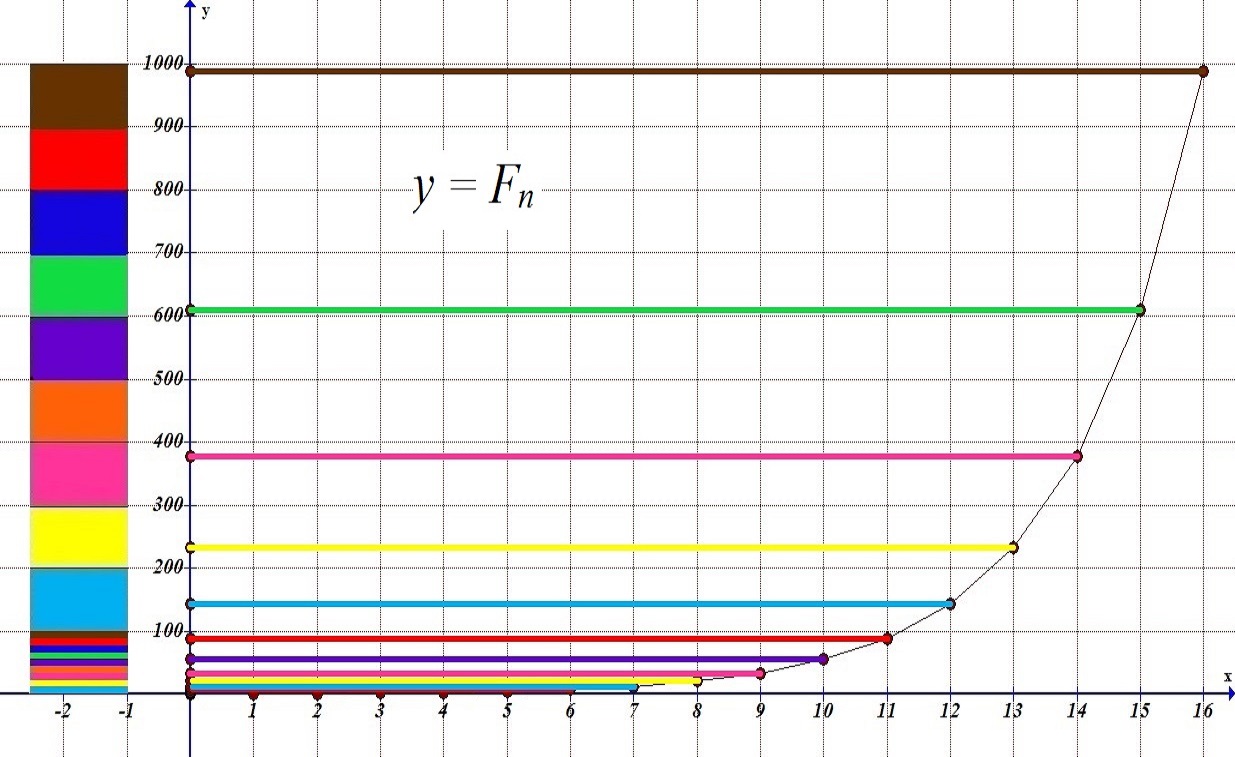} 
\caption{Fibonacci numbers in the domain $[1, 16]$}
\label{Fibonacci_Medium}
\end{figure}

\begin{figure}[h]
\centering 
\includegraphics[width=1\columnwidth]{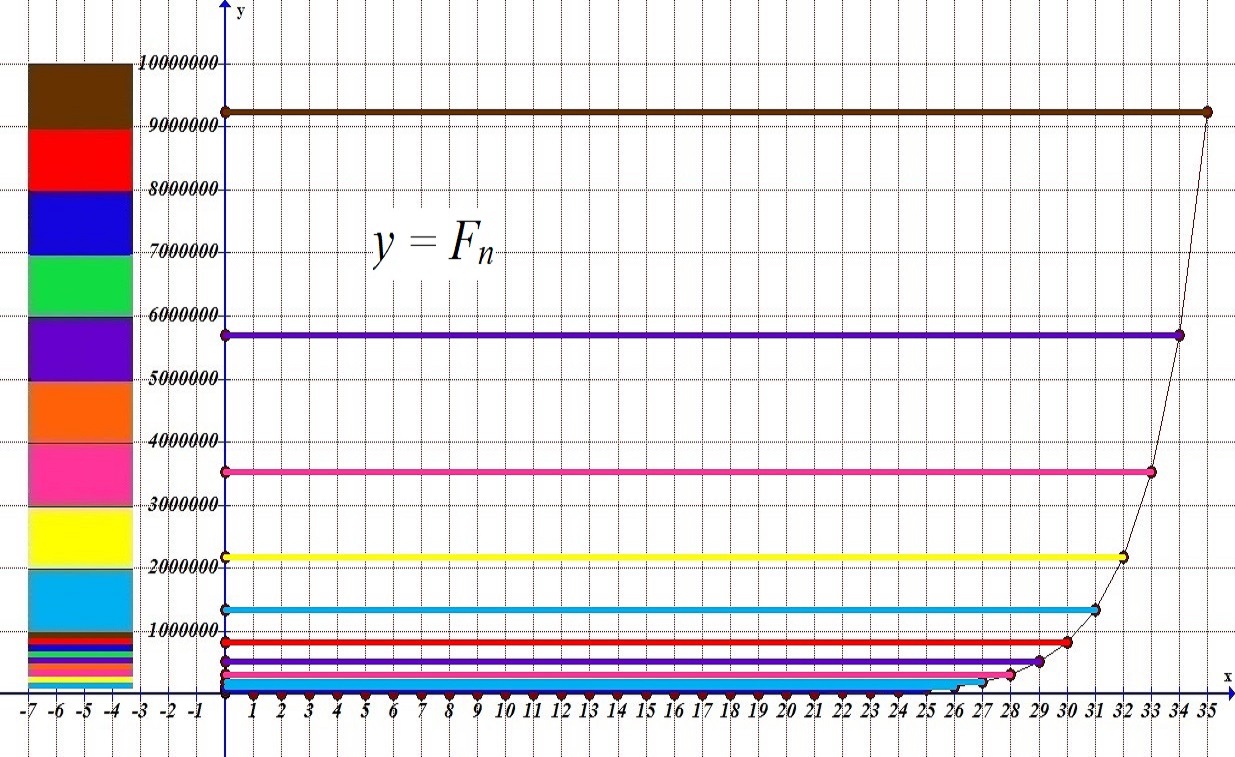} 
\caption{Fibonacci numbers in the domain $[1, 35]$}
\label{Fibonacci_Big}
\end{figure}

Let us look at three graphs that show Fibonacci numbers. The graph in Figure~\vref{Fibonacci_Small} shows the numbers in the domain $[1, 11]$, the graph in Figure~\vref{Fibonacci_Medium} shows the numbers in the domain $[1, 16]$, and the graph in Figure~\vref{Fibonacci_Big} shows the same function in the domain $[1, 35]$.

As we can see, the Fibonacci function looks like a discrete approximation to the exponential one. There is a reason for this and we will discuss it later.  So far, it looks like smaller digits have higher probabilities than bigger ones.

Let us prove it. We will calculate the probability of each first digit of the Fibonacci function using our list of the first $500$ Fibonacci numbers. Table~\vref{F_n_analysis} shows the Fibonacci $fdp$ and the exponential $fdp$.

\begin{table}[hbt]
\caption{Analysis of $F_n$ values first digits}
\label{F_n_analysis}
\centering
\begin{tabular}{cccc}
\toprule
Digit & Count & $fdp$ for $y(n) = F_n$ & $fdp$ for $y(x) = a^x$ \\  
      &       & (Count / Sum) &  \\ 
\midrule
1	&	151	&	0.30200000	&	0.30103000	\\
2	&	88	&	0.17600000	&	0.17609126	\\
3	&	63	&	0.12600000	&	0.12493874	\\
4	&	47	&	0.09400000	&	0.09691001	\\
5	&	40	&	0.08000000	&	0.07918125	\\
6	&	33	&	0.06600000	&	0.06694679	\\
7	&	29	&	0.05800000	&	0.05799195	\\
8	&	27	&	0.05400000	&	0.05115252	\\
9	&	22	&	0.04400000	&	0.04575749	\\
\bottomrule
Sum	&	500	&	1.00000000	&	1.00000000	\\
\bottomrule
\end{tabular}
\end{table}

Table~\vref{F_n_analysis} shows that smaller first digits have higher probabilities than bigger ones. 

\begin{figure}[h]
\centering 
\includegraphics[width=1\columnwidth]{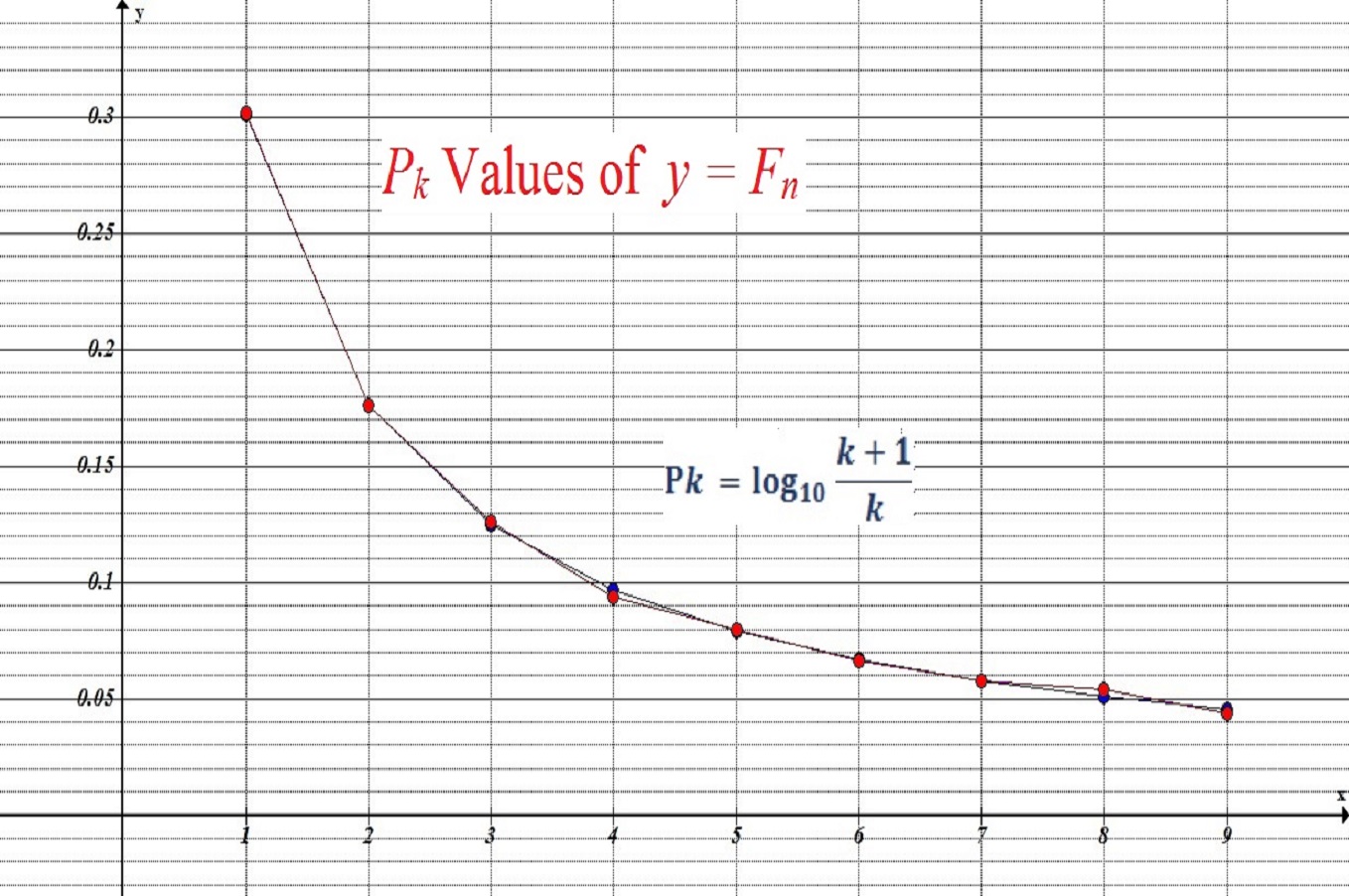} 
\caption{$P_k$ values of the Fibonacci and exponential functions}
\label{Fibonacci_Probability}
\end{figure}

Let us look at Figure~\vref{Fibonacci_Probability}, which shows the $P_k$ values of the Fibonacci function and the exponential one for a comparison. It is obvious that every $P_k$ value of the Fibonacci function is really close to the corresponding one of the exponential function. 

\FloatBarrier
$$ $$

\clearpage
\subsection{$P_k$ values of the factorial function}
Let us recall that $x! = 1 \cdot 2 \cdot 3 \cdot ... \cdot (x-1) \cdot x$, where $x$ is a natural number. As it was mentioned earlier, we are not using any formulas for discrete functions in our calculations, only tables containing numbers.

\begin{figure}[h]
\centering 
\includegraphics[width=1\columnwidth]{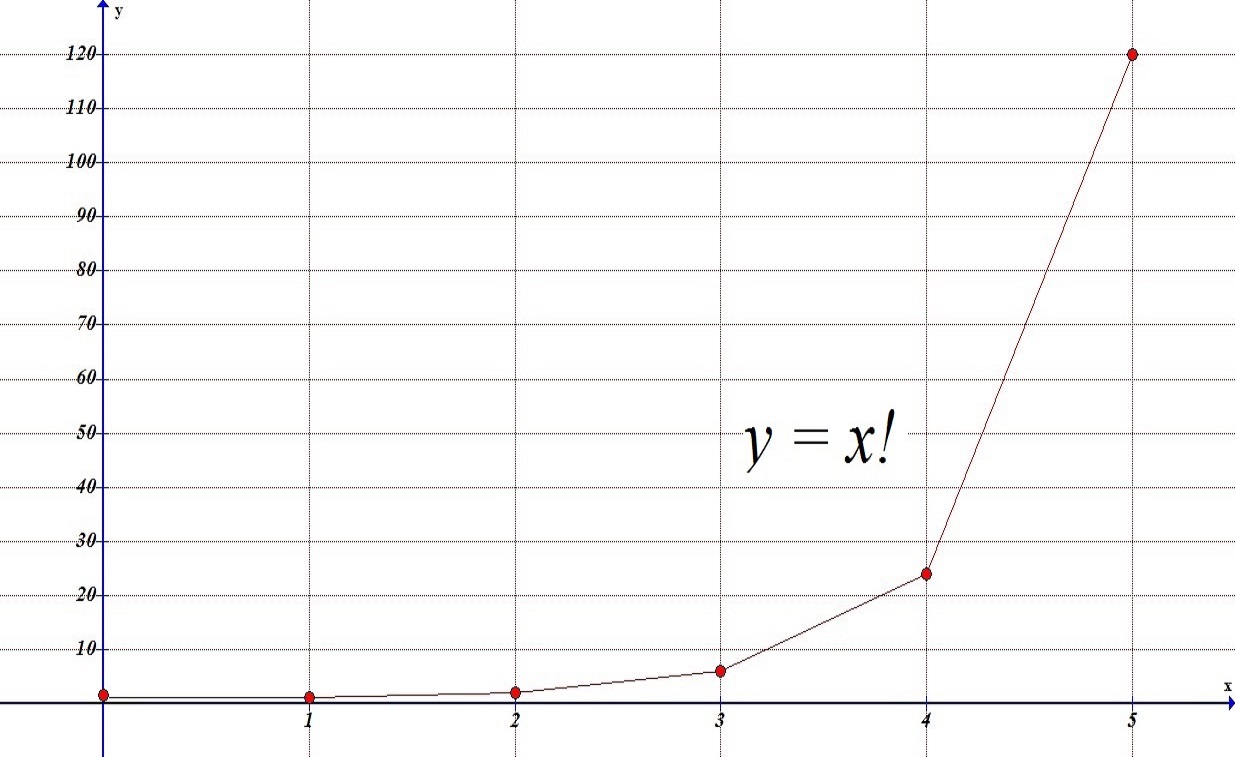} 
\caption{Factorials in the domain $[1, 5]$}
\label{Factorial_Small}
\end{figure}

\begin{figure}[h]
\centering 
\includegraphics[width=1\columnwidth]{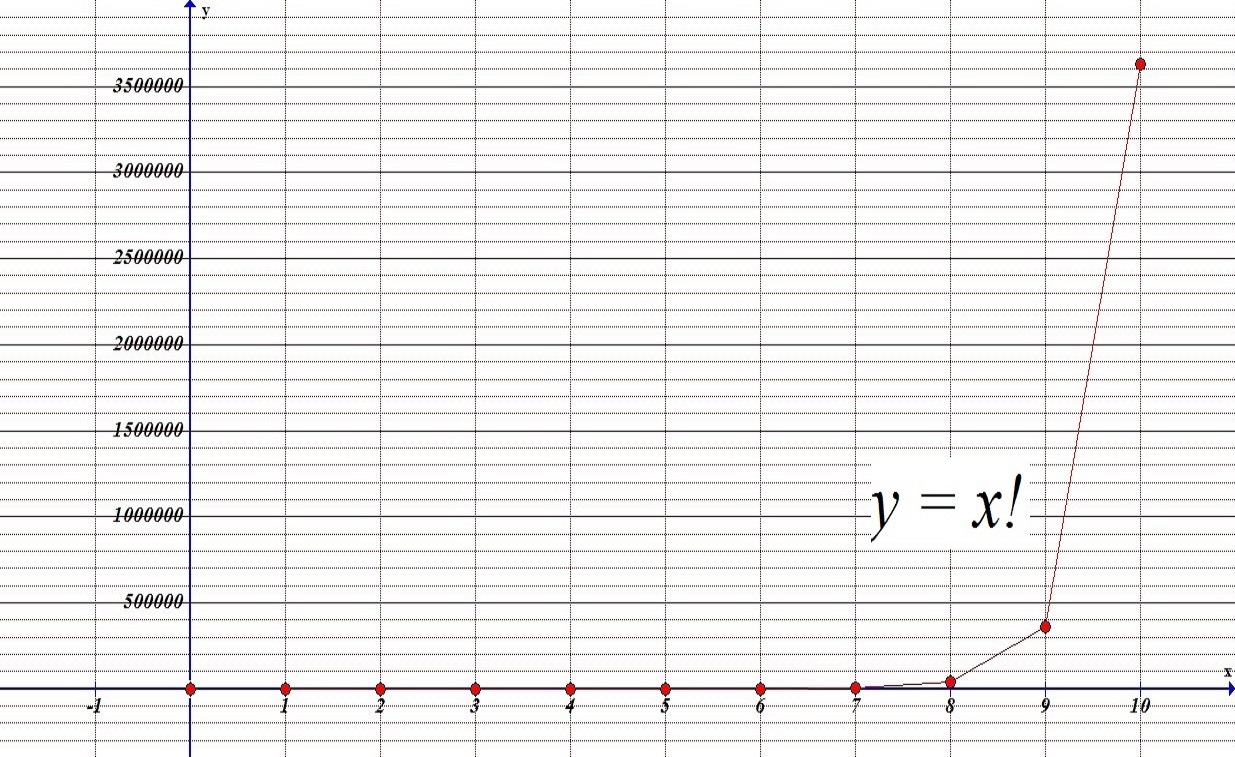} 
\caption{Factorials in the domain $[1, 10]$}
\label{Factorial_Medium}
\end{figure}

\begin{figure}[h]
\centering 
\includegraphics[width=1\columnwidth]{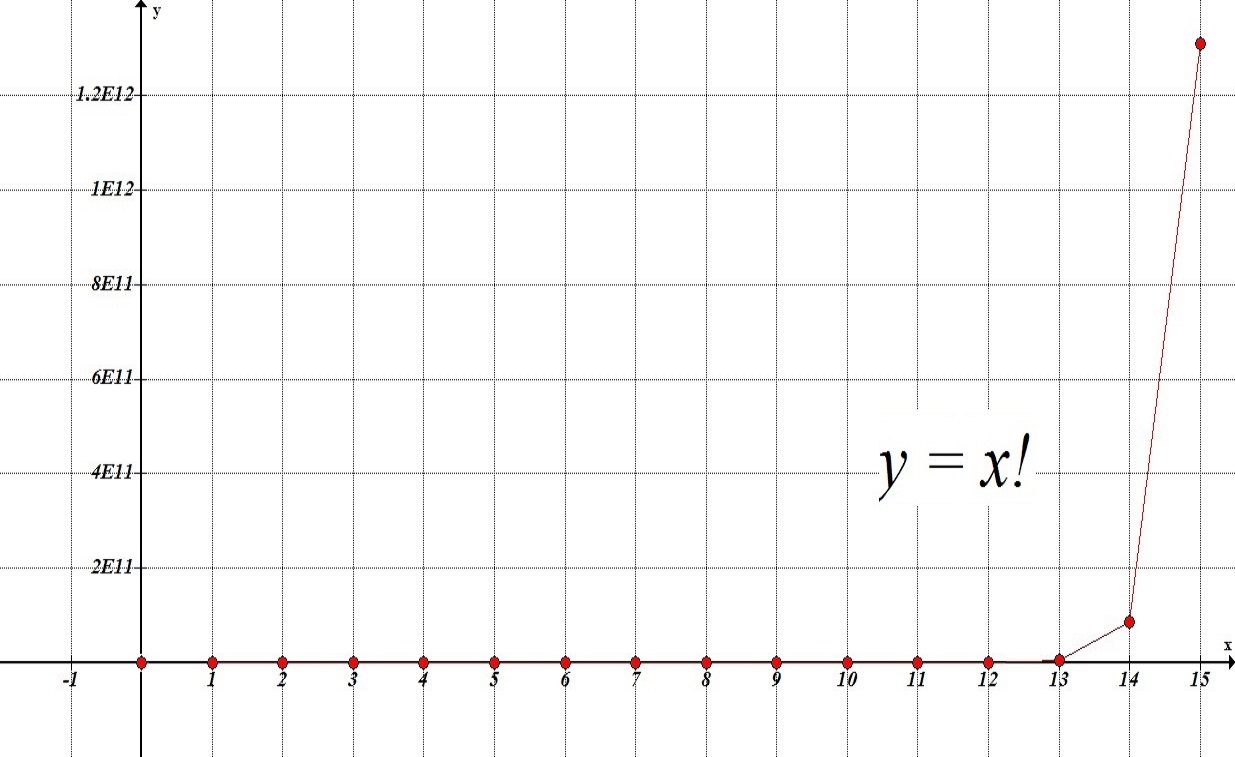} 
\caption{Factorials in the domain $[1, 15]$}
\label{Factorial_Big}
\end{figure}

\FloatBarrier
Let us look at three graphs that show factorials. The graph in Figure~\vref{Factorial_Small} shows factorials in the domain $[1, 5]$, the graph in Figure~\vref{Factorial_Medium} shows factorials in the domain $[1, 10]$, and the graph in Figure~\vref{Factorial_Big} shows the same function in the domain $[1, 15]$.

As we can see, the factorial function grows much faster than all the previous ones. In addition, even though it would be hard to use our ruler for the $fdp$ comparison for this particular function, it looks like smaller digits have higher probabilities than bigger ones and no ruler is needed for this conclusion.

\begin{longtable}[c]{lllllllll} 
\caption{Factorials and their first digits}
\label{Factorials}\\
\toprule
$x$- & $y = x!$ & First & $x$- & $y = x!$ & First & $x$- & $y = x!$ & First \\
value &  & digit & value & & digit & value & & digit \\ 
\midrule
\endfirsthead 

\toprule
$x$- & $y = x!$ & First & $x$- & $y = x!$ & First & $x$- & $y = x!$ & First \\
value &  & digit & value & & digit & value & & digit \\ 
\midrule
\endhead 
1	&	1	&	1	&	15	&	1.31E+12	&	1	&	29	&	8.84E+30	&	8	\\
2	&	2	&	2	&	16	&	2.09E+13	&	2	&	30	&	2.65E+32	&	2	\\
3	&	6	&	6	&	17	&	3.56E+14	&	3	&	31	&	8.22E+33	&	8	\\
4	&	24	&	2	&	18	&	6.40E+15	&	6	&	32	&	2.63E+35	&	2	\\
5	&	120	&	1	&	19	&	1.22E+17	&	1	&	33	&	8.68E+36	&	8	\\
6	&	720	&	7	&	20	&	2.43E+18	&	2	&	34	&	2.95E+38	&	2	\\
7	&	5040	&	5	&	21	&	5.11E+19	&	5	&	35	&	1.03E+40	&	1	\\
8	&	40320	&	4	&	22	&	1.12E+21	&	1	&	36	&	3.72E+41	&	3	\\
9	&	362880	&	3	&	23	&	2.59E+22	&	2	&	37	&	1.38E+43	&	1	\\
10	&	3628800	&	3	&	24	&	6.20E+23	&	6	&	38	&	5.23E+44	&	5	\\
11	&	39916800	&	3	&	25	&	1.55E+25	&	1	&	39	&	2.04E+46	&	2	\\
12	&	479001600	&	4	&	26	&	4.03E+26	&	4	&	40	&	8.16E+47	&	8	\\
13	&	6.22E+9 	&	6	&	27	&	1.09E+28	&	1	&	...	&	...	&	...	\\
14	&	8.71E+10	&	8	&	28	&	3.05E+29	&	3	&	...	&	...	&	...	\\
...	&	...	&	...	&	...	&	...	&	...	&	...	&	...	&	...	\\
1961	&	8.76E+206	&	8	&	1975	&	1.15E+253	&	1	&	1989	&	1.66E+299	&	1	\\
1962	&	1.72E+210	&	1	&	1976	&	2.27E+256	&	2	&	1990	&	3.31E+302	&	3	\\
1963	&	3.37E+213	&	3	&	1977	&	4.49E+259	&	4	&	1991	&	6.59E+105	&	6	\\
1964	&	6.63E+216	&	6	&	1978	&	8.88E+262	&	8	&	1992	&	1.31E+109	&	1	\\
1965	&	1.30E+220	&	1	&	1979	&	1.76E+266	&	1	&	1993	&	2.62E+112	&	2	\\
1966	&	2.56E+223	&	2	&	1980	&	3.48E+269	&	3	&	1994	&	5.22E+115	&	5	\\
1967	&	5.03E+226	&	5	&	1981	&	6.89E+272	&	6	&	1995	&	1.04E+119	&	1	\\
1968	&	9.91E+229	&	9	&	1982	&	1.37E+276	&	1	&	1996	&	2.08E+122	&	2	\\
1969	&	1.95E+233	&	1	&	1983	&	2.71E+279	&	2	&	1997	&	4.15E+125	&	4	\\
1970	&	3.84E+236	&	3	&	1984	&	5.37E+282	&	5	&	1998	&	8.29E+128	&	8	\\
1971	&	7.57E+239	&	7	&	1985	&	1.07E+286	&	1	&	1999	&	1.66E+132	&	1	\\
1972	&	1.49E+243	&	1	&	1986	&	2.12E+289	&	2	&	2000	&	3.32E+135	&	3	\\
1973	&	2.95E+246	&	2	&	1987	&	4.21E+292	&	4	&	...	&	...		&	...	\\
1974	&	5.82E+249	&	5	&	1988	&	8.37E+295	&	8	&	...	&	...		&	...	\\
\bottomrule \\
\end{longtable}

\newpage
A list of $2000$ consecutive factorials was created \cite{r_Factorial_Const}. Since the number $2000!$ contains $5736$ digits and the previous factorials are not much smaller, the exponents of the calculated factorials were truncated every time they got too big.

The beginning and the ending portions of the modified list were printed in Table~\vref{Factorials}.  In addition, the probability of each first digit of the factorial function using our list of the first $2000$ factorials was calculated. The results and the exponential $fdp$ values are shown in Table~\vref{Factorial_analysis}.

\begin{table}[hbt]
\caption{Analysis of factorials first digits}
\label{Factorial_analysis}
\centering
\begin{tabular}{cccc}
\toprule
Digit & Count & $fdp$ for $y(x) = x!$ & $fdp$ for $y(x) = a^x$ \\
      &       & (Count / Sum) &  \\ 
\midrule
1	&	591	&	0.29550	&	0.30103000	\\
2	&	335	&	0.16750	&	0.17609126	\\
3	&	250	&	0.12500	&	0.12493874	\\
4	&	204	&	0.10200	&	0.09691001	\\
5	&	161	&	0.08050	&	0.07918125	\\
6	&	156	&	0.07800	&	0.06694679	\\
7	&	107	&	0.05350	&	0.05799195	\\
8	&	102	&	0.05100	&	0.05115252	\\
9	&	94	&	0.04700	&	0.04575749	\\
\bottomrule
Sum	&	2000	&	1.00000	&	1.00000000	\\
\bottomrule
\end{tabular}
\end{table}

Table~\vref{Factorial_analysis} shows that smaller first digits have higher probabilities than bigger ones. 

\begin{figure}[h]
\centering 
\includegraphics[width=1\columnwidth]{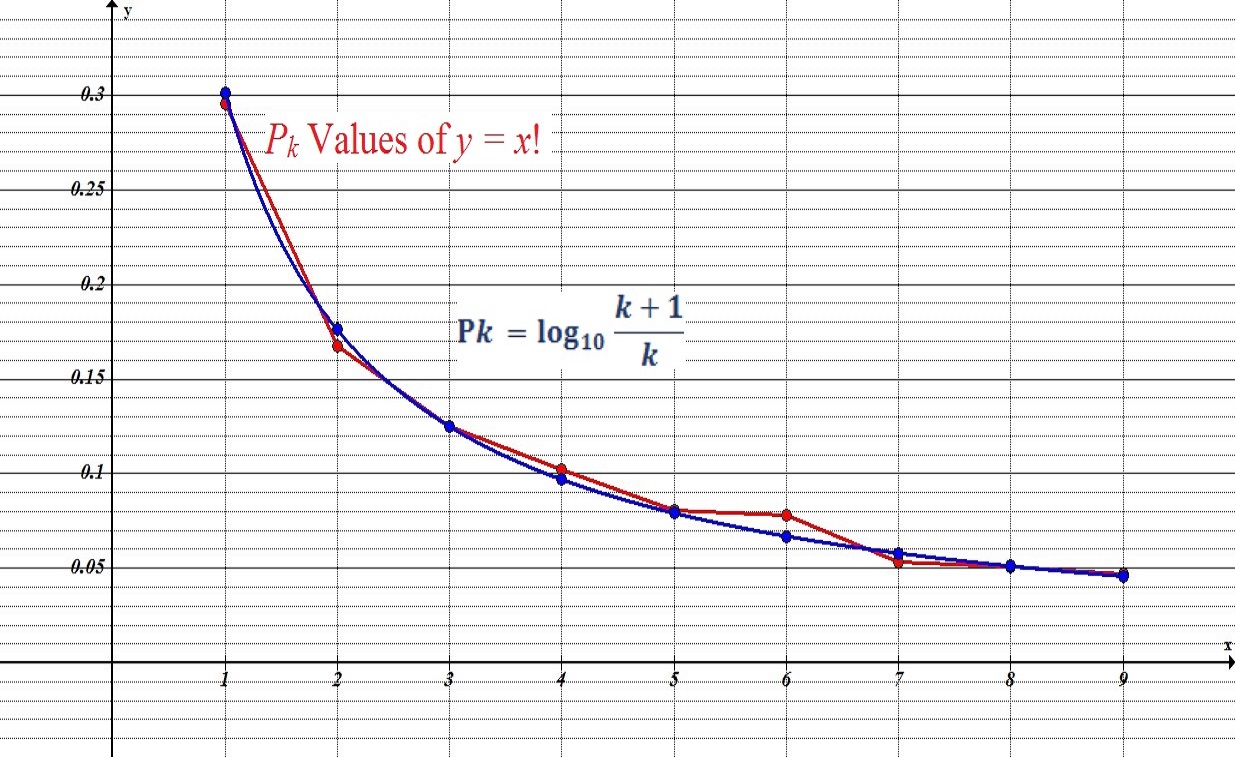} 
\caption{$P_k$ values of $y = x!$ and $y = a^x$}
\label{Factorial_Probability}
\end{figure}

Let us look at Figure~\vref{Factorial_Probability}, which shows $P_k$ values of the function $y = x!$ and the exponential function for a comparison. 

It is obvious that every $P_k$ value of the factorial function is really close to the corresponding one of the exponential function. It looks like that if our domain were bigger, the $P_k$ values would be even closer together. We will discuss this hypothesis later.

\newpage
\subsection{Summary of $P_k$ values of discrete functions}

Let us talk about $P_k$ values of two discrete functions discussed in this chapter. 

We will talk about the Fibonacci function first. Based on our calculations and graphs, we made a hypothesis that the function follows Benford's distribution. Let us prove it by taking the limit of the ratio $\frac{f(x+1)}{f(x)}$, where $f(x) = F_n$.

\begin{proof}
As it was mentioned above, $F_n = \frac{\Phi^n - (- \Phi)^{-n}}{\sqrt 5}$, where $\Phi = \frac{1 + \sqrt 5}{2} \approx 1.618$.

$$\lim_{n \to \infty} \frac{F_{n+1}}{F_n} = 
\lim_{n \to \infty} \frac {\frac{\Phi^{n+1} - (- \Phi)^{-(n+1)}}{\sqrt 5}}{\frac{\Phi^n - (- \Phi)^{-n}}{\sqrt 5}} = 
\lim_{n \to \infty} \frac {\Phi^{n+1} - (- \Phi)^{-(n+1)}}{\Phi^n - (- \Phi)^{-n}}$$

$$= \lim_{n \to \infty} \frac {\Phi^{n+1} -  (-1)^{-(n+1)}\Phi^{-(n+1)}}{\Phi^n - 
(-1)^{-n}\Phi^{-n}} =
\lim_{n \to \infty} \frac {\Phi^{n+1} -  (-1)^{-(n+1)}\cdot \frac{1}{\Phi^{n+1}}}{\Phi^n - 
(-1)^{-n}\cdot \frac{1}{\Phi^n}}$$

$$= \lim_{n \to \infty} \frac {\Phi^n(\Phi -  (-1)^{-(n+1)}\cdot \frac{1}{\Phi^{2n+1}})}{\Phi^n(1 - (-1)^{-n}\cdot \frac{1}{\Phi^2n})} =
\lim_{n \to \infty} \frac {\Phi -  (-1)^{-(n+1)}\cdot \frac{1}{\Phi^{2n+1}}}{1 - (-1)^{-n}\cdot \frac{1}{\Phi^2n}}$$

$$= \frac{\Phi - 0}{1 - 0} = \Phi;$$

Thus, the ratio of the Fibonacci function approaches to $\Phi$ and the function can be approximated to a continuous exponential function $y = \Phi^x$.
\end{proof}

First, it explains why the Fibonacci function looks like an exponential one if we graph it on a large domain. In addition, it explains why their corresponding $P_k$ numbers are equal to each other.

Now, we should talk about the function $y = x!$. Unfortunately, the limit of the ratio $\frac{f(x+1)}{f(x)}$ will not lead us to any meaningful results.

$$\lim_{x \to \infty} \frac{(x+1)!}{x!} = \lim_{x \to \infty} \frac{x! \cdot (x+1)}{x!} = 
 \lim_{x \to \infty} (x+1) = \infty;$$

It may look like $P_k$ numbers of the factorial function should be really different from those of the exponential one; however, we know that this is not true. All of our calculations show that $P_k$ numbers of the factorial function approach and get really close to $P_k$ numbers of the exponential function on a large domain. Thus, more work needs to be done to explain why the function $y = x!$ has Benford's probability.

According to Stirling's formula \cite{r_Stirling}, 

\begin{equation}
x! \approx \frac{x^x}{e^x} \cdot \sqrt {2\pi x} = \sqrt {2\pi x}\cdot (\frac{x}{e})^x
\label{Stirling_Formula}
\end{equation}

Thus, the factorial function can be approximated to a continuous function, which looks like an exponential one with a particular coefficient. In other words, 

$$x! \approx ma^x, m = \sqrt {2\pi x}, a = \frac{x}{e}$$ 

However, the coefficient and the base of the function are not constant values, which makes the function different from a regular exponential one.

Let us do a couple of limits. First, we will take a limit of the ratio $\frac{m(x+1)}{m(x)}$ and 
then we will do a limit of $\frac{a(x+1)}{a(x)}$ where $x$ approaches to infinity.

$$\lim_{x \to \infty} \frac{m(x+1)}{m(x)} = 
\lim_{x \to \infty} \frac{\sqrt {2\pi (x+1)}}{\sqrt {2\pi x}} =
\lim_{x \to \infty} \frac{\sqrt {x+1}}{\sqrt x} =
\lim_{x \to \infty} \sqrt {\frac{x+1}{x}}$$ 

$$= \lim_{x \to \infty} \sqrt {1 + \frac{1}{x}} = \sqrt {1 + 0} = 1;$$

$$\lim_{x \to \infty} \frac{a(x+1)}{a(x)} = 
\lim_{x \to \infty} \frac{\frac{x+1}{e}}{\frac{x}{e}} =
\lim_{x \to \infty} \frac{x+1}{x} =
\lim_{x \to \infty} (1 + \frac{1}{x}) = 1 + 0 = 1;$$

Thus, even though the base and the coefficient depend on $x$, their ratios approach to $1$. It means that as $x$ approaches to infinity, the base and the coefficient become more stable.

Again, the function $f(x) = x! \approx \sqrt {2\pi x}\cdot (\frac{x}{e})^x$ does not belong to a group of exponential ones; however, it has the same $P_k$ numbers as the exponential functions do. There is a chance that the two ratios mentioned above shed some light on this puzzle. 


\clearpage
\section{References}
\renewcommand{\refname}{\spacedlowsmallcaps{}} 

\bibliographystyle{unsrt}


\end{document}